\documentclass[11pt]{article}
\usepackage[utf8]{inputenc}

\usepackage{amsthm, amsfonts, amssymb, amsmath,enumitem, natbib, bbm, mathabx}
\usepackage{cases, amsthm}
\usepackage{fullpage}
\usepackage{wasysym} 
\usepackage{fancyhdr}
\usepackage{graphicx}
\usepackage{bbm}
\usepackage{caption}
\usepackage{subcaption}
\usepackage{authblk}
 \usepackage{setspace}
\usepackage[colorlinks=true, pdfstartvie w=FitV,
linkcolor=blue, citecolor=blue, urlcolor=blue]{hyperref}
\usepackage{comment}
\usepackage{mathtools}

\newtheorem{theorem}{Theorem}
\newtheorem{prop}{Proposition}

\newtheorem{corollary}{Corollary}

\newtheorem{remark}{Remark}

\usepackage[linesnumbered,ruled,vlined]{algorithm2e}
\newcommand{\indep}{\perp \!\!\! \perp}

\usepackage{mathrsfs}

\usepackage{xcolor}
\usepackage{cleveref}

\usepackage[top=1in, bottom=1in, left=1in, right=1in]{geometry}
\title{On a Probability Inequality for Order Statistics with Applications to Bootstrap, Conformal Prediction, and more}

\author[1]{Manit Paul}
\author[2]{Arun Kumar Kuchibhotla}
\affil[1]{Department of Statistics \& Data Science, University of Pennsylvania}
\affil[2]{Department of Statistics \& Data Science, Carnegie Mellon University}
\date{}

\begin{document}

\maketitle

\begin{abstract}
``Behind every limit theorem, there is an inequality'' said Kolmogorov. We say ``for every inequality, there is an approximate inequality under approximate regularity conditions.'' Suppose $X, X'$ are independent and identically distributed random variables. Then $X \le X'$ with a probability of at least $1/2$, irrespective of the underlying (common) distribution. One can ask what happens to the probability if $X, X'$ are independent but not identically distributed. It should be approximately $1/2$ if the distributions are approximately equal. Similarly, what if the random variables are dependent? It should, again, be approximately $1/2$ if the random variables are approximately independent. We explore an extension of this probability inequality involving order statistics and develop approximate versions of such an inequality under violations of independence and identical distribution assumptions. We further show that this inequality can be used as a basis to prove asymptotic validity of bootstrap/subsampling, finite-sample validity of conformal prediction, permutation tests, and asymptotic validity of rank tests without group invariance. Specifically, in the context of resampling inference, our results can be seen as a finite-sample instantiation of some results by Peter Hall and yield an alternative ``cheap bootstrap'' that applies to high-dimensional data.
\end{abstract}

\section{Introduction and Motivation}
\label{sec:introduction}
Suppose we are given $B$ independent and identically distributed (IID) real-valued random variables $W_1, \ldots, W_B$ with a common distribution function $F$. Let $W_{(1)} \le W_{(2)} \le \ldots \le W_{(B)}$ be the order statistics. If $W$ is an independent copy of $W_1$ and $F(\cdot)$ is continuous, then $\mathbb{P}(W \in (W_{(k_1)}, W_{(k_2)}])$ is independent of $F$ and can be evaluated explicitly as a function of $B, k_1,$ and $k_2$. This fact follows from the well-known distribution-free property of order statistics~\citep[Chapter 1]{reiss1989approximate}. Relaxing the assumption of continuity of $F(\cdot)$, it can be proved that there exists a function $(B, k_1, k_2)\mapsto h(B, k_1, k_2)\in[0, 1]$ such that 
\begin{equation}\label{eq:order-statistics-domination}
    \mathbb{P}_{F}(W \in (W_{(k_1)},\, W_{(k_2)}]) ~\ge~ h(B, k_1, k_2), \quad\mbox{for all}\quad F, B, 1\le k_1 \le k_2 \le B.  
\end{equation}
It is important to note that $h(B, k_1, k_2)$ depends only on $B, k_1, k_2$ and not on $F$.
Inequality~\eqref{eq:order-statistics-domination} becomes an equality whenever $F$ is continuous, implying a distribution-free property. By choosing $k_1, k_2$ such that $h(B, k_1, k_2) = 1-\alpha$, one can construct a prediction interval for a ``test point'' $W$. 
This forms the basis of the construction of finite-sample valid prediction sets and have garnered a lot of interest in the recent times through the conformal prediction literature~\citep{lei2014distribution,lei2018distribution,angelopoulos2023conformal,angelopoulos2024theoretical}. 
Beyond this simple setting of IID random variables, several works within conformal prediction has considered the question of distribution-shift (non-identical distributions) and non-independent random variables. The considerations here can be divided into two parts: what if $W_1, \ldots, W_B$ are IID from $F$ and $W$ is from a (potentially) different distribution $G$? and What if $(W_1, \ldots, W_B, W)$ is an arbitrary random vector (with an arbitrary dependence structure)? In the first case, it is easy to see that inequality~\eqref{eq:order-statistics-domination} can be corrected as
\begin{equation}\label{eq:corrected-order-statistics-domination}
\mathbb{P}_{G, F}(W \in (W_{(k_1)}, W_{(k_2)}]) \ge h(B, k_1, k_2) - \Delta(F, G),
\end{equation}
for some metric on the set of distributions such that when $F = G$, $\Delta(F, G) = 0$ and we recover~\eqref{eq:order-statistics-domination}. When $(W_1, \ldots, W_B, W)$ is an arbitrary random vector,~\cite{barber2023conformal} provide an extension of~\eqref{eq:order-statistics-domination} involving the total variation distance. 

Going back to inequality~\eqref{eq:order-statistics-domination}, it should be mentioned that for a given value of $B \ge 1$ and $\alpha\in(0, 1)$, there may not exist $k_1, k_2$ such that $h(B, k_1, k_2) \ge 1-\alpha$. It can be proved that such $k_1, k_2$ exist only if $B \ge 1/\alpha - 1$. In the context of prediction set construction, $B$ represents the number of observations available which is typically assumed large enough such that this constraint on $B$ is satisfied. Although theoretically interesting, the practical utility of the corrected results under approximate similarity of distributions or independence is debatable (in the context of conformal prediction). For example, for~\eqref{eq:corrected-order-statistics-domination} to be useful, one has to have $\Delta(F, G)$ to ``converge'' to zero which as $B\to\infty$ means that the distributions $F$ and $G$ are changing with $B$. Interestingly enough, in the context of bootstrap or resampling, such a requirement of $\Delta(F, G) \to 0$ makes sense without letting $B\to\infty.$

Let us briefly describe the bootstrap methodology to illustrate how inequality~\eqref{eq:order-statistics-domination} plays a crucial role. Here we shall restrict to bootstrap confidence interval for the mean. Suppose $X_1, \ldots, X_m$ are IID real-valued random variables from some distribution $P$ with mean $\mu$. Let $W = m^{1/2}(\overline{X}_m - \mu)$, where $\overline{X}_m = m^{-1}\sum_{i=1}^m X_i$. Fix any $B \ge 1$ and for $1\le b\le B$, let $X_1^{*b}, \ldots, X_m^{*b}$ denote a simple random sample with replacement from $X_1, \ldots, X_m$. Then set
\[
W_b = m^{1/2}(\overline{X}_m^{*b} - \overline{X}_m)\quad\mbox{for}\quad 1\le b\le B.
\]
The usual bootstrap confidence interval is given by
\[
\widehat{\mathrm{CI}}_{m,\alpha}^{\mathrm{boot}} := \left\{\theta\in\mathbb{R}:\, m^{1/2}(\overline{X}_m - \theta) \in \left(W_{(B\alpha/2)},\, W_{(B(1-\alpha/2))}\right]\right\}.
\]
(Assume, for simplicity, that $B\alpha/2$ is an integer.) The coverage probability of this bootstrap confidence set is given by $\mathbb{P}(W \in (W_{(B\alpha/2)}, W_{(B(1-\alpha/2))}])$, which is precisely the left hand side of~\eqref{eq:order-statistics-domination}, albeit without the independence and identical distributions for $W_1, \ldots, W_B, W$. Typically, the validity of bootstrap is argued in two steps: (1) the conditional distribution of $W_1$ given $X_1, \ldots, X_m$ asymptotically (as $m\to\infty$) matches the distribution of $W$; and (2) as $B\to\infty$, the empirical distribution of $W_1, \ldots, W_B$ converges to the conditional distribution which implies that $W_{(B\alpha/2)}, W_{(B(1-\alpha/2))}$ converge to the corresponding quantiles of the conditional distribution of $W_1$ given $X_1, \ldots, X_m$. While this argument is correct, it requires both $B, m\to\infty$. Relating to inequality~\eqref{eq:order-statistics-domination}, in the context of bootstrap, we have that $W_1, \ldots, W_B$ are IID conditional on $X_1, \ldots, X_m$ and the distribution of $W$ is close to the conditional distribution of $W_b$'s conditional on $X_1, \ldots, X_m$ (as $m\to\infty$). Now considering~\eqref{eq:corrected-order-statistics-domination}, one must have correct coverage as $m\to\infty$ but not necessarily $B \to \infty$. This was precisely the message of~\cite{hall1986number} where a slight modification of $\widehat{\mathrm{CI}}_{m,\alpha}^{\mathrm{boot}}$ is used to ensure $h(B, k_1, k_2) = 1-\alpha$. His results imply that as long as $B \ge 1/\alpha - 1$, the one-sided bootstrap confidence intervals have asymptotically correct coverage. (Note that we treat bootstrap confidence interval as a randomized confidence interval and consider marginal coverage. Some authors discuss conditional coverage of bootstrap confidence interval, i.e., $\mathbb{P}(\mu\in\widehat{\mathrm{CI}}_{m,\alpha}^{\mathrm{boot}}|X_1, \ldots, X_m)$~\citep{andrews2000three}. We do not consider this viewpoint in this paper.)

In addition to bootstrap, several other popular methods in statistics rely on~\eqref{eq:order-statistics-domination} and its extensions. All resampling methods, including subsampling~\citep{romano1999subsampling,politis1994large} and permutation testing~\citep{fisher1935design,ramdas2023permutation}, rely on~\eqref{eq:corrected-order-statistics-domination}. In addition, the extensions of~\eqref{eq:order-statistics-domination} alluded to in the discussion above also imply that certain rank tests and permutation tests remain valid under approximate notions of the underlying invariance conditions, thus potentially recovering some of the results of~\cite{romano1989bootstrap,romano1990behavior}.
Apart from connecting all these areas, our results reiterate that the number of Monte Carlo replicates needed in all these methods need not converge with the sample size, as long as they are sufficiently large when compared to a function of miscoverage level $\alpha$. This should not be interpreted as recommending the use of the smallest possible number of Monte Carlo replications. For example, in the context of bootstrap, the traditional method of proof requires $B$ to diverge in order to ensure that the bootstrap quantiles converge to the true quantiles. This ensures that the resulting confidence interval has a ``low variance'' and matches the classical Wald interval in width asymptotically if $B, m \to\infty$. Such a conclusion cannot be drawn and in fact, does not hold if $B$ is fixed as $m\to\infty$. The message of our paper is that validity is not lost if one chooses to use a fixed number of Monte Carlo replications.

\subsection{Contributions}

Our main contributions are as follows.

\begin{enumerate}[leftmargin=*, itemsep=0.35em]
    \item We establish an extension of~\eqref{eq:order-statistics-domination} for the coverage of a real-valued random variable $\psi(Z)$, a function of any random variable $Z$, by the interval $(W_{(k_1)},\, W_{(k_2)}]$ formed by the order statistics of random variable $W_1, \ldots, W_B$ in three regimes: $W_1, \ldots, W_B$ are (i) independent and identically distributed conditional on $Z$, (ii) independent but non-identically distributed conditional on $Z$, or (iii) arbitrarily dependent. The resulting bounds are non-asymptotic. Moreover, the coverage slack is made explicit through some metrics that measure the ``distance'' between $(W_1, \ldots, W_B, \psi(Z))$ from the idealized case of IID real-valued random variables. The coverage slack in cases (i) and (ii) can converge to zero even if $B$ does not diverge to $\infty$.

    \item We apply the master theorems to construct valid confidence intervals and tests, whether finite-sample or asymptotic, in a range of settings including the nonparametric bootstrap, subsampling, multiplier bootstrap inference for stochastic gradient descent, permutation testing, conformal prediction, and randomization tests without exact group invariance. In each of these settings, a minor modification of the conventional order-statistic cutoffs removes the standard condition that the number of Monte Carlo draws must diverge to ensure validity.

    \item We provide further refinements in addition to the validity results for confidence intervals and tests. Specifically, we propose a method based on external randomization for constructing modified confidence intervals that attain exact asymptotic coverage in nonparametric bootstrap, subsampling, and related settings using only a finite number of Monte Carlo resamples. We also obtain sharper coverage analyses for certain conformal prediction problems in non-exchangeable settings.
\end{enumerate}

The simulation studies in Section~\ref{sec:applications} show that the proposed modifications achieve the nominal coverage or level with small fixed Monte Carlo budgets, at the cost of only a modest increase in interval width relative to the conventional procedures.
\paragraph{\textbf{Organization}.}
Section~\ref{sec:master_thms} develops the main theorems that characterize the coverage of the target statistic by confidence intervals constructed from test statistics computed using resampled data. Section~\ref{sec:applications} applies these results to nonparametric bootstrap, subsampling, inference based on stochastic gradient descent, permutation testing, conformal prediction, and randomization tests without exact group invariance, and also reports evidence from simulation studies. Section~\ref{sec:extension} concludes with a discussion. Proofs of all main results, propositions, and corollaries are collected in the appendix.



\section{Master theorems}
\label{sec:master_thms}
Consider $B$ univariate data points $W_1,\ldots,W_B$ and a random variable (or random vector) $Z$ defined on a common probability space $(\Omega, \mathscr F, \mathbb P)$ where $W_i : \Omega \mapsto \mathcal W$ for $i \in [B]$ (where $\mathcal W \subset \mathbb R$) and $Z : \Omega \mapsto \mathcal Z$. Let $\mathcal B(\mathcal W^B)$ be the Borel $\sigma-$algebra on the product space $\mathcal W^B$ and let $P_{W^{(B)}|Z}(z, \cdot): \mathcal B(\mathcal W^B) \mapsto [0,1]$ be a regular conditional distribution of $W^{(B)} = (W_1,\ldots,W_B)$ given $Z$ i.e.\ $(i)$ for each $z \in \mathcal Z$, $A \mapsto P_{W^{(B)}|Z}(z, A)$ is a probability measure on $(\mathcal W^B, \mathcal B(\mathcal W^B))$ and $(ii)$ for each measurable $A \subset \mathcal W^B$, $z \mapsto P_{W^{(B)}|Z}(z, A)$ is measurable. We present some probability inequalities in this section regarding the coverage of $\psi(Z)$ (where $\psi: \mathcal{Z} \mapsto \mathbb R$) by intervals based on the order statistics $(W_{(1)}, \ldots, W_{(B)})$ under suitable assumptions on $P_{W^{(B)}|Z}(z, \cdot)$. 

For any two probability measures $\nu_1(\cdot), \nu_2(\cdot)$ and for any collection of measurable sets $\mathcal{I}$ we define $d_{\mathcal{I}}(\nu_1,\nu_2) = \sup_{A \in \mathcal{I}}|\nu_1(A) - \nu_2(A)|$. Let $\mathcal{I}_0 = \{(\infty, b]| b \in \mathbb R\}$, $\mathcal{I}_1 = \{(a, b]| a,b\in \mathbb R\}$ and $\mathcal{I}_2 = \{I| I \mbox{ is measurable}\}$ be three collections of measurable sets. We define $d_{\mathrm{KS}}(\nu_1,\,\nu_2) = d_{\mathcal I_0}(\nu_1,\,\nu_2)$, $\widetilde d_{\mathrm{KS}}(\nu_1,\,\nu_2) = d_{\mathcal I_1}(\nu_1,\,\nu_2)$ and $d_{\mathrm{TV}}(\nu_1,\,\nu_2) = d_{\mathcal I_2}(\nu_1,\,\nu_2)$. We note that $d_{\mathrm{KS}}(\nu_1,\,\nu_2)$ and $d_{\mathrm{TV}}(\nu_1,\,\nu_2)$ are, respectively, the standard Kolmogorov-Smirnov and total variation distances between $\nu_1$ and $\nu_2$. The modified Kolmogorov-Smirnov distance satisfies 
\[
d_{\mathrm{KS}}(\nu_1,\,\nu_2) ~\leq~ \widetilde d_{\mathrm{KS}}(\nu_1,\,\nu_2) ~\leq~ \min \{2d_{\mathrm{KS}}(\nu_1,\,\nu_2), d_{\mathrm{TV}}(\nu_1,\,\nu_2)\}.
\]
If $Y_1\sim\nu_1$ and $Y_2\sim \nu_2$, then we also use the notation $d_{\mathrm{KS}}(Y_1, Y_2)$ to denote $d_{\mathrm{KS}}(\nu_1, \nu_2)$. The same is also done for the modified KS and the TV distances. 

The next theorem describes the coverage of $\psi(Z)$ under the assumption that $W_1,\ldots,W_B$ are distributed independently and identically upon conditioning on $Z$ i.e.\ $P_{W^{(B)}|Z}(z, \cdot) = (P_{W_1|Z}(z, \cdot))^{\otimes B}$ for all $z \in \mathcal Z$.
\begin{theorem}
\label{thm:main_result}
Suppose $W_1,\ldots,W_B$ are $B$ random variables which conditional on another random vector $Z$ are distributed identically and independently. Define $F_0(z) = \mathbb{P}(W_1 \leq \psi(Z)| Z = z)$, $\widetilde F_0(z) = \mathbb{P}(W_1 < \psi(Z)| Z = z)$. Let $U(0,1)$ be a standard uniform random variable. Set
\[
\Delta ~:=~ \widetilde{d}_{\mathrm{KS}}(F_0(Z),\, U(0, 1))\quad\mbox{and}\quad \widetilde{\Delta} ~:=~ \widetilde{d}_{\mathrm{KS}}(\widetilde{F}_0(Z),\, U(0, 1)).
\]
Fix any two integers $a, b$ satisfying $0 \leq a < B-b \leq B$. Then 
\begin{equation*}
\begin{split}
  - \Delta  ~\leq~ &\mathbb{P}\left(\psi(Z) \in [W_{(a)}, W_{(B-b)}] \right) -\left(  1- \frac{a+b +1}{B+1} \right) ~\leq~  \frac{1}{B+1} + \Delta,\\
    - \Delta ~\leq~ &\mathbb{P}\left(\psi(Z) \in [W_{(a)}, W_{(B-b)})  \right) -\left(  1- \frac{a+b +1}{B+1} \right) ~\leq~   \Delta, \\
     - \widetilde{\Delta} ~\leq~ & \mathbb{P}\left(\psi(Z) \in (W_{(a)}, W_{(B-b)}] \right) -\left(  1- \frac{a+b +1}{B+1} \right) ~\leq~  \widetilde{\Delta},
\end{split}
\end{equation*}
We set the order statistic $W_{(0)} =\min\{w|w \in \mathcal W\}$. 
\end{theorem}
\Cref{thm:main_result} characterizes the coverage of $\psi(Z)$ by a confidence interval based on the order statistics $(W_{(1)}, \ldots, W_{(B)})$ when $W_1,\ldots W_B$ are IID conditioned on $Z$. We note that \cite{hall1986number} investigated the coverage properties of confidence interval $(W_{(a)}, W_{(B-b)}]$ and established a comparable result in the special case where $a = 0$. The detailed proof of \Cref{thm:main_result} is provided in \Cref{app:proof_main_res}. To contain $\psi(Z)$ with a $(1-\alpha)$ probability, we can set $a = \lfloor (B+1) (\alpha/2) \rfloor$ and $b = \lfloor (B+1)\alpha) \rfloor - \lfloor (B+1) (\alpha/2) \rfloor - 1$ and use the second statement of \Cref{thm:main_result} to get, 
\begin{equation*}
    \begin{split}
  \mathbb{P}\left(\psi(Z) \in [W_{(\lfloor (B+1) (\alpha/2) \rfloor)}, W_{(\lceil (B+1)(1 -\alpha) \rceil + \lfloor (B+1) (\alpha/2) \rfloor)}) \right)  ~\geq &~ 1 - \frac{\lfloor (B+1) \alpha \rfloor}{B+1} - \Delta   \\
  ~\geq &~ 1 - \alpha - \Delta  .      
    \end{split}
\end{equation*}
The slack $\Delta$ characterizes the difference between the conditional distribution of $W_1|Z=z$ and the distribution of $\psi(Z)$. In particular, it can be easily checked that if $(W_1,\ldots,W_B) \indep Z$, $(W_1,\ldots, W_B, \psi(Z))$ are IID and the distribution of $\psi(Z)$ is continuous then $\Delta = 0$. 

The next proposition connects the slack $ d_{\mathrm{KS}}(F_0(Z), U(0,1))$ to the Kolmogorov-Smirnov distance between the conditional distribution of $W_1|Z= z$ and the distribution of $\psi(Z)$. Define
\begin{equation}\label{eq:bootstrap-consistency}
\Delta^*(Z) = \sup_{x \in \mathbb{R}} | \mathbb P(W_1 \leq x|Z) - \mathbb P(\psi(Z) \leq x)|.
\end{equation}
For any real-valued random variable $Y$, define the L{\'e}vy concentration function as
\[
Q_Y(\varepsilon) := \sup_{a, b:\, b - a \le \varepsilon}\,\mathbb{P}(Y \in (a, b])\quad\mbox{for all}\quad \varepsilon > 0.
\]
The concentration function measures the discreteness of a distribution. If $Y$ has a bounded Lebesgue density, then $Q_Y(\varepsilon) \le C\varepsilon$ for some constant $C$ and if $Y$ has an atom, then $Q_Y(\varepsilon)$ is bounded away from zero as $\varepsilon \to 0.$
\begin{prop}
    \label{prop:ks_connection}
 Suppose there exists $c, \eta > 0$ such that for all small enough $\zeta > 0$, $Q_{\psi(Z)}(\zeta) \le c\zeta + \eta$. Then we have the following bound for any $\epsilon > 0$,
 \[
 d_{\mathrm{KS}}(F_0(Z), U(0,1)) \leq \epsilon + \mathbb{P}( \Delta^*(Z) > \epsilon) + \eta. 
 \]
\end{prop}
Refer to \Cref{app:proof_ks_connection} for the proof of Proposition~\ref{prop:ks_connection}. From the proof, it is easy to see that we only need
\[
\lim_{\varepsilon\to 0}\, Q_{\psi(Z)}(\varepsilon) \le \eta\quad\equiv \quad \sup_{x\in\mathbb{R}}\,\mathbb{P}(\psi(Z) = x) \le \eta.
\]
Because $\Delta^*(Z)\in[0, 1]$, we can apply Markov's inequality to get a simple corollary:
\begin{align*}
 d_{\mathrm{KS}}(F_0(Z), U(0,1)) &\leq \epsilon + \mathbb{P}( \Delta^*(Z) > \epsilon) + \eta\\ 
 &\leq \inf_{p \geq 1} \left\{\epsilon + \frac{\mathbb E[\Delta^*(Z)^p]}{\epsilon^p} \right\}+ \eta \\ 
 &= \inf_{p \geq 1} \left\{(p+1) \left( \frac{\|\Delta^*(Z) \|_p}{p}\right)^{\frac{p}{p+1}} \right\} + \eta,
\end{align*}
where $\|\Delta^*(Z) \|_p = ( \mathbb E[\Delta^*(Z)^p])^{1/p} $ is the $L^p$ norm of $\Delta^*(Z)$. Proposition~\ref{prop:ks_connection} also implies that if $\Delta^*(Z) \stackrel{P}{\rightarrow} 0$, then $d_{\mathrm{KS}}(F_0(Z), U(0,1)) \leq \eta$. Note that if $\Delta = 0$, there is a minimum $B$ that yields a non-trivial one-sided $(1-\alpha)$-valid interval for $\psi(Z)$. This minimum can be obtained by setting $a = b = 0$ in \Cref{thm:main_result} and observing that we need the following inequality to hold true for $(1-\alpha)$ coverage,
\[
 1 - \frac{1}{B+1} \geq 1 - \alpha \iff B \geq \frac{1}{\alpha}-1.
\]
In order words, $B$ must be at least $\lceil (1/\alpha) -1 \rceil$ for a non-trivial one-sided interval based on the order statistics  $(W_{(1)}, \ldots, W_{(B)})$ to exist that has $(1-\alpha)$ coverage. It can be easily checked that the minimal $B$ for Theorem~\ref{thm:main_result} to yield a non-trivial two-sided interval with $(1 - \alpha)$ coverage (when $\Delta = 0$) is $(2/\alpha) - 1$. 
\begin{remark}[Sharper coverage bounds]
    \label{rem:tighter_cov_bound}
The slack $\Delta$ appearing in the coverage bounds in \Cref{thm:main_result} can be made sharper. Let $\mathcal H = \{H | H :\Omega \mapsto \mathbb R \mbox{ is measurable and } \mathbb E[H^r] = 1/(r +1) \mbox{ for } r \in \{1,\ldots,B\}\} $ be the class of univariate random variables whose first $B$ raw moments match with that of $U(0,1)$. It can be shown that under the setting of \Cref{thm:main_result}, the following holds for any $0 \leq a < B-b \leq B$, 
\[
    - \inf_{H \in \mathcal H}\widetilde d_{\mathrm{KS}}(F_0(Z), H)  \leq \mathbb{P}\left(\psi(Z) \in [W_{(a)}, W_{(B-b)}] \right) -\left(  1- \frac{a+b +1}{B+1} \right) \leq  \frac{1}{B+1} + \inf_{H \in \mathcal H}\widetilde d_{\mathrm{KS}}(F_0(Z), H).
\]
Refer to \Cref{app:proof_tighter_cov} for a proof of this statement. The bounds on coverage can be analogously improved for the other statements in \Cref{thm:main_result}. 
\end{remark}
We can extend the coverage guarantees to the case when $W_1,\ldots,W_B$ are distributed independently upon conditioning on $Z$ i.e.\ $P_{W^{(B)}|Z}(z, \cdot) = P_{W_1|Z}(z, \cdot) \otimes \cdots \otimes P_{W_B|Z}(z, \cdot) $ for all $z \in \mathcal Z$. The following theorem obtains both lower and upper bounds on the coverage in terms of the conditional distributions $\{P_{W_i|Z}(z, \cdot)\}_{i =1}^B$. 
\begin{theorem}
\label{thm:main_result_extension}
Suppose $W_1,\ldots,W_B$ are $B$ random variables which conditional on another random vector $Z$ are distributed independently and $\psi(Z) \sim F$. For $1\le i\le B$, define $F_i(z) = \mathbb{P}(W_i \leq \psi(Z)| Z = z)$. Then we have the following bounds on coverage for any $0 \leq a < B-b \leq B$,  
\begin{equation*}
  -\delta  \leq \mathbb{P}(\psi(Z) \in [W_{(a)}, W_{(B-b)}]) - \left(  1- \frac{a+b +1}{B+1} \right) \leq \delta + \frac{1}{B+1}, 
\end{equation*}
where
\begin{equation*}
    \delta ~=~ \widetilde d_{\mathrm{KS}}(\overline F(Z), U(0,1)) +\left(\sum_{i = 1}^B \kappa_i^2\right)\left[I_B + \widetilde d_{\mathrm{KS}}(\overline F(Z), U(0,1))\right],
\end{equation*}
where $\kappa_i = \sup_{u \in \mathcal Z} | F(\psi(u)) - F_i(u)|$ for $i\in [B]$, $\overline F(z) = (1/B) \sum_{i = 1}^B F_i(z)$ and $I_B$ is defined as follows,
\[
 I_B = \begin{cases}
        1 & \mbox{if} \quad B \leq 4, \\
       1 - \sqrt{1 - \frac{4}{B}} + \frac{2}{B}\log \left( \frac{1 + \sqrt{1 - (4/B)}}{1 -\sqrt{1 - (4/B)}}\right) & \mbox{if} \quad B > 4.
       \end{cases}
\]
We set the order statistic $W_{(0)} = \min\{w|w \in \mathcal W\}$. 
\end{theorem}
\Cref{thm:main_result_extension} characterizes the coverage of $\psi(Z)$ by a confidence interval based on the order statistics $(W_{(1)}, \ldots, W_{(B)})$ when $W_1,\ldots W_B$ are independent conditioned on $Z$. The proof of \Cref{thm:main_result_extension} relies on a bound (\cite{ehm1991binomial}) on the total variation distance between $\mbox{Poi-Bin}(B, p_1,\ldots,p_B)$ and $\mbox{Bin}(B, \overline p_B)$ where $\overline p_B = (1/B) \sum_{i = 1}^B p_i$. To see the detailed proof of \Cref{thm:main_result_extension}, refer to \Cref{app:proof_main_res_ext}. We can prove versions of \Cref{thm:main_result_extension} similar to the second and third statements of \Cref{thm:main_result} by following the same argument as in the proof of \Cref{thm:main_result}. 
\begin{remark}[Behavior of the coverage slack in \Cref{thm:main_result_extension}]
\label{rem:cov_slackIn}
It can easily verified that $BI_B = 2\log B + 2 + o(1)$ as $B \rightarrow \infty$. In other words $I_B = O(\log B /B)$ as $B \rightarrow \infty$. Therefore $ \mathbb{P}(\psi(Z) \in [W_{(a)}, W_{(B-b)}])  \geq 1 - (a+b+1)/(B+1) -o(1)$ holds as $B \rightarrow \infty$ provided $\widetilde d_{\mathrm{KS}}(\overline F(Z), U(0,1)) \rightarrow 0$ and $\sum_{i = 1}^B \kappa_i^2 < \infty$ as $ B \rightarrow \infty$. In the special case when $W_i \indep Z$ we have $\kappa_i = d_{\mathrm{KS}}(W_i, \psi(Z))$. Therefore $\sum_{i = 1}^B d_{\mathrm{KS}}(W_i, \psi(Z))^2 < \infty$ is a necessary and sufficient condition for $\sum_{i = 1}^B \kappa_i^2 < \infty$ under the independence assumption $(W_1,\ldots,W_B) \indep Z$. 
\end{remark}
Based on \Cref{rem:cov_slackIn}, if $ d_{\mathrm{KS}}(\overline F(Z), U(0,1)) \rightarrow 0$ and $(\sum_{i = 1}^B \kappa_i^2)(\log B/B) \rightarrow 0$ then $[W_{(\lfloor (B+1) (\alpha/2) \rfloor)},$ $ W_{(\lceil (B+1)(1 -\alpha/2) \rceil)}]$ is an asymptotically valid $(1-\alpha)$ interval of $\psi(Z)$. However if $d_{\mathrm{KS}}(\overline F(Z), U(0,1)) \rightarrow 0$ and $(\sum_{i = 1}^B \kappa_i^2)(\log B/B) $ is bounded away from $0$, \Cref{thm:main_result_extension} does not provide a method to choose an asymptotically valid $(1-\alpha)$ confidence interval of $\psi(Z)$ when $W_1,\ldots W_B$ are independent conditioned on $Z$. The next theorem solves this problem using stochastic ordering results between $\mbox{Poi-Bin}(B, p_1,\ldots,p_B)$ and $\mbox{Bin}(B, \overline p_B)$ random variables from \cite{hoeffding1956distribution}; see also \cite{tang2023poisson}. 
\begin{theorem}
    \label{thm:finer_main_thm_ext}
Suppose $W_1,\ldots,W_B$ are $B$ random variables which conditional on another random vector $Z$ are distributed independently. Then we have the following lower bound on coverage for any $0 \leq a < B-b \leq B$,
\begin{equation*}
     \mathbb{P}\left(\psi(Z) \in [W_{(a)}, W_{(B-b)}] \right) \geq 1 - \frac{3(a+b+1)}{2B} - 6 d_{\mathrm{KS}}(\overline F(Z), U(0,1)) ,
\end{equation*}
where $F_i(z) = \mathbb{P}(W_i \leq \psi(Z)| Z = z)$ for $i\in [B]$ and $\overline F(z) = (1/B) \sum_{i = 1}^B F_i(z)$. We set the order statistic $W_{(0)} = \min\{w|w \in \mathcal W\}$. 
\end{theorem}
\Cref{thm:finer_main_thm_ext} implies that if we set $a = b= \lfloor (B\alpha/3) - (1/2) \rfloor $, the confidence interval $[W_{(\lfloor (B\alpha/3) - (1/2) \rfloor )},$ $ W_{(B - \lfloor (B\alpha/3) - (1/2) \rfloor )}]$ is an asymptotically valid $(1-\alpha)$ confidence interval of $\psi(Z)$ when $W_1,\ldots W_B$ are independent conditioned on $Z$ and $ d_{\mathrm{KS}}(\overline F(Z), U(0,1)) \rightarrow 0$. The proof of \Cref{thm:finer_main_thm_ext} is discussed in \Cref{app:proof_mainthm_finer}. We note that extensions of Bernstein-Hoeffding method (see Theorem $1.8$ in \cite{pelekis2015bernstein}) can also be used to obtain a lower bound on the coverage of the confidence interval $[W_{(a)}, W_{(B-b)}]$ in the set-up of \Cref{thm:main_result_extension}. 

The coverage guarantees of the confidence interval based on the order statistics $(W_{(1)}, \ldots, W_{(B)})$ can be studied in more general settings. Setting $V =(W_1,\ldots, W_B, \psi(Z))$, we let 
\[
V^i ~=~ (W_1,\ldots,W_{i-1},\psi(Z),W_{i+1},\ldots,W_B,W_i),
\]
where $\psi(Z)$ is swapped with $W_i$. We set $V^{B+1} = V$. We define $\Gamma(W_1,\ldots,W_B,\psi(Z))$ as,
\begin{equation}
\label{eq:gamma_defn}
    \Gamma(W_1,\ldots,W_B,\psi(Z)) = \sup_{A } \left|\mathbb P(V \in A) -  \frac{1}{B+1} \sum_{i = 1}^{B+1} \mathbb P(V^i \in A) \right|,
\end{equation} 
where $V = (W_1,\ldots, W_B, \psi(Z))$ and the supremum is taken over all measurable sets. This quantity $\Gamma(\cdot)$ measures how far is the vector $V$ from an exchangeable random vector in terms of the total variation distance.

The following theorem characterizes the coverage of $\psi(Z)$ using $[W_{(\lfloor (B+1)(\gamma/2)  \rfloor - 1)}, W_{(\lceil (B+1)(1 - (\beta/2))\rceil})]$ for any $0< \gamma, \beta < 1$. 
\begin{theorem}
    \label{thm:beyond_exchangeability}
For any $0 < \gamma, \beta < 1$,
\begin{equation*}
   \mathbb P\left(\psi(Z) \in \left[W_{(\lfloor (B+1)(\gamma/2)  \rfloor - 1)}, W_{(\lceil (B+1)(1 - (\beta/2))\rceil})\right] \right)  - (1 - (\gamma + \beta)/2)\geq  -\delta,
\end{equation*}
where $ \delta = \Gamma(W_1,\ldots,W_B,\psi(Z))$. We set $W_{(r)} = \min\{w|w \in \mathcal W\}$ for $r \leq 0$ and $W_{(r)} = \max\{w|w \in \mathcal W\}$ for $r \geq B+1$. Additionally if the distributions of the random variables $\psi(Z), \{W_i\}_{i = 1}^B$ are continuous then we have,
\[
P\left(\psi(Z) \in \left[W_{(\lfloor (B+1)(\gamma/2)  \rfloor - 1)}, W_{(\lceil (B+1)(1 - (\beta/2))\rceil})\right] \right)  - (1 - (\gamma + \beta)/2) \leq  \delta +  \frac{4}{B+1}.
\]
\end{theorem}
\Cref{thm:beyond_exchangeability} implies that the confidence interval $ [W_{(\lfloor (B+1)(\alpha/2)  \rfloor - 1)}, W_{(\lceil (B+1)(1 - (\alpha/2))\rceil})]$ is an asymptotically valid $(1-\alpha)$ confidence interval of $\psi(Z)$ if $\Gamma(W_1,\ldots,W_B,\psi(Z)) \rightarrow 0$. A sufficient condition for the slack $\Gamma(W_1,\ldots,W_B,\psi(Z))$ to vanish is that the $(B+1)$ dimensional vector $(W_1,\ldots, W_B, \psi(Z))$ be exchangeable. The proof of \Cref{thm:beyond_exchangeability} is discussed in \Cref{app:proof_beyond_exchangeability}.   

\section{Applications}
\label{sec:applications}
This section illustrates several applications of the master theorems introduced in \Cref{sec:master_thms}, with examples drawn from bootstrap, permutation tests, conformal prediction, and subsampling. The key message is that one can obtain asymptotically valid confidence intervals for the target functional even when the number of resamples does not tend to infinity.
\paragraph{Code availability.}
All code for this project, including the implementation of the methods, the simulation scripts, and the code used to generate the figures in the paper, is publicly available at \url{https://github.com/manitpaul/cheap_bootstrap_project}. The repository contains all materials needed to reproduce the computational results and figures reported in this manuscript.
\subsection{Non-parametric bootstrap}
\label{subsec:boot}
Consider $m$ IID observations $D_m = (X_1,\ldots, X_m)$ drawn from some distribution $P$. The goal is to construct an asymptotically valid $(1-\alpha)$ confidence interval of a functional $\theta_0 = \phi(P) \in \mathbb R^d$. Moreover let $S:\mathbb R^d \mapsto \mathbb R$ be a function that maps from $\mathbb R^d$ to $\mathbb R$. Suppose $\widehat \theta_m = \phi(P_m) $ is the observed sample estimate of $\theta_0$ (where $P_m$ is the empirical distribution of $(X_1,\ldots,X_m)$) and suppose $S_m(\theta) = S(\tau_m(\widehat \theta_m - \theta))$ is a specified root (\cite{beran1988prepivoting, beran2025bootstrap}), a real valued function of $(X_1,\ldots,X_m)$ and $\theta$ (where $\tau_m$ is the rate of convergence of $\widehat \theta_m$). We draw $B$ bootstrap samples $\{(X_1^{*b}, \ldots, X_m^{*b})\}_{b = 1}^B$ of size $m$ from $P_m$. Let $\widehat \theta_m^{*b} = \phi(P_m^{*b})$ (for $b \in [B]$) be the estimate of $\theta_0$ based on $P_m^{*b}$, the empirical distribution of the bootstrap sample $(X_1^{*b}, \ldots, X_m^{*b})$ and let $W_b = S(\tau_m(\widehat \theta_m^{*b}  - \widehat \theta_m))$. Let $(W_{(1)}, \ldots, W_{(B)})$ be the order statistics of the estimates $\{W_b\}_{b = 1}^B$ from the bootstrap samples. The vanilla non-parametric bootstrap confidence interval $\mathrm{CI}^{\mathtt{vanilla-boot}}_{m,B,\alpha} =\{\theta: S_m(\theta) \in [W_{(\lceil B(\alpha/2)\rceil)}, W_{(\lceil B(1 - (\alpha/2))\rceil)}) \}$ (discussed in \cite{hall1986number}, \cite{hall1988theoretical}, \cite{hall2013bootstrap}) has been shown to attain the following coverage guarantee under bootstrap consistency assumption,
\[
\mathbb P\left(\theta_0 \in \mathrm{CI}^{\mathtt{vanilla-boot}}_{m,B,\alpha}\right) \geq 1- \alpha - O(1/B) - o(1), \quad \mbox{as} \quad m \rightarrow \infty. 
\]
Therefore we need the number of bootstrap samples $B$ to go to infinity to ensure an asymptotic coverage of $(1-\alpha)$ for $\mathrm{CI}^{\mathtt{vanilla-boot}}_{m,B,\alpha}$. To remedy this issue, we introduce the modified non-parametric bootstrap confidence interval $\mathrm{CI}^{\mathtt{mod-boot}}_{m,B,\alpha} = \{\theta: S_m(\theta) \in [W_{(\lfloor (B+1)(\alpha/2)\rfloor)},  W_{(\lceil (B+1)(1 - \alpha)\rceil + \lfloor (B+1)(\alpha/2)\rfloor)}) \}$. An application of the theorems proved in \Cref{sec:master_thms} provides the following corollary.
\begin{corollary}
    \label{cor:boot}
    Using the stated notations, the modified non-parametric bootstrap confidence interval $\mathrm{CI}^{\mathtt{mod-boot}}_{m,B,\alpha}$ satisfies the following coverage guarantee,
    \[
    \left|P( \theta_0 \in \mathrm{CI}^{\mathtt{mod-boot}}_{m,B,\alpha} )  - \frac{\lceil(B+1)(1-\alpha) \rceil}{B+1} \right| \leq \widetilde d_{\mathrm{KS}}(F_0(D_m), U(0,1)),
    \]
where $F_0(D_m) = \mathbb P(W_b  \leq S_m(\theta_0) | D_m)$. The bounds on coverage can be restated as,
    \[
    - \widetilde d_{\mathrm{KS}}(F_0(D_m), U(0,1)) \leq \mathbb P( \theta_0 \in \mathrm{CI}^{\mathtt{mod-boot}}_{m,B,\alpha} ) - (1 - \alpha) \leq \frac{1}{B+1} +  \widetilde d_{\mathrm{KS}}(F_0(D_m), U(0,1)) ,
    \]
\end{corollary}
\Cref{cor:boot} establishes the lower and upper bounds on the coverage of the confidence interval $\mathrm{CI}^{\mathtt{mod-boot}}_{m,B,\alpha}$. Since the slack term $\widetilde d_{\mathrm{KS}}(F_0(D_m), U(0,1))$ does not have any dependence on the number of bootstrap samples, the asymptotic validity of $\mathrm{CI}^{\mathtt{mod-boot}}_{m,B,\alpha}$ holds for any $B \geq (2/\alpha) - 1$ (for a non-trivial confidence interval) if $\widetilde d_{\mathrm{KS}}(F_0(D_m), U(0,1)) \rightarrow 0$ as $m \rightarrow \infty$. The first statement of \Cref{cor:boot} also implies that if $(B+1)(1-\alpha) \in \mathbb N$ is an integer, then we have exact coverage of $(1-\alpha)$ when $\widetilde d_{\mathrm{KS}}(F_0(D_m), U(0,1)) \rightarrow 0$ as $m \rightarrow \infty$. The proof of \Cref{cor:boot} is provided in \Cref{app:proof_cor_boot}. We note that similar coverage guarantees can likewise be derived for alternative weighted bootstrap procedures (see \cite{chatterjee2005generalized}) that use weights other than the usual multinomial ones. In light of Proposition~\ref{prop:ks_connection}, if the distribution of $S(\tau_m(\widehat \theta_m - \theta_0))$ is continuous, one sufficient condition for $\widetilde d_{\mathrm{KS}}(F_0(D_m), U(0,1)) \rightarrow 0$ as $m \to \infty$ is the standard bootstrap consistency assumption \cite{hall2013bootstrap},
\[
\Delta^*(D_m) = \sup_{x \in \mathbb R} | \mathbb P(S(\tau_m(\widehat \theta_m^{*1} - \widehat \theta_m)) \leq x| D_m) - \mathbb P(S(\tau_m(\widehat \theta_m - \theta_0 ))\leq x)| \stackrel{P}{\rightarrow} 0, \quad \mbox{as} \quad m \rightarrow \infty.
\] 
\begin{remark}[Randomized modified bootstrap confidence interval]
    \label{rem:random_boot}We can define a randomized modified bootstrap confidence interval that guarantees exact asymptotic coverage of $(1-\alpha)$ for any $B \geq (2/\alpha) - 1$ provided $\widetilde d_{\mathrm{KS}}(F_0(D_m), U(0,1)) \rightarrow 0$ as $m \to \infty$. Generate a Uniform $(0,1)$ random variable $U$ and set,
  \begin{equation*}
  \begin{split}
      \mathrm{CI}^{\mathtt{rand-mod-boot}}_{m,B,\alpha} =&\begin{cases} \{\theta: S_m(\theta) \in [W_{(\lfloor (B+1)(\alpha/2)\rfloor)},  W_{(\lceil (B+1)(1 - \alpha)\rceil + \lfloor (B+1)(\alpha/2)\rfloor)}) \} \quad & \mbox{if} \quad U \leq \tau_{\alpha, B} \\
      \{\theta: S_m(\theta) \in [W_{(\lfloor (B+1)(\alpha/2)\rfloor)},  W_{(\lfloor (B+1)(1 - \alpha)\rfloor + \lfloor (B+1)(\alpha/2)\rfloor)}) \} \quad & \mbox{if} \quad U > \tau_{\alpha, B}
      \end{cases} \quad \mbox{where}, \\
  \tau_{\alpha, B} =& \begin{cases}
        1 \quad & \mbox{if} \quad (B+1)(1-\alpha) \in \mathbb N, \\
        \frac{(1-\alpha) - (\lfloor (B+1)(1-\alpha) \rfloor/ (B+1))}{ (\lceil (B+1)(1-\alpha) \rceil - \lfloor (B+1)(1-\alpha) \rfloor)/ (B+1)} \quad  &\mbox{otherwise}.
    \end{cases}
    \end{split}
  \end{equation*}  
We can show that $\mathrm{CI}^{\mathtt{rand-mod-boot}}_{m,B,\alpha}$ enjoys the following coverage guarantee,
 \[
    \left|P( \theta_0 \in \mathrm{CI}^{\mathtt{rand-mod-boot}}_{m,B,\alpha} )  - (1-\alpha)\right| \leq \widetilde d_{\mathrm{KS}}(F_0(D_m), U(0,1)).
    \]
Similar randomized confidence intervals that guarantee exact asymptotic coverage of $(1-\alpha)$ can be defined for other examples (such as subsampling (\Cref{subsec:subsample})) studied in this paper. The proof of the coverage guarantee for the randomized confidence interval can be seen in \Cref{app:proof_rand_boot}.
\end{remark}
We compare the performance of the modified bootstrap confidence interval $\mathrm{CI}^{\mathtt{mod-boot}}_{m,B,\alpha}$ and the randomized modified bootstrap confidence interval $\mathrm{CI}^{\mathtt{rand-mod-boot}}_{m,B,\alpha}$ with that of the vanilla bootstrap confidence interval $\mathrm{CI}^{\mathtt{vanilla-boot}}_{m,B,\alpha}$ and the cheap bootstrap confidence interval (introduced in \cite{lam2022cheap}) in terms of both coverage and mean-width in the following two settings:
\paragraph{Setting $1$}:\label{para:setting:1} $X_1,\ldots, X_m \stackrel{iid}{\sim} \mbox{Exp}(5)$, $\theta_0 = \mathbb E[\mbox{Exp}(5)] = 1/5$, $\widehat \theta_m = (1/m)\sum_{i = 1}^m X_i$, and $\widehat \theta_m^{*b} = (1/m)\sum_{i = 1}^m X_i^{*b}$ for $b \in [B]$, $\tau_m = \sqrt{m}$, and $S:\mathbb R \mapsto \mathbb R$ is defined as $S( x) = x$. 
\paragraph{Setting $2$}:\label{para:setting:2} $X_1,\ldots, X_m \stackrel{iid}{\sim} X \boldsymbol{Y}$ where $X \sim t_5$ and $\boldsymbol{Y} = (Y_1,\ldots, Y_{100})$ where $Y_1,\ldots Y_{100} \stackrel{iid}{\sim} \mathcal X^2_1$, $\theta_0 = \mathbb E[(X \boldsymbol{Y})_j] = \boldsymbol 0$, $\widehat \theta_m = (1/m)  \sum_{i = 1}^m X_{i} $, $\widehat \theta_m^{*b} = (1/m) \sum_{i = 1}^m X_{i}^{*b} $ for $b \in [B]$, $\tau_m = \sqrt{m}$, and $S:\mathbb R^{100} \mapsto \mathbb R$ is defined as $S(\boldsymbol x) = \| \boldsymbol x \|_{\infty}$. 
\begin{figure}[h]
    \centering
    \begin{minipage}[t]{0.48\textwidth}
        \centering
        \includegraphics[width=\linewidth]{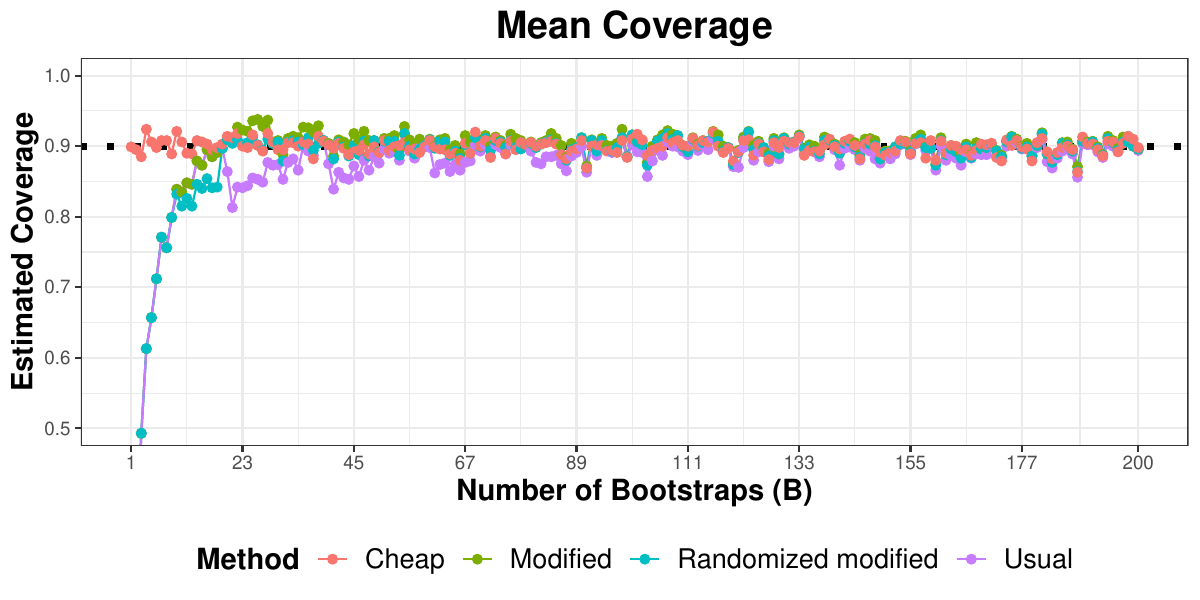}
    \end{minipage}
    \hfill
    \begin{minipage}[t]{0.48\textwidth}
        \centering
        \includegraphics[width=\linewidth]{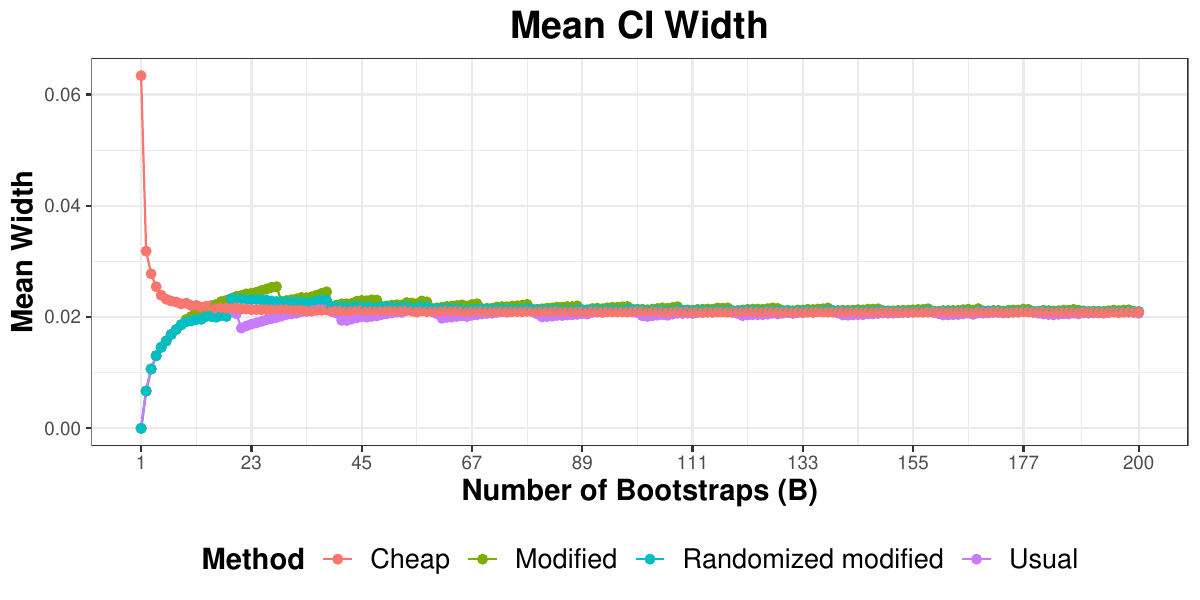}
    \end{minipage}
    \caption{Comparison of the coverage and the mean width between {usual}, {modified}, randomized modified, and {cheap} non-parametric bootstrap confidence interval in Setting~1 at level $1 - \alpha = 0.9$. The black dotted horizontal line represents the nominal level of $0.9$.}
    \label{fig:univariate_boot}
\end{figure}
\begin{figure}[h]
    \centering
    \begin{minipage}[t]{0.48\textwidth}
        \centering
        \includegraphics[width=\linewidth]{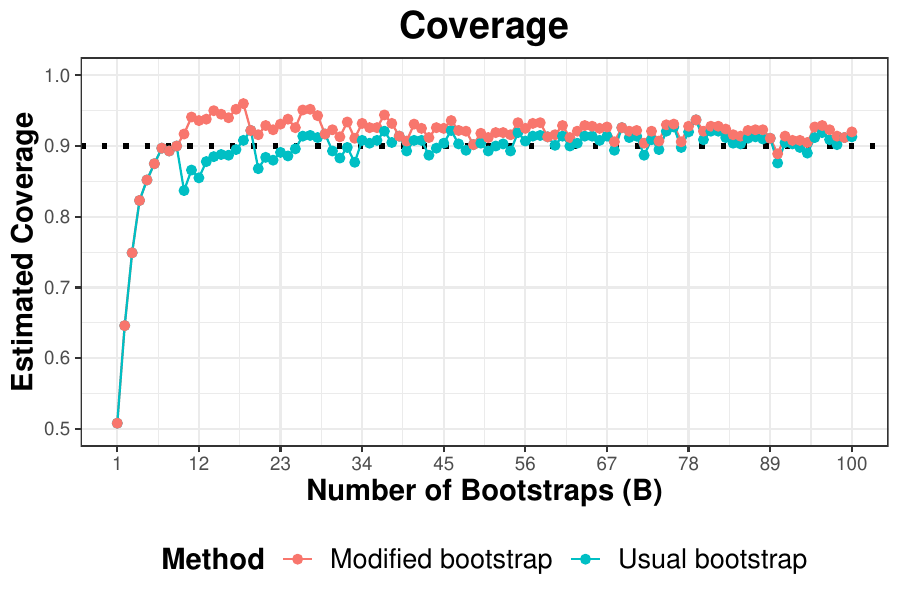}
    \end{minipage}
    \hfill
    \begin{minipage}[t]{0.48\textwidth}
        \centering
        \includegraphics[width=\linewidth]{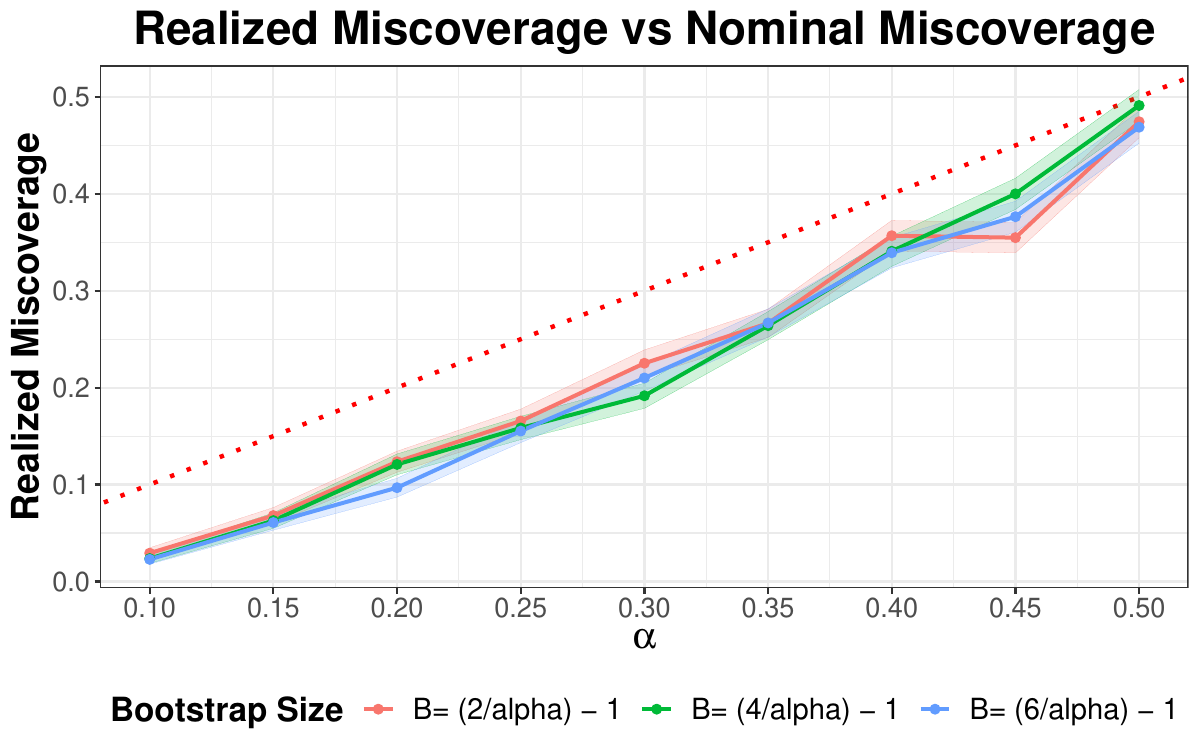}
    \end{minipage}
     \caption{The left plot shows the comparison of the coverage between {usual} and {modified} non-parametric bootstrap confidence interval in the high dimensional setting~2 at level $1 - \alpha = 0.9$ for sample size $m = 1000$. The black dotted horizontal line represents the nominal level of $0.9$. The right plot shows the realized miscoverage of the modified non-parametric bootstrap confidence interval against the nominal miscoverage $(\alpha)$ for sample size $m = 100$ in the high dimensional setting~2 as $\alpha \in [0,1, 0.5]$. The performance of the modified bootstrap confidence interval is shown for three different bootstrap budgets (depending on $\alpha$), $B \in\{\lceil 2 /\alpha\rceil -1, \lceil 4 /\alpha\rceil -1, \lceil 6 /\alpha\rceil -1 \} $. The red dotted line depicts the ideal miscoverage level. }
    \label{fig:multivariate_boot}
\end{figure}

For the univariate simulation study, we vary the number of bootstrap samples $B$ in the set $\{1,\ldots,200\}$ and for the multivariate simulation study (left plot) we vary the number of bootstrap samples $B$ in the set $\{1,\ldots,100\}$. We set the sample size $m = 100$ for the univariate case, and set the sample size $m = 1000$ for the multivariate case (left plot). We compute the coverage and mean width of the confidence interval over $1000$ replications. We can see from both figures~\ref{fig:univariate_boot} and \ref{fig:multivariate_boot} that the usual non-parametric bootstrap confidence interval attains the coverage of $0.9$ as $B$ increases, but does not have valid coverage for small $B$. The modified bootstrap confidence interval attains the coverage of $0.9$ for $B \geq (2/\alpha) - 1 = 19$ in both figures~\ref{fig:univariate_boot} and \ref{fig:multivariate_boot}. The cheap bootstrap confidence interval attains the required coverage of $0.9$ for $B\geq 1$ in setting~1 when the asymptotic distribution of the test statistic is gaussian, but is undefined for $B \leq 100$ in setting~2 as the $(1-\alpha)$ quantile of the Hotelling $T^2$ distribution with parameters $100$ and $B$ is not defined for $B \leq 100$ (refer to section $A.2$ in \cite{lam2022cheap}). More generally even if $\tau_m(\widehat \theta_m - \theta_0)$ and $\tau_m(\widehat \theta_m^{*b} - \widehat \theta_m)$ (conditioned on the data) converge to mean zero gaussian, the cheap bootstrap confidence interval is not defined in the regime $B \leq d$. 

The simulation study in setting~2 goes on to illustrate that the modified bootstrap confidence interval $\mathrm{CI}^{\mathtt{mod-boot}}_{m,B,\alpha}$ provides asymptotic validity for any $B \geq (2/\alpha) -1$ even in high-dimensional regimes as long as we are dealing with a univariate projection of the test statistic and the standard bootstrap consistency condition holds (\cite{chernozhukov2013gaussian}). We observe from the right plot in figure~\ref{fig:multivariate_boot} that even when $m = 100$ (in setting~2), $\mathrm{CI}^{\mathtt{mod-boot}}_{m,B,\alpha}$ guarantees valid coverage over the entire range of $\alpha \in [0.1, 0.5]$ for all three bootstrap budgets. It can be seen from both the figures (and also from the definitions) that the modified bootstrap confidence interval is slightly wider than the vanilla bootstrap confidence interval. We observe from figure~\ref{fig:univariate_boot} that the randomized modified bootstrap confidence interval (\Cref{rem:random_boot}) attains exact coverage of $(1 -\alpha) = 0.9$ for $B \geq (2/\alpha) - 1$. Because of the randomization, the width of $\mathrm{CI}^{\mathtt{rand-mod-boot}}_{m,B,\alpha}$ is smaller than that of $\mathrm{CI}^{\mathtt{mod-boot}}_{m,B,\alpha}$ and is only slightly wider than $\mathrm{CI}^{\mathtt{vanilla-boot}}_{m,B,\alpha}$. The cheap bootstrap confidence interval can be observed to be significantly wider than the other confidence intervals for small $B$ and of comparable width for larger values of $B$. 

\subsection{Subsampling}
\label{subsec:subsample}
As in the earlier subsection, we consider $m$ IID observations $D_m = (X_1,\ldots, X_m)$ drawn from some distribution $P$. The goal is to construct an asymptotically valid $(1-\alpha)$ confidence interval of a functional $\theta_0 = \phi(P)$. Let $S:\mathbb R^d \mapsto \mathbb R$ be a function that maps from $\mathbb R^d$ to $\mathbb R$. Suppose $\widehat \theta_m = \phi(P_m) $ is the observed sample estimate of $\theta_0$ (where $P_m$ is the empirical distribution of $(X_1,\ldots,X_m)$). Let $\tau_m$ be the normalizing sequence such that $S(\tau_m(\widehat \theta_m - \theta_0))$ converges in distribution to $J$ (say) as $m \rightarrow \infty$ and let $S_m(\theta) = S(\tau_m(\widehat \theta_m - \theta))$. We draw $B$ subsamples $\{(X_1^{*b}, \ldots, X_k^{*b})\}_{b = 1}^B$ of size $k$ from $D_m$. Let $\widehat \theta_k^{*b} = \phi(P_k^{*b})$ (for $b \in [B]$) be the estimate of $\theta_0$ based on $P_k^{*b}$, the empirical distribution of the subsample $(X_1^{*b}, \ldots, X_k^{*b})$ and with some abuse of notation let $W_b = S(\tau_k(\widehat \theta_k^{*b} - \widehat \theta_m))$. Let $(W_{(1)}, \ldots, W_{(B)})$ be the order statistics of the estimates $\{W_b\}_{b = 1}^B$ from the subsamples. The vanilla subsampling confidence interval $\mathrm{CI}^{\mathtt{vanilla-subsample}}_{m,k,B,\alpha} = \{\theta: S_m(\theta) \in [W_{(\lceil B(\alpha/2)\rceil)}, W_{(\lceil B(1 - (\alpha/2))\rceil)}) \}$ (discussed in \cite{romano1999subsampling}, \cite{politis1999subsampling}, \cite{politis2001asymptotic}) requires the number of subsamples $B$ to go to infinity to ensure an asymptotic coverage of $(1-\alpha)$ under the subsampling consistency condition $\mathbb P(S(\tau_k(\widehat \theta_k^{*(b)} - \widehat \theta_m)) \leq x | D_m) \stackrel{P}{\rightarrow} J(x)$ at all continuity points of $J$. 

To remedy this issue, we introduce the modified subsampling confidence interval $\mathrm{CI}^{\mathtt{mod-subsample}}_{m,k, B,\alpha} = \{\theta: S_m(\theta) \in [W_{(\lfloor (B+1)(\alpha/2)\rfloor)},  W_{(\lceil (B+1)(1 - \alpha)\rceil + \lfloor (B+1)(\alpha/2)\rfloor)}) \}$. An application of the theorems proved in \Cref{sec:master_thms} provides the following corollary.
\begin{corollary}
    \label{cor:subsample}
    Using the stated notations, the modified subsampling confidence interval $\mathrm{CI}^{\mathtt{mod-subsample}}_{m,k,B,\alpha} $ satisfies the following coverage guarantee,
      \[
    \left|P( \theta_0 \in \mathrm{CI}^{\mathtt{mod-subsample}}_{m,k,B,\alpha} )  - \frac{\lceil(B+1)(1-\alpha) \rceil}{B+1} \right| \leq \widetilde d_{\mathrm{KS}}(F_0(D_m), U(0,1)),
    \]
    where $F_0(D_m) = \mathbb P(W_b  \leq S_m(\theta_0) | D_m)$. The bounds on coverage can be restated as,
    \[
    - \widetilde d_{\mathrm{KS}}(F_0(D_m), U(0,1)) \leq \mathbb P( \theta_0 \in \mathrm{CI}^{\mathtt{mod-subsample}}_{m,k,B,\alpha} ) - (1 - \alpha) \leq \frac{1}{B+1} +  \widetilde d_{\mathrm{KS}}(F_0(D_m), U(0,1)) .
    \]
\end{corollary}
\Cref{cor:subsample} establishes the lower and upper bounds on the coverage of the confidence interval $\mathrm{CI}^{\mathtt{mod-subsampling}}_{m,k,B,\alpha}$. Since the slack term $\widetilde d_{\mathrm{KS}}(F_0(D_m), U(0,1))$ does not have any dependence on the number of subsamples, the asymptotic validity of $\mathrm{CI}^{\mathtt{mod-subsample}}_{m,k,B,\alpha}$ holds for any $B \geq (2/\alpha) - 1$ (for a non-trivial confidence interval) if $\widetilde d_{\mathrm{KS}}(F_0(D_m), U(0,1)) \rightarrow 0$ as $k,m \rightarrow \infty$. The first statement of \Cref{cor:subsample} also implies that if $(B+1)(1-\alpha) \in \mathbb N$ is an integer, then we have exact coverage of $(1-\alpha)$ when $\widetilde d_{\mathrm{KS}}(F_0(D_m), U(0,1)) \rightarrow 0$ as $m \rightarrow \infty$. The proof of \Cref{cor:subsample} is analogous to that of \Cref{cor:boot} and is therefore omitted for brevity. In light of Proposition~\ref{prop:ks_connection}, if the distribution of $S(\tau_m(\widehat \theta_m - \theta_0))$ is continuous, one sufficient condition for $\widetilde d_{\mathrm{KS}}(F_0(D_m), U(0,1)) \rightarrow 0$ as $k,m \to \infty$ is the standard subsampling consistency assumption,
\[
\Delta^*(D_m) = \sup_{x \in \mathbb R} | \mathbb P(S(\tau_k(\widehat \theta_k^{*(b)} - \widehat \theta_m))\leq x| D_m) - \mathbb P(S(\tau_m(\widehat \theta_m - \theta_0)) \leq x)| \stackrel{P}{\rightarrow} 0, \quad \mbox{as} \quad k, m \rightarrow \infty.
\] 
We note that it is not necessary to know the rate of convergence of the estimator to construct the subsampling confidence interval. There are a number of works (see for example \cite{bertail1999subsampling}) in the literature that provide methods to obtain an estimate $\widehat \tau_m$ of the unknown rate of convergence $\tau_m$. Similar coverage guarantees as in \Cref{cor:subsample} can be shown for the modified subsampling confidence interval $\mathrm{CI}^{\mathtt{mod-subsample}}_{m,k,B,\alpha}$ with $\tau_k, \tau_m$ replaced by the estimated rates $\widehat \tau_k, \widehat \tau_m$. We compare the performance of the modified subsampling confidence interval $\mathrm{CI}^{\mathtt{mod-subsample}}_{m,k,B,\alpha}$ with that of the vanilla subsampling confidence interval $\mathrm{CI}^{\mathtt{vanilla-subsample}}_{m,k,B,\alpha}$ in terms of both coverage and mean-width in setting~$2$ (\ref{para:setting:2}) and the following setting:
\paragraph{Setting $3$}:\label{para:setting:3} $X_1,\ldots, X_m \stackrel{iid}{\sim} \mbox{Uniform}(0, 1)$, $\theta_0 = 1$, $\widehat \theta_m = \max\{ X_i: i \in [m]\}$, and $\widehat \theta_k^{*b} = \max\{ X_i^{*b}: i \in [k] \}$ for $b \in [B]$, $\tau_m = \sqrt{m}$, $\tau_m = m$, and $S:\mathbb R \mapsto \mathbb R$ is defined as $S( x) = x$. . 
\begin{figure}[h]
    \centering
    \begin{minipage}[t]{0.48\textwidth}
        \centering
        \includegraphics[width=\linewidth]{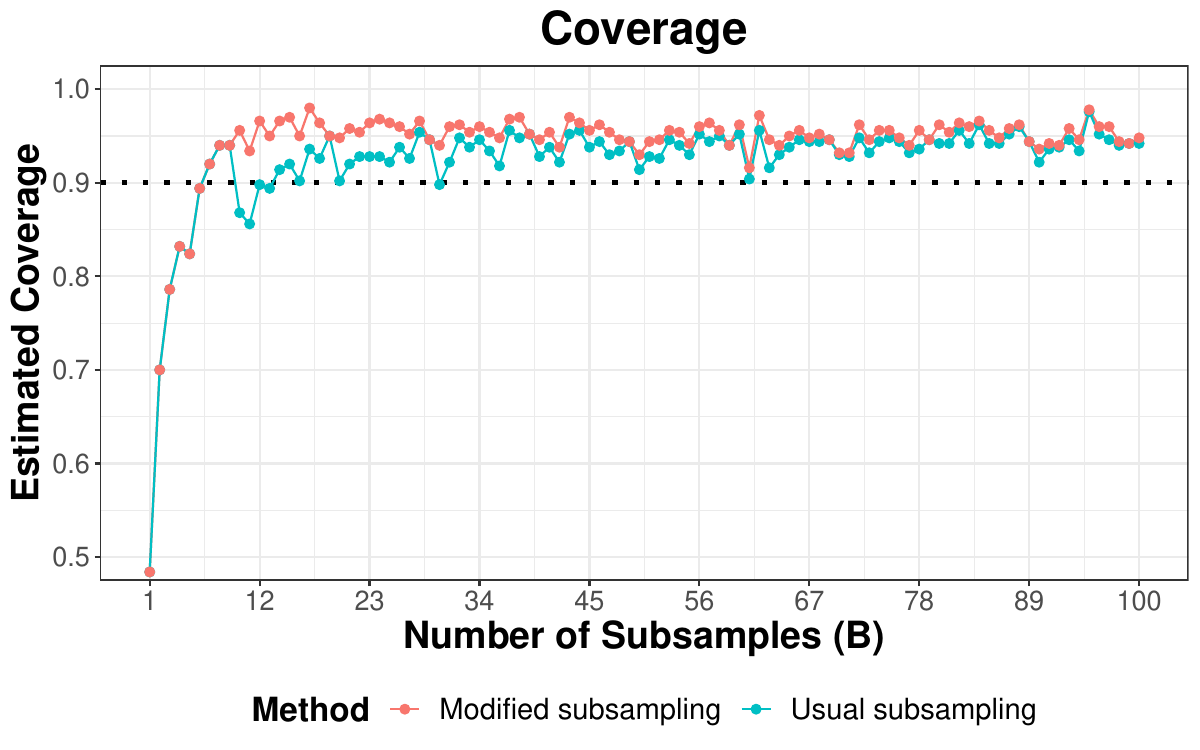}
    \end{minipage}
    \hfill
    \begin{minipage}[t]{0.48\textwidth}
        \centering
        \includegraphics[width=\linewidth]{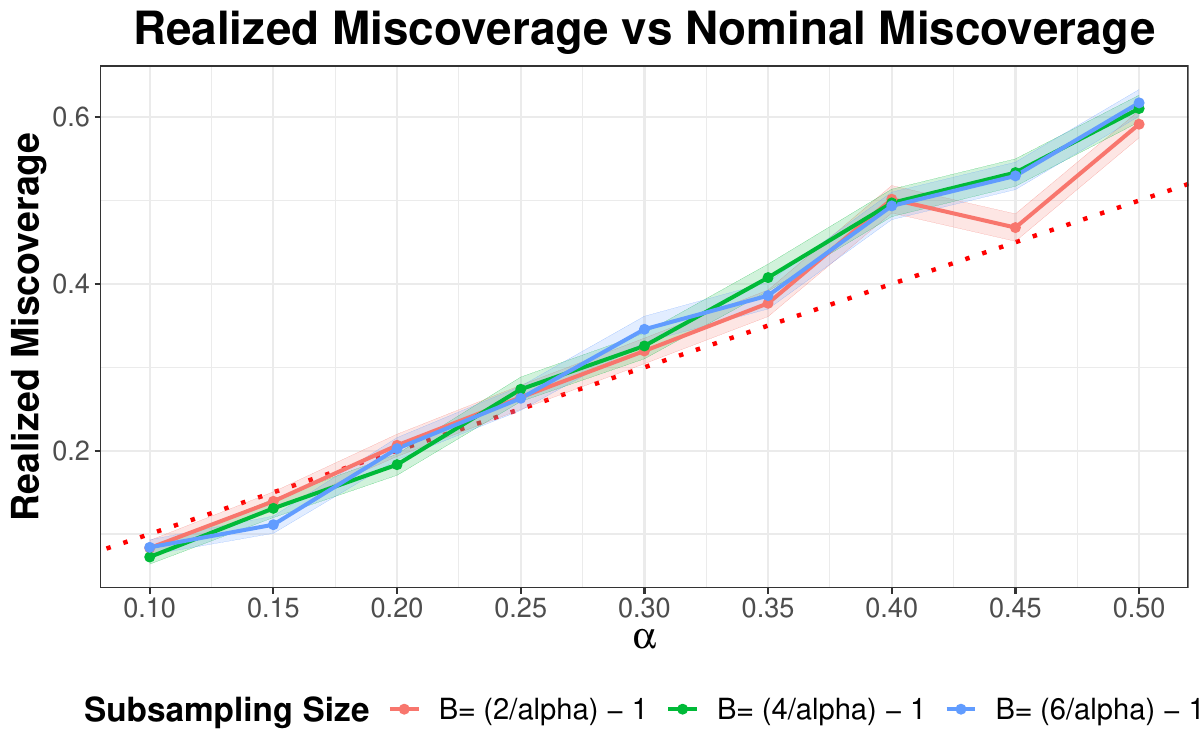}
    \end{minipage}
    \caption{The left plot shows the comparison of the coverage between {usual} and {modified} subsampling confidence interval in the high dimensional setting~2 at level $1 - \alpha = 0.9$ for sample size $m = 1000$. The black dotted horizontal line represents the nominal level of $0.9$. The right plot shows the realized miscoverage of the modified subsampling confidence interval against the nominal miscoverage $(\alpha)$ for sample size $m = 100$ in the high dimensional setting~2 as $\alpha \in [0,1, 0.5]$. The performance of the modified subsampling confidence interval is shown for three different subsampling budgets (depending on $\alpha$), $B \in\{\lceil 2 /\alpha\rceil -1, \lceil 4 /\alpha\rceil -1, \lceil 6 /\alpha\rceil -1 \} $. The red dotted line depicts the ideal miscoverage level. }
    \label{fig:multivariate_subsample}
\end{figure}
\begin{figure}[h]
    \centering
    \begin{minipage}[t]{0.48\textwidth}
        \centering
        \includegraphics[width=\linewidth]{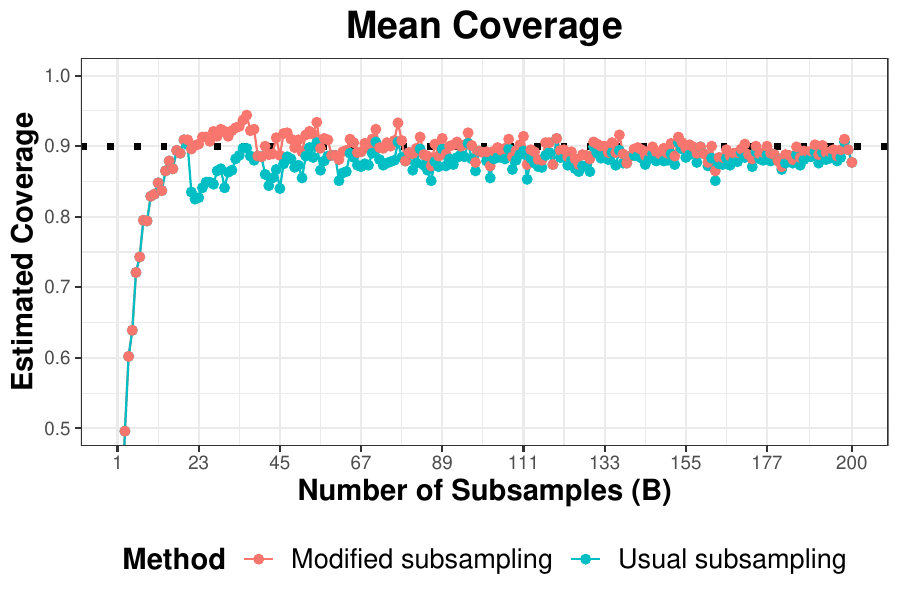}
    \end{minipage}
    \hfill
    \begin{minipage}[t]{0.48\textwidth}
        \centering
        \includegraphics[width=\linewidth]{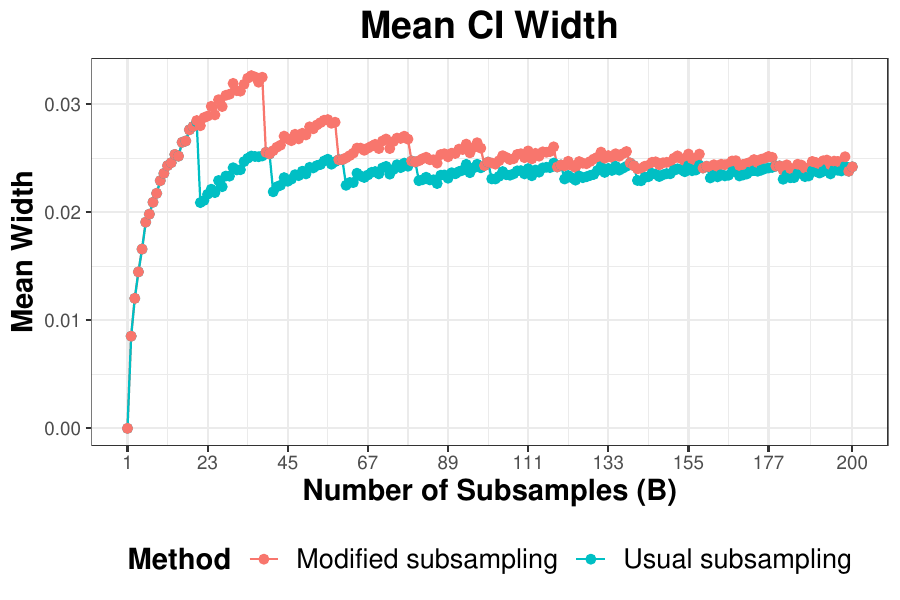}
    \end{minipage}
    \caption{Comparison of the coverage and the mean width between {usual} and modified subsampling confidence interval in Setting~3 at level $1 - \alpha = 0.9$. The dotted horizontal line represents the nominal level of $0.9$. }
    \label{fig:univariate_subsample}
\end{figure}

For the univariate simulation study (setting~$3$), we vary the number of subsamples $B$ in the set $\{1,\ldots,200\}$ and for the multivariate simulation study (setting~$2$) we vary the number of subsamples $B$ in the set $\{1,\ldots,100\}$. We set the sample size $m = 100$ for the univariate case, and set the sample size $m = 1000$ for the multivariate case (left plot of figure~\ref{fig:multivariate_subsample}). The size of each subsample is set at $k = m^{2/3}$, and we compute the coverage and mean width of the confidence interval over $1000$ replications. The rate of convergence in setting~$2$ is $\tau_m= \sqrt{m}$. The setting~$3$ is of particular interest as it is known (\cite{loh1984estimating}, \cite{politis1994large}) that although classical bootstrap is invalid in this case, subsampling works here under weaker conditions with $\tau_m = m$. We can see from both figures~\ref{fig:multivariate_subsample} and \ref{fig:univariate_subsample} that the usual subsampling confidence interval attains the coverage of $0.9$ as $B$ increases, but does not have valid coverage for small $B$. The modified subsampling confidence interval attains the coverage of $0.9$ for $B \geq (2/\alpha) - 1 = 19$ in both figures~\ref{fig:multivariate_subsample} and \ref{fig:univariate_subsample}. As in \Cref{subsec:boot}, the simulation study in setting~2 goes on to illustrate that the modified subsampling confidence interval $\mathrm{CI}^{\mathtt{mod-subsample}}_{m,k,B,\alpha}$ provides asymptotic validity for any $B \geq (2/\alpha) -1$ even in high-dimensional regimes as long as we are dealing with a univariate projection of the test statistic and the standard subsampling consistency condition holds. We observe from the right plot in figure~\ref{fig:multivariate_subsample} that even when $m = 100$ (in setting~2), $\mathrm{CI}^{\mathtt{mod-subsample}}_{m,B,\alpha}$ guarantees valid coverage over the entire range of $\alpha \in [0.1, 0.5]$ for all three values of $B$. It can be seen from both the figures (and also from the definitions) that the modified subsampling confidence interval is slightly wider than the vanilla subsampling confidence interval.
\subsection{Inference using Stochastic Gradient Descent (SGD) estimator}
\label{subsec:SGD}
Stochastic Gradient Descent (SGD) is an important method for estimating target functionals from large datasets in a computationally fast and memory-efficient manner. The algorithm updates the estimate iteratively by computing the gradient of the objective function using one observation at a time (see \cite{wang2016statistical} for a brief review). Owing to this property, SGD has become especially popular in online settings where data are received sequentially. Although the asymptotic properties of SGD have been studied in earlier works (\cite{ruppert1988efficient}, \cite{polyak1992acceleration}), the problem of obtaining asymptotically valid confidence intervals for SGD has received less attention. Since standard inference methods like the usual non-parametric bootstrap does not apply here, prior works (\cite{fang2018online}) have studied versions of weighted bootstrap (\cite{rubin1981bayesian}) to construct an online bootstrap process to obtain asymptotically valid confidence interval. 

\begin{figure}[h]
    \centering
    \begin{minipage}[t]{0.48\textwidth}
        \centering
        \includegraphics[width=\linewidth]{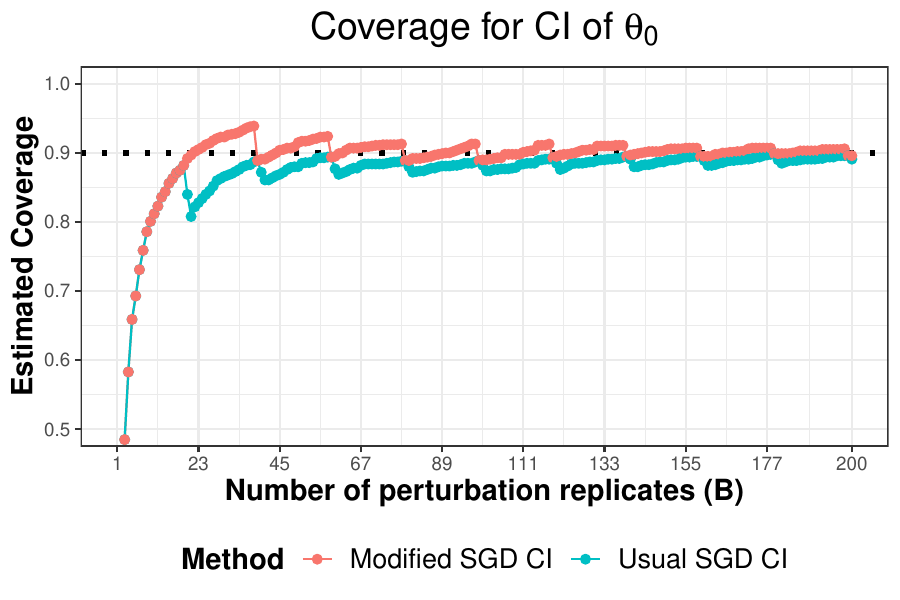}
    \end{minipage}
    \hfill
    \begin{minipage}[t]{0.48\textwidth}
        \centering
        \includegraphics[width=\linewidth]{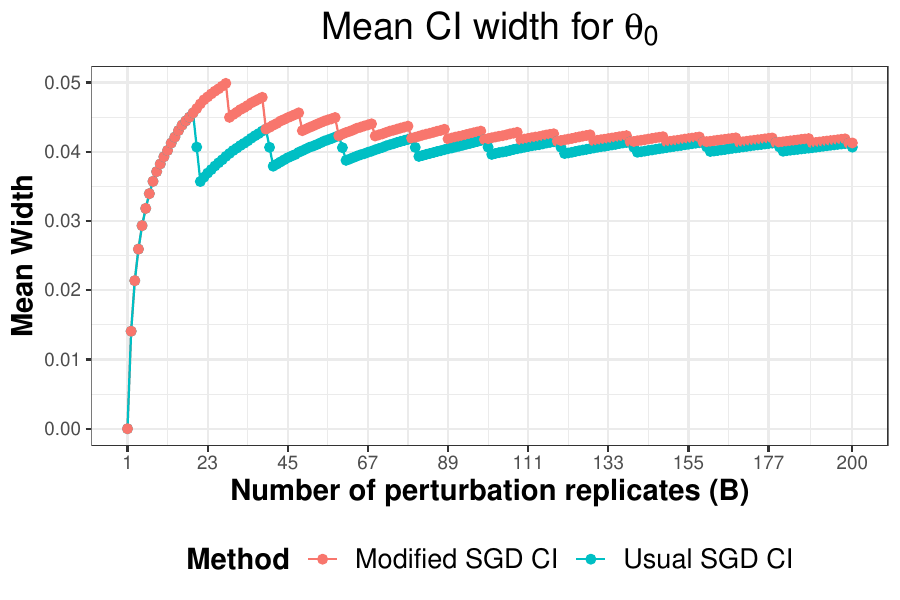}
    \end{minipage}
    \caption{Comparison of the coverage and the mean width between {usual} and modified SGD confidence interval in Setting~4 at level $1 - \alpha = 0.9$. The dotted horizontal line represents the nominal level of $0.9$. }
    \label{fig:sgd}
\end{figure}

In this subsection, we discuss the confidence interval proposed in \cite{fang2018online} and we propose a simple modification that yields confidence intervals with valid coverage without the need for the number of bootstrap samples to be too large. Consider the problem of obtaining an asymptotically valid confidence interval for $\theta_0 = \mbox{arg min}_{\theta \in \Theta} \mathbb E[\ell(\theta, Z)]$ $ \in \mathbb R^p$. SGD starts from an initial estimate $\widehat \theta_0$ and recursively updates the estimate using the point $Z_n$ obtained at the $n$-th stage ($n \in [m]$),
\[
\widehat \theta_n = \widehat \theta_{n-1} - \gamma_n \nabla \ell(\widehat \theta_{n-1},Z_n), \quad \overline \theta_n = \frac{1}{n}\sum_{i = 1}^n \widehat \theta_i,
\]
where the learning rate is generally taken as $\gamma_n = \gamma_1 n^{-\tau}$ for some constant $\gamma_1 > 0$ and $\tau \in (0.5, 1)$. The mean estimate $\overline \theta_n$ is generally considered for greater stability (\cite{ruppert1988efficient}, \cite{polyak1992acceleration}). For each $b \in [B]$, we generate a set $\{W_{i,b}: i \in [m]\}$ of iid non-negative random variables with both mean and variance equal to $1$ and construct perturbed SGD estimates as,
\[
\widehat \theta_n^{*b} = \widehat \theta_{n-1}^{*b} - \gamma_n W_{n,b}\nabla \ell(\widehat \theta_{n-1}^{*b},Z_n), \quad \overline \theta_n^{*b} = \frac{1}{n}\sum_{i = 1}^n \widehat \theta_i^{*b}.
\]
With some abuse of notation let $(\overline \theta^*_{(1)},\ldots, \overline \theta^*_{(B)} )$ be the order statistics of the weighted bootstrap estimates $\{\overline \theta_m^{*b}\}_{b = 1}^B$. The random weighting (RW-Q) method discussed in \cite{fang2018online} considers the confidence interval $\mathrm{CI}^{\mathtt{vanilla-SGD}}_{m,B,\alpha} = (2 \overline \theta_m - \overline \theta^*_{(\lceil B(1 - (\alpha/2))\rceil)}, 2 \overline \theta_m - \overline \theta^*_{(\lceil B(\alpha/2)\rceil)}]$. It can be checked that $\mathrm{CI}^{\mathtt{vanilla-SGD}}_{m,B,\alpha}$ is an asymptotically valid $(1-\alpha)$ confidence interval of $\theta_0$ provided $B \rightarrow \infty$ and the bootstrap consistency condition $\sup_{x \in \mathbb R} |\mathbb P(\overline \theta_m^{*b} - \overline \theta_m \leq x | D_m) - \mathbb P(\overline \theta_m - \theta_0 \leq x) | \stackrel{P}{\rightarrow} 0$ where $D_m = \{Z_1,\ldots,Z_m\}$. \cite{fang2018online} discusses several assumptions under which the bootstrap consistency condition holds (Theorem $3$ in their paper). We introduce the modified confidence interval $\mathrm{CI}^{\mathtt{mod-SGD}}_{m,B,\alpha} =  (2 \overline \theta_m - \overline \theta^*_{(\lceil (B + 1)(1 - \alpha)\rceil + \lfloor (B+1)(\alpha/2)\rfloor )}, 2 \overline \theta_m - \overline \theta^*_{(\lfloor (B+1)(\alpha/2)\rfloor)}]$ to remove the requirement of the number of bootstrap samples to go to infinity for the confidence interval to be asymptotically valid. The following corollary quantifies the slack in coverage of the modified confidence interval $\mathrm{CI}^{\mathtt{mod-SGD}}_{m,B,\alpha}$ as a function of $m, B$.
\begin{corollary}
    \label{cor:sgd}
    Using the stated notations, the modified SGD confidence interval $\mathrm{CI}^{\mathtt{mod-SGD}}_{m,B,\alpha}$ satisfies the following coverage guarantee,
     \[
    \left|P( \theta_0 \in \mathrm{CI}^{\mathtt{mod-SGD}}_{m,B,\alpha} )  - \frac{\lceil(B+1)(1-\alpha) \rceil}{B+1} \right| \leq \widetilde d_{\mathrm{KS}}(F_0(D_m), U(0,1)),
    \]
    where $F_0(D_m) = \mathbb P(\overline \theta_m^{*b}  \leq 2\overline \theta_m - \theta_0 | D_m)$. The bounds on coverage can be restated as,
    \[
    - \widetilde d_{\mathrm{KS}}(F_0(D_m), U(0,1)) \leq \mathbb P( \theta_0 \in \mathrm{CI}^{\mathtt{mod-SGD}}_{m,B,\alpha} ) - (1 - \alpha) \leq \frac{1}{B+1} +  \widetilde d_{\mathrm{KS}}(F_0(D_m), U(0,1)).
    \]

\end{corollary}
Since the slack term $\widetilde d_{\mathrm{KS}}(F_0(D_m), U(0,1))$ does not have any dependence on the number of bootstrap samples, the asymptotic validity of $\mathrm{CI}^{\mathtt{mod-SGD}}_{m,B,\alpha}$ holds for any $B \geq (2/\alpha) - 1$ (for a non-trivial confidence interval) if $\widetilde d_{\mathrm{KS}}(F_0(D_m), U(0,1)) \rightarrow 0$ as $m \rightarrow \infty$. The first statement of \Cref{cor:sgd} also implies that if $(B+1)(1-\alpha) \in \mathbb N$ is an integer, then we have exact coverage of $(1-\alpha)$ when $\widetilde d_{\mathrm{KS}}(F_0(D_m), U(0,1)) \rightarrow 0$ as $m \rightarrow \infty$. The proof of \Cref{cor:sgd} is analogous to that of \Cref{cor:boot} and is therefore omitted for brevity. In light of Proposition~\ref{prop:ks_connection}, if the distribution of $(\overline \theta_m - \theta_0)$ is continuous, one sufficient condition for $\widetilde d_{\mathrm{KS}}(F_0(D_m), U(0,1)) \rightarrow 0$ as $m \to \infty$ is the standard bootstrap consistency assumption for SGD,
\[
\Delta^*(D_m) = \sup_{x \in \mathbb R} |\mathbb P(\overline \theta_m^{*b} - \overline \theta_m \leq x | D_m) - \mathbb P(\overline \theta_m - \theta_0 \leq x) | \stackrel{P}{\rightarrow} 0, \quad \mbox{as} \quad m \rightarrow \infty. 
\]
We compare the performance of the modified SGD confidence interval $\mathrm{CI}^{\mathtt{mod-SGD}}_{m,B,\alpha}$ with that of the vanilla SGD confidence interval $\mathrm{CI}^{\mathtt{vanilla-SGD}}_{m,B,\alpha}$ for computing the confidence interval of $\theta_0$ (the first coordinate of $\theta$) in the following quantile regression problem.
\paragraph{Setting $4$}:\label{para:setting:quantile_regression} $Z_i = (X_i, Y_i)$ is the data-point observed at $i$-th time point for $i\in [N]$. The data points are iid with $X_{i,j} \stackrel{iid}{\sim}  N(0,1)$ for $j \in \{1,2,3\}$ and $\epsilon_i = Y_i - \theta^\top X_i \sim \mbox{Double Exponential}(0,1) $, $\theta = (\theta_0, \theta_1, \theta_2) = \mbox{arg min}_{\theta} \mathbb E[\rho(Y - X^\top \theta)] = (0.2, -0.2, 0)$ where $\rho(u) = u((1/2) - \textbf{1}\{u < 0\})$, the SGD and the perturbed SGD iterates are computed as follows for $n \in [N]$,
\begin{equation*}
    \begin{split}
        \widehat \theta_n =& \widehat \theta_{n-1} + \gamma_n \{ (1/2) - \textbf{1}\{Y_n - X_n^\top \widehat \theta_{n-1} <0 \} \}X_n, \\
        \widehat \theta_n^{*b} =& \widehat \theta_{n-1}^{*b} + \gamma_n \{ (1/2) - \textbf{1}\{Y_n - X_n^\top \widehat \theta_{n-1}^{*b} <0 \} \}X_n. 
    \end{split}
\end{equation*}
The mean iterates $\overline \theta_n$ and $\{\overline \theta_m^{*b}\}_{b = 1}^B$ are computed as described earlier in this subsection. We set the learning rate $\gamma_n = n^{-2/3}$, $\mbox{burn-in} = 2000$, total number of SGD iterates $N = 10000$, and $m =N - \mbox{burn-in}  = 8000$. We compute the coverage and mean width of the confidence interval over $1000$ replications. We observe from figure~\ref{fig:sgd} that the RW-Q method (\cite{fang2018online}) attains the coverage of $0.9$ as $B$ increases, but does not have valid coverage for small $B$. On the other hand, the modified confidence interval $\mathrm{CI}^{\mathtt{mod-SGD}}_{m,B,\alpha}$ based on the SGD estimator attains the coverage of $0.9$ for $B \geq (2/\alpha) - 1 = 19$ in figure~\ref{fig:sgd}. It can be seen from the figure (and also from the definitions) that the modified SGD-based confidence interval is slightly wider than the RW-Q confidence interval.
\subsection{Permutation testing}
\label{subsec:permute}
Permutation tests trace back to Fisher's randomization based approach to inference and are now a standard tool for distribution-free hypothesis testing under invariance or exchangeability assumptions (\cite{fisher1935design,lehmann2005testing}). In permutation testing approach, one observes data $D_m = \{X_1,\dots,X_m\}$ and tests the given null hypothesis by comparing a test statistic $T(X)$ (where $X = (X_1,\ldots,X_m)$) to the same statistic evaluated on permuted versions of the data. The permutation test works when the statistics evaluated on all permuted versions of the data are exchangeable. As emphasized by \cite{ramdas2023permutation}, this basic idea remains one of the most widely used nonparametric inferential tools because it avoids the need to specify the null distribution of $T(X)$ exactly, while still yielding valid finite-sample $p$-values under exchangeability.

A central issue of permutation testing is the computational complexity. Evaluating $T(X_\sigma)$ over all $m!$ permutations quickly becomes infeasible, and this motivates restricting attention to a smaller collection of permutations. Classical theory therefore often works with a subgroup $G \subseteq S_m$ ($S_m$ is the full permutation group of cardinality $m!$), chosen to balance computational tractability and statistical power, and then either enumerates all elements of $G$ or samples from $G$ uniformly at random (\cite{hemerik2018exact,koning2022faster}). A canonical example where permutation testing is used often is for testing independence in IID pairs $(X_i,Y_i)$. Under the null hypothesis $X \perp\!\!\!\perp Y$, the coordinates of $X=(X_1,\dots,X_m)$ are exchangeable conditional on $Y=(Y_1,\dots,Y_m)$ and hence one may compare an observed association statistic such as
\[
T(X,Y)=\left|\mathrm{Corr}\big((X_1,\dots,X_m),(Y_1,\dots,Y_m)\big)\right|
\]
to its values on permuted data,
\[
T(X_\sigma,Y)=\left|\mathrm{Corr}\big((X_{\sigma(1)},\dots,X_{\sigma(m)}),(Y_1,\dots,Y_m)\big)\right|.
\]
The key contribution of \cite{ramdas2023permutation} is to show that the usual emphasis on subgroup structure and uniform sampling can be relaxed by introducing an auxiliary random permutation. Their generalized construction draws $\sigma_0 \sim q$, where $q$ is an arbitrary distribution on the set of all permutations $S_m$, and then forms a $p$-value by comparing $T(X)$ to statistics of the form $T(X_{\sigma \circ \sigma_0^{-1}})$,
\begin{equation}
    \label{eq:ramdas_P}
    P_{\mbox{gen}} = \sum_{\sigma \in S_n} q (\sigma) \textbf{1} \left\{T(X_{\sigma \circ \sigma_0^{-1}}) \geq T(X) \right\}.
\end{equation}
In the special case where $q$ is uniform on a subgroup $G$, the closure of $G$ under inversion and composition implies that $\{\sigma \circ \sigma_0^{-1} : \sigma \in G\} = G$ for every $\sigma_0 \in G$. Therefore the generalized $p$-value reduces exactly to the classical subgroup-based permutation $p$-value. Thus the subgroup method appears as a special case of the broader framework. This shows that the subgroup structure is sufficient, but not fundamentally necessary, for valid permutation-based inference (\cite{hemerik2018exact, ramdas2023permutation}). 

We can recover the results of \cite{ramdas2023permutation} using the master theorems proposed in \Cref{sec:master_thms}. Consider an arbitrary subset $G \subseteq S_m$. For any $\sigma \in G$, we denote the permuted vector $(X_{\sigma(1)}, \ldots, X_{\sigma(m)})$ by $X_{\sigma}$. We draw $B$ permutations $\{\sigma^*_1,\ldots,\sigma^*_B\}$ uniformly (with replacement) from $G$ where $B \leq |G|$. With some abuse of notation, we denote the order statistics of $\{T(X_{\sigma^*_1}), \ldots, T(X_{\sigma^*_B})\}$ by $\{T^*_{(1)},\ldots,T^*_{(B)}\}$. The rejection rule is $\phi^{\mathtt{mod-permute}}_{m,B,\alpha}(X) = \textbf{1}\{T(X) \geq T^*_{(\lceil B(1-\alpha) \rceil +2} \}$ if $B  = |G|$ and $\phi^{\mathtt{mod-permute}}_{m,B,\alpha}(X) = \textbf{1}\{T(X) \geq T^*_{(\lceil (B+1)(1-\alpha) \rceil +1} \}$ if $B  < |G|$. The next corollary shows that the test $\phi^{\mathtt{mod-permute}}_{m,B,\alpha}(X)$ enjoys valid type-I error control. 
\begin{corollary}
    \label{cor:permute}
    We assume that under the null hypotheses $H_0$, the statistics based on permuted data $\{T(X_{\sigma})\}_{\sigma \in G}$ are exchangeable. We also assume that the statistics $\{T(X_{\sigma})\}_{\sigma \in G}$ do not have any ties almost surely i.e.\ $\mathbb P(T(X_{\sigma}) = T(X_{\sigma'}) \mbox{ for any } \sigma \neq \sigma', \sigma \in G, \sigma' \in G) = 0$.
    Using the stated notations, the modified permutation test $\phi^{\mathtt{mod-permute}}_{m,B,\alpha}(X)$ satisfies the following calibration guarantee for any subset $G \subseteq S_m$ and for any $B \leq |G|$,
    \[
    \alpha -\frac{2 + \alpha \textbf{1}\{B = |G|\}}{B+1} -\frac{1}{|G|} \leq \mathbb E_{H_0}[\phi^{\mathtt{mod-permute}}_{m,B,\alpha}(X)] \leq \alpha .
    \]
\end{corollary}
\Cref{cor:permute} guarantees the calibration of $\phi^{\mathtt{mod-permute}}_{m,B,\alpha}(X)$ for any $B \leq |G|$. The proof of the corollary is provided in \Cref{app:proof_permute}. The key step in the proof is to show that the Kolmogorov Smirnov distance $d_{\mathrm{KS}}(F_0(D_m), U(0,1))  \leq 1 / |G|$ is bounded by the reciprocal of cardinality of the subset, where $F_0(D_m) = \mathbb P( T(X_{\sigma_b^*}) \leq T(X) | D_m)$. \Cref{cor:permute} echoes the message of \cite{ramdas2023permutation} that we do not need any subgroup structure on $G$ to ensure the validity of the permutation test as long as we pick the permutations randomly from $G$. 

\subsection{Conformal prediction}
\label{subsec:CP}
In conformal prediction (\cite{vovk2005algorithmic}, \cite{lei2014distribution}, \cite{lei2018distribution}), our objective is to obtain a valid $(1-\alpha)$ prediction set of $Y_{m+1}$ given a calibration data $D_{\mbox{calib}} = \{(X_1,Y_1),\ldots, (X_m, Y_{m})\}$ and a test point $X_{m+1}$ where the $(m+1)$ points $D_{\mbox{calib}}  \cup(X_{m+1}, Y_{m+1})$ are assumed to be exchangeable. Given a pre-trained non-conformity score function $R(\cdot, \cdot)$ that maps the data points to non-negative real values, the split conformal prediction approach (\cite{papadopoulos2002inductive}, \cite{lei2014distribution}) computes the non-conformity scores $R_i = R(X_i, Y_i)$ for all $i \in [m]$ and reports the prediction set $\widehat C_{m,\alpha}(X_{m+1}) = \{y: R(X_{m+1}, y) \leq R_{(\lceil (m+1)(1-\alpha) \rceil)}\}$ for $Y_{m+1}$ (here $\{R_{(1)}, \ldots, R_{(m)} \}$ are the order statistics of $\{R_1,\ldots, R_m \}$). It can be shown using suitable exchangeability arguments that the prediction set $\widehat C_{m,\alpha}(X_{m+1})$ enjoys valid coverage,
\[
\mathbb P( Y_{m+1} \in \widehat C_{m,\alpha}(X_{m+1})) \geq 1- \alpha. 
\]
The above coverage guarantee no longer holds if the data points in $D_{\mbox{calib}}$ and the test point $(X_{m+1}, Y_{m+1})$ are not exchangeable. Let $S = (R_1,\ldots,R_{m+1})$ (where $R_{m+1} = R(X_{m+1}, Y_{m+1})$) and for any $i \in [m]$ let $S^i = (R_1,\ldots,R_{i-1},R_{m+1}, R_{i+1},\ldots,R_m,R_i)$ be the sequence obtained after swapping $R_{m+1}$ with $R_i$ in the sequence $S$. Theorem $2$ and Theorem $3$ of \cite{barber2023conformal} characterize the slack in coverage of the usual prediction set $\widehat C_{m,\alpha}(X_{m+1})$ when exchangeability of data points is violated. In particular, the following coverage guarantees hold if the distributions of the scores $R_1,\ldots,R_{m+1}$ are continuous (continuity of distributions is needed only for the upper bound on coverage), 
\begin{equation}
    \label{eq:barber}
    -\frac{1}{m+1} \sum_{i = 1}^m d_{\mathrm{TV}}(S,S^i) \leq \mathbb P( Y_{m+1} \in \widehat C_{m,\alpha}(X_{m+1})) -(1 - \alpha) \leq \frac{1}{m+1} + \frac{1}{m+1} \sum_{i = 1}^m d_{\mathrm{TV}}(S,S^i). 
\end{equation}
In particular if the scores $R_1,\ldots,R_{m+1}$ are independent, then the coverage bounds can be simplified as,
\begin{equation}
    \label{eq:barber_ind}
    -\frac{2}{m+1} \sum_{i = 1}^m d_{\mathrm{TV}}(R_i,R_{m+1}) \leq \mathbb P( Y_{m+1} \in \widehat C_{m,\alpha}(X_{m+1})) -(1 - \alpha) \leq \frac{1}{m+1} + \frac{2}{m+1} \sum_{i = 1}^m d_{\mathrm{TV}}(R_i, R_{m+1}). 
\end{equation}
The following corollary of the master theorems in \Cref{sec:master_thms} provides an alternate slack in coverage in terms of the Kolmogorov-Smirnov distances between $R_i$'s ($i \in [m]$) and $R_{m+1}$ when the scores $R_1,\ldots,R_{m+1}$ are independent and possibly have different distributions. 
\begin{corollary}
    \label{cor:indep_CP}
Suppose the scores $R_1,\ldots,R_{m+1}$ are independent continuous random variables and $F_i(z) = \mathbb P(R_i < z)$ for $i \in [m+1]$ and for all $z \in \mathbb R$. Let $\kappa_i = d_{\mathrm{KS}}(F_i,F_{m+1})$ for $i \in [m]$ and $\overline F(z) = (1/m) \sum_{i = 1}^m F_i(z)$ for $z \in \mathbb R$. Then we have the following coverage guarantee,
\[
-\delta \leq \mathbb P( Y_{m+1} \in \widehat C_{m,\alpha}(X_{m+1})) -(1 - \alpha) \leq \frac{1}{m+1} + \delta \quad \mbox{where},
 \]
 \begin{equation*}
    \delta = d_{\mathrm{KS}}(\overline F(R_{m+1}), U(0,1)) +\left(\sum_{i = 1}^m \kappa_i^2\right)\left[I_m +  d_{\mathrm{KS}}(\overline F(R_{m+1}), U(0,1))\right].
\end{equation*}
Here $I_m$ is as defined in \Cref{thm:main_result_extension}. Moreover, the following lower bound on coverage holds for the modified prediction interval $\widehat C_{m,\alpha}^{\mathtt{mod-pred}}(X_{m+1}) = \{y: R(X_{m+1}, y) \leq R_{(m+1 - \lfloor 2m \alpha/3 \rfloor )}\}$,
\begin{equation*}
    \mathbb P( Y_{m+1} \in \widehat C_{m,\alpha}^{\mathtt{mod-pred}}(X_{m+1})) -(1 - \alpha)  \geq -3 d_{\mathrm{KS}}(\overline F(R_{m+1}), U(0,1)).
\end{equation*}
\end{corollary}
The proof of \Cref{cor:indep_CP} is provided in \Cref{app:proof_indep_CP}. \Cref{cor:indep_CP} shows that even when the distributions of the test point and the data points in the calibration set are non-identical, it is possible to characterize the loss in coverage of the usual and the modified prediction interval in terms of Kolmogorov-Smirnov distances between the distributions of data points in the calibration set and the distribution of the test point. Since Kolmogorov-Smirnov distances are often much smaller than the corresponding total variation distances, the bounds in \Cref{cor:indep_CP} are generally tighter than \eqref{eq:barber_ind}. 

Consider the following example where $R_i \sim \delta((i - (1/2))/m)$ for all $i \in [m]$ (a random variable distributed as $\delta(x)$ takes the value $x$ with probability $1$) and $R_{m+1} \sim \mbox{Uniform}(0,1)$. Our goal is to design a prediction interval of $Y_{m+1}$ that has coverage of at least $0.9$ asymptotically. Since the distribution of $R_i$ is discrete and the distribution of $R_{m+1}$ is continuous, we have $d_{\mathrm{TV}}(R_i, R_{m+1}) = 1$ for all $i \in [m]$. Therefore the slack in coverage in \eqref{eq:barber_ind} converges to $2$,
\[
\frac{2}{m+1} \sum_{i = 1}^m d_{\mathrm{TV}}(R_i, R_{m+1}) = \frac{2m}{m+1} \rightarrow 2 \quad \mbox{as} \quad m \rightarrow \infty. 
\]
This implies that it is not possible to show asymptotic coverage guarantee of $0.9$ of the prediction interval $\widehat C_{m,\alpha}(X_{m+1})$ for any $\alpha \in (0,1)$ using the bounds in \eqref{eq:barber_ind}. However, the tighter Kolmogorov-Smirnov bounds in \Cref{cor:indep_CP} allow us to design a prediction interval that is guaranteed to have an asymptotic coverage of $0.9$. In this example, we have $F_i (z) = \textbf{1}\{z > (i - (1/2))/m\}$ and $\overline F(z) = (1/m) | \{ i: (i - (1/2))/m < z \}|$ for $z \in \mathbb R$. It can be shown that distribution of $\overline F(R_{m+1}) $ is $H_m$ where,
\[
H_m(t) = \begin{cases}
    0 \quad & \mbox{if} \quad t <0, \\
    \frac{\lfloor mt \rfloor}{m} + \frac{1}{2m} \quad & \mbox{if} \quad t \in [0,1), \\
    1 \quad & \mbox{if} \quad t \geq 1. 
\end{cases}
\]
The Kolmogorov-Smirnov distance appearing in the slack in the second statement of \Cref{cor:indep_CP} is $3 d_{\mathrm{KS}}(\overline F(R_{m+1}), U(0,1)) = 3 \sup_{t \in \mathbb R}|H_m(t) - t| = 3/(2m)$. Therefore the modified prediction interval $\widehat C_{m,0.1}^{\mathtt{mod-pred}}(X_{m+1}) = \{y: R(X_{m+1}, y) \leq R_{(m+1 - \lfloor 0.2m /3 \rfloor )}\}$ has the required asymptotic coverage guarantee,
\[
\mathbb P(Y_{m+1} \in \widehat C_{m,0.1}^{\mathtt{mod-pred}}(X_{m+1})) \geq 0.9 - \frac{3}{2m}. 
\]
The bounds in \eqref{eq:barber_ind} can be simplified if all the scores $\{R_1,\ldots,R_m\}$ in the calibration set are identically distributed. The distribution of the test score $R_{m+1}$ may be different from the common distribution of the scores in the calibration set.
\begin{equation}
    \label{eq:barber_ind_common}
    -\frac{2m}{m+1} d_{\mathrm{TV}}(R_1,R_{m+1}) \leq \mathbb P( Y_{m+1} \in \widehat C_{m,\alpha}(X_{m+1})) -(1 - \alpha) \leq \frac{1}{m+1} + \frac{2m}{m+1} d_{\mathrm{TV}}(R_1, R_{m+1}). 
\end{equation}
An application of the results in \Cref{sec:master_thms} provides a much sharper coverage guarantee in this scenario. The next corollary shows that the slack in coverage of the prediction interval $ \widehat C_{m,\alpha}(X_{m+1})$ can be improved to $d_{\mathrm{KS}}(R_1, R_{m+1})$ when the scores in the calibration set are identically distributed. 
\begin{corollary}
    \label{cor:indep_CP_same}
Suppose the scores $R_1,\ldots,R_{m+1}$ are independent continuous random variables and the scores $R_1,\ldots,R_m$ are identically distributed. Then we have the following coverage guarantee,
\[
-d_{\mathrm{KS}}(\widetilde F_0(R_{m+1}), U(0,1)) \leq \mathbb P( Y_{m+1} \in \widehat C_{m,\alpha}(X_{m+1})) -(1 - \alpha) \leq \frac{1}{m+1} +  d_{\mathrm{KS}}(\widetilde F_0(R_{m+1}), U(0,1)). 
\]
where $\widetilde F_0(z) = \mathbb P(R_1 < z)$ for all $z \in \mathbb R$. 
\end{corollary}
The proof of \Cref{cor:indep_CP_same} is provided in \Cref{app:proof_indep_CP_same}. If the distribution of $R_1$ is continuous, then the slack in coverage simplifies as $ d_{\mathrm{KS}}(\widetilde F_0(R_{m+1}), U(0,1)) = d_{\mathrm{KS}}(F_0(R_{m+1}), U(0,1)) = d_{\mathrm{KS}}(R_1, R_{m+1}) < d_{\mathrm{TV}}(R_1, R_{m+1})$ where $ F_0(z) = \mathbb P(R_1 \leq z)$ (for $z \in \mathbb R$). Therefore we can establish that the prediction set $\widehat C_{m,\alpha}(X_{m+1})$ retains stronger coverage guarantees even under distribution shift in the test data. We can improve the slack in coverage in \eqref{eq:barber} even for a general vector of scores $(R_1,\ldots,R_{m+1})$ where $\{R_i\}_{i = 1}^{m+1}$ are possibly not independent. 
\begin{corollary}
    \label{cor:exchangeable_CP}
The following coverage guarantee is valid for any $(m+1)$ dimensional vector of scores $(R_1,\ldots,R_m,R_{m+1})$ and for any $0 < \alpha< 1$,
\begin{equation*}
  \mathbb P( Y_{m+1} \in \widehat C_{m,\alpha}(X_{m+1})) - (1-\alpha) \geq  -\delta,
\end{equation*}
where $ \delta = \Gamma(R_1,\ldots,R_m,R_{m+1})$ (see \eqref{eq:gamma_defn}). Additionally if the distributions of the random variables $\{R_i\}_{i = 1}^{m+1}$ are continuous then we have,
\[
\mathbb P( Y_{m+1} \in \widehat C_{m,\alpha}(X_{m+1})) - (1-\alpha)  \leq  \delta +  \frac{2}{m+1}.
\]
\end{corollary}
\Cref{cor:exchangeable_CP} quantifies the slack in coverage of the classical split conformal prediction interval when applied to a general score vector $(R_1,\ldots,R_m,R_{m+1})$ that may not be exchangeable or independent. The proof of \Cref{cor:exchangeable_CP} follows from \Cref{thm:beyond_exchangeability} by observing that the $(m+1)$ dimensional vector $(R_1,\ldots,R_m,R_{m+1})$ plays the role of the $(B+1)$ dimensional vector $(W_1,\ldots,W_B,\psi(Z))$. The proof is skipped for brevity.
\begin{remark}[Sharper coverage bounds for general non-exchangeable score vector]
    \label{rem:sharper_exchangeability}
\Cref{cor:exchangeable_CP} provides sharper coverage bounds than that provided in \eqref{eq:barber} (\cite{barber2023conformal}) i.e.\ we have,
\[
\Gamma(R_1,\ldots,R_m,R_{m+1}) \leq \frac{1}{m+1}\sum_{i = 1}^m d_{\mathrm{TV}}(S,S^i). 
\]
The proof of the above inequality is provided in \Cref{app:proof_rem_exchangeavility}.
\end{remark}

\subsection{Randomization tests without Group Invariance}
\label{subsec:group_invariance}
Consider $m$ IID data points $X_1,\ldots,X_m \stackrel{iid}{\sim} P$. Let $X = (X_1,\ldots,X_m)$ be the $m$ dimensional vector of the data points. We wish to test the null hypothesis $H_0: P \in \mathcal P_0$. Let $G$ be a finite group of transformations of the sample points such that for any $g \in G$ we have $gX \stackrel{d}{=} X$ whenever $ P \in \mathcal P_0$. In other words, under the null hypothesis the distribution of the sample $X$ is invariant under transformations from the group $G$. Let $T(X)$ be any given statistic for testing $H_0$ where $X = (X_1,\ldots,X_m)$. With slight abuse of notation, we let $\{T_{(1)},\ldots,T_{(|G|)}\}$ be the order statistics of $\{T(gX) \}_{g \in G}$. Consider the vanilla test $\phi^{\mathtt{vanilla}}(X) = \textbf{1}\{ T(X) > T_{(|G| - \lfloor |G| \alpha \rfloor) }\}$. Under the group invariance condition, it can be shown that (\cite{romano1990behavior}) the test $\phi^{\mathtt{vanilla}}(X)$ has valid type-I error control,
\[
\mathbb E_{H_0} [\phi^{\mathtt{vanilla}}(X)] \stackrel{(i)}{=} \frac{1}{|G|} \sum_{g \in G} \mathbb E_{H_0} [\phi^{\mathtt{vanilla}}(gX)] = \frac{\lfloor |G| \alpha \rfloor}{|G|} \leq \alpha. 
\]
The equality $(i)$ holds because $G$ is a finite group. To illustrate this setting, let $X_1,\ldots,X_m \stackrel{\mathrm{iid}}{\sim} P_\theta$ with \(P_\theta\) symmetric around \(\theta\), and consider testing the null hypothesis $H_0:\theta=\theta_0$. Define \(G\) to be the finite sign-flip group on \(\mathbb R^m\), namely the set of transformations
\[
gx=\bigl((-1)^{i_1}x_1,\ldots,(-1)^{i_m}x_m\bigr),
\]
with \(i_j\in\{0,1\}\) for each \(j\in[m]\). This group has cardinality \(2^m\). Under \(H_0\), the centered vector \((X_1-\theta_0,\ldots,X_m-\theta_0)\) is invariant in distribution under every transformation \(g\in G\). Hence \(\phi^{\mathtt{vanilla}}_{\theta_0}(X)\), constructed from the statistic $T_{\theta_0}(X)=m^{-1/2}\sum_{i=1}^m (X_i-\theta_0)$, yields a finite-sample level-\(\alpha\) test of \(H_0\).

Having considered this example, we now turn back to the general framework. Structural assumptions such as symmetry of the underlying distribution may be too restrictive in many practical settings. When the group invariance condition fails under the null hypothesis, the test \(\phi^{\mathtt{vanilla}}(X)\) is no longer finite-sample valid. \cite{romano1990behavior} investigates the asymptotic validity of randomization tests in precisely such settings, where the group invariance assumption is violated. When the size of the group $G$ is too large, one may sample transformations $g_1^*,\ldots,g_B^*$ uniformly (with replacement) from the group $G$ and consider the randomization test based on these $B$ transformations. With slight abuse of notation, we let $(T_{(1)},\ldots,T_{(B)}\}$ be the order statistics of $\{T(g^*_bX) \}_{b \in [B]}$. The next result establishes the asymptotic validity of the modified test $\phi^{\mathtt{mod}}(X) = \textbf{1}\{ T(X) \geq T_{( \lceil (B+1)(1 -\alpha) \rceil) }\}$ when the group invariance condition fails to hold in finite samples and the number of Monte Carlo resamples $B \geq (1/\alpha) - 1$. 

\begin{corollary}
    \label{cor:randomized_test}
Using the stated notations, the modified testing rule $\phi^{\mathtt{mod}}(X) $ satisfies the following guarantee,
\[
\left|\mathbb E_{H_0}[\phi^{\mathtt{mod}}(X)]  - \frac{\lfloor (B+1)\alpha \rfloor}{B+1}\right| \leq d_{\mathrm{KS}}(F_0(D_m),U(0,1)),
\]
where $D_m= \{X_1,\ldots,X_m\}$ and $F_0(D_m) = \mathbb P_{H_0}(T(g^*_bX)  \leq T(X)|D_m)$. The type-I error bounds can be restated as,
\[
\alpha - \frac{1}{B+1} -d_{\mathrm{KS}}(F_0(D_m),U(0,1))\leq \mathbb E_{H_0}[\phi^{\mathtt{mod}}(X)] \leq \alpha + d_{\mathrm{KS}}(F_0(D_m),U(0,1)).
\]
\end{corollary}
\Cref{cor:randomized_test} implies that the asymptotic type-I error of the test is bounded above by $\alpha$ for any $B \geq (1/\alpha) - 1$ if $d_{\mathrm{KS}}(F_0(D_m),U(0,1)) \rightarrow 0$ as $m \rightarrow \infty$. The first statement of \Cref{cor:randomized_test} also implies that if $(B+1)\alpha \in \mathbb N$, then the test $\phi^{\mathtt{mod}}(X) $ has exact asymptotic type-I error of $(1-\alpha)$ when $d_{\mathrm{KS}}(F_0(D_m),U(0,1)) \rightarrow 0$ as $m \rightarrow \infty$. The proof of \Cref{cor:randomized_test} is provided in \Cref{app:proof_rand_test_invariance}. \Cref{cor:randomized_test} conveys the same message as in \cite{romano1989bootstrap, romano1990behavior} that it is possible to guarantee asymptotic validity of the randomization test for finitely many Monte Carlo resamples $B$ even when the exact group invariance condition does not hold. In light of proposition~\ref{prop:ks_connection}, if the distribution of $T(X)$ is continuous, one sufficient condition for $d_{\mathrm{KS}}(F_0(D_m),U(0,1)) \rightarrow 0$ as $m \rightarrow \infty$ is the following consistency condition of the randomization distribution under the null hypothesis,
\[
\Delta^*(D_m) = \sup_{x \in \mathbb R} |\mathbb P_{H_0}(T(g^*_b X) \leq x |D_m) - \mathbb P_{H_0}(T(X) \leq x)  | \stackrel{P}{\rightarrow} 0, \quad \mbox{as} \quad m \rightarrow \infty. 
\]
The above consistency condition has been shown to hold in \cite{romano1990behavior} for some one-sample problems, such as testing for mean or median, even when the underlying distribution does not satisfy the group invariance condition. 
\section{Discussion}
\label{sec:extension}

This paper develops a general probability inequality for order statistics and shows that it serves as a common basis for a wide range of inference procedures built on resampling, randomization, and ranking. The main conclusion is that validity of these methods need not depend on the number of Monte Carlo replicates diverging to infinity. As long as the gap between the resampling (or randomization) distribution and the target law vanishes, the order-statistic arguments developed in this paper yield valid confidence intervals and tests, either in finite samples or asymptotically, even when the number of Monte Carlo draws is fixed. This viewpoint unifies a number of apparently distinct methods (including bootstrap, subsampling, conformal prediction, testing under permutation or approximate group invariance) through a single probabilistic mechanism. Our work also clarifies that while increasing the number of Monte Carlo replicates may improve stability or reduce variability, it is not inherently required for ensuring validity. In this sense, the present results may be viewed as extending and systematizing earlier insights, particularly those of \cite{hall1986number}, in a considerably broader framework.

There are several directions for further work. First, it would be useful to sharpen the bounds in the general master theorems, especially in the conditionally independent but non-identically distributed setting, where the current results are likely not optimal. Second, it would be of interest to move beyond order-statistic intervals and study probability inequalities of the form
\[
\mathbb{P}\bigl(W \in [g_\ell(W_1,\ldots,W_B),\, g_u(W_1,\ldots,W_B)]\bigr) \geq h_{g_\ell, g_u}(B) - \Delta(W_1,\cdots,W_B,W),
\]
for general measurable functions $g_\ell(\cdot)$ and $g_u(\cdot)$ when $(W_1,\ldots,W_B,W)$ has an arbitrary dependence structure. Here $h_{g_\ell,g_u}(B)$ denotes the ideal lower bound under an IID or exchangeable benchmark, while $\Delta(W_1,\ldots,W_B,W)$ measures the extent to which the joint law of $(W_1,\ldots,W_B,W)$ departs from that benchmark. Such an extension could substantially enlarge the scope of the theory and lead to new methods of inference. Third, it would be worthwhile to investigate applications of the master theorems to other resampling-based procedures not considered here. More broadly, understanding the tradeoff among validity, efficiency, and computational cost under fixed Monte Carlo budgets remains an important problem for future research.

\section*{Acknowledgments}
We thank Woonyoung Chang for a careful reading of the manuscript and for suggesting some corrections.

\bibliography{references}
\bibliographystyle{plainnat}
\newpage
\appendix
\setcounter{section}{0}
\setcounter{equation}{0}
\setcounter{figure}{0}
\setcounter{remark}{0}
\renewcommand{\thesection}{S.\arabic{section}}
\renewcommand{\theequation}{E.\arabic{equation}}
\renewcommand{\thefigure}{A.\arabic{figure}}
\renewcommand{\theremark}{R.\arabic{remark}}
  \begin{center}
  \Large {\bf Appendix to ``On a Probability Inequality for Order Statistics with Applications to Bootstrap, Conformal Prediction, and more''}
  \end{center}
\section{Proof of Theorem~\ref{thm:main_result}}
\label{app:proof_main_res}
We wish to bound the probability:
\[
\mathbb{P}\left(\psi(Z) \in [W_{(a)}, W_{(B-b)}] \right), \quad \mbox{for } 0 \leq a < B-b \leq B. 
\]
Following \cite{hall1986number} we observe that
\[
\mathbb{P}\left(\psi(Z) \in [W_{(a)}, W_{(B-b)}] \right) = \mathbb{E}_Z \left[ \mathbb{P}\left( \psi(Z) \in [W_{(a)}, W_{(B-b)}] \mid Z \right) \right].
\]
We consider the following definition,
\[
S_B(Z) := \sum_{i=1}^B \mathbf{1}(W_i \leq \psi(Z)) \sim \text{Bin}(B, F_0(Z)) \quad \mbox{conditioned on }Z.
\]
We also note the following equivalence,
\begin{align}
\label{eq:equiv}
W_{(a)} \leq \psi(Z) \leq W_{(B-b)} \begin{cases}
   &\Rightarrow a \leq S_B(Z) \leq n - b, \\ 
   & \Leftarrow a \leq S_B(Z) \leq n-b-1.
\end{cases} 
\end{align}

Hence:
\begin{align*}
    \mathbb{P}(\psi(Z) \in [W_{(a)}, W_{(B-b)}]) = & \mathbb{E}_Z \left[ \mathbb{P}\left( \psi(Z) \in [W_{(a)}, W_{(B-b)}] \mid Z \right) \right]\\
    \geq & \mathbb{E}_Z \left[ \mathbb{P}\left( S_B(Z) \in [a, B-b-1] \mid Z \right) \right] \\
    = & \mathbb{E}_Z\left[ \mathbb{P}\left( a \leq \text{Bin}(B, F_0(Z)) \leq B - b -1 \right) \right] \\
    = & \mathbb{E}_Z\left[ \mathbb{P}\left(\text{Bin}(B, F_0(Z))\geq  a  \right) + \mathbb{P}\left(\text{Bin}(B, F_0(Z))\leq  B-b-1  \right) - 1\right] \\
    = & \mathbb{E}_Z\left[   \mathbb{P}\left(\text{Bin}(B, F_0(Z))\leq  B-b-1  \right) - \mathbb{P}\left(\text{Bin}(B, F_0(Z))\leq  a-1  \right) \right] \\
    = & \mathbb{E}_Z \left[ \sum_{j = 0}^{B-b-1} \binom{B}{j} F_0(Z)^j (1 - F_0(Z))^{B-j} - \sum_{j = 0}^{a-1} \binom{B}{j} F_0(Z)^j (1 - F_0(Z))^{B-j} \right] \\
    = & \int_{0}^1\left\{  \sum_{j = 0}^{B-b-1} \binom{B}{j} y^j (1 - y)^{B-j} - \sum_{j = 0}^{a-1} \binom{B}{j} y^j (1 - y)^{B-j} \right\} d\mathbb P(F_0(Z) \leq y).
\end{align*}
If we let $R(y) = \mathbb{P}(F_0(Z) \leq y) - y$ we can decompose the above lower bound on coverage into the two following terms,
\begin{equation}
\label{eq:1,2,defn}
    \begin{split}
       \mathbb{P}(\psi(Z) \in [W_{(a)}, W_{(B-b)}]) \geq & \int_{0}^1\left\{  \sum_{j = 0}^{B-b-1} \binom{B}{j} y^j (1 - y)^{B-j} - \sum_{j = 0}^{a-1} \binom{B}{j} y^j (1 - y)^{B-j} \right\} d\mathbb P(F_0(Z) \leq y) \\
    = & \int_{0}^1\left\{  \sum_{j = 0}^{B-b-1} \binom{B}{j} y^j (1 - y)^{B-j} - \sum_{j = 0}^{a-1} \binom{B}{j} y^j (1 - y)^{B-j} \right\} dy \\
    &+ \int_{0}^1\left\{  \sum_{j = 0}^{B-b-1} \binom{B}{j} y^j (1 - y)^{B-j} - \sum_{j = 0}^{a-1} \binom{B}{j} y^j (1 - y)^{B-j} \right\} dR(y) \\
    \coloneqq &  \mathbf{I} + \mathbf{II}.  
    \end{split}
\end{equation}
We observe that the term $\mathbf{I}$ can be simplified as follows, 
\begin{align*}
    \mathbf{I} = & \int_{0}^1\left\{  \sum_{j = 0}^{B-b-1} \binom{B}{j} y^j (1 - y)^{B-j} - \sum_{j = 0}^{a-1} \binom{B}{j} y^j (1 - y)^{B-j} \right\} dy \\
    = &   \sum_{j = 0}^{B-b-1} \binom{B}{j} \int_{0}^1  y^j (1 - y)^{B-j} dy - \sum_{j = 0}^{a-1} \binom{B}{j} \int_{0}^1  y^j (1 - y)^{B-j} dy\\
    = &   \sum_{j = 0}^{B-b-1}  \binom{B}{j} \frac{j! (B-j)!}{(B+1)!} - \sum_{j = 0}^{a-1}  \binom{B}{j} \frac{j! (B-j)!}{(B+1)!}  \\
    = &    \sum_{j = 0}^{B-b-1} \frac{1}{B+1} - \sum_{j = 0}^{a-1} \frac{1}{B+1}\\
    = & \frac{B-a-b}{B+1}. 
\end{align*}
Consider the following definition of $G_k(y)$ for any $k \in [B]$,
\[
G_k(y) \coloneqq \sum_{j = 0}^{k} \binom{B}{j} y^j (1 - y)^{B-j} \quad \mbox{for all } y \in [0,1].
\]
Note that $y \mapsto G_k(y) = \mathbb{P}(\mbox{Bin}(B, y) \leq k )$ is a strictly decreasing function as $y$ increases in $(0,1)$. Therefore there exists a unique inverse function $\tau_k(\cdot)$ of $G_k(\cdot)$, such that for any $x \in (0,1)$, $\tau_k(x)$ is the unique solution of $G_k(\tau_k) = x$. We observe the following by Fubini's theorem for any $k \in [n]$,
\begin{equation}
\label{eq:r_simplify}
    \begin{split}
        \int_{0}^1  G_k(y) dR(y) = & \int_{0}^1 dR(y) \int_{0}^{G_k(y)} dx \\
        = & \int_0^1 dx \int_{\tau_k(x)}^{1} dR(y) \\
        \stackrel{(i)}{=} & - \int_0^1 R(\tau_k(x)) dx.
    \end{split}
\end{equation}
The step-$(i)$ follows from the fact that $R(1) = \mathbb{P}(F_0(Z) \leq 1 ) - 1 = 1 - 1 = 0$. Therefore we can simplify the expression for $\mathbf{II}$ as follows,
\begin{equation*}
    \begin{split}
        \mathbf{II} = &  \int_{0}^1\left\{  \sum_{j = 0}^{B-b-1} \binom{B}{j} y^j (1 - y)^{B-j} - \sum_{j = 0}^{a-1} \binom{B}{j} y^j (1 - y)^{B-j} \right\} dR(y) \\
        = & \int_{0}^1 \left\{  G_{B-b-1}(y) - G_{a-1}(y) \right\} dR(y) \\
        \stackrel{\eqref{eq:r_simplify}}{=} & \int_{0}^1 \left[ R(\tau_{a-1}(x)) - R(\tau_{B-b-1}(x))\right] dx \\
        = & \int_0^1 \left[\{\mathbb{P}(F_0(Z) \leq \tau_{a-1}(x)) - \tau_{a-1}(x)  \} - \{ \mathbb{P}(F_0(Z) \leq \tau_{B-b-1}(x)) - \tau_{B-b-1}(x) \}\right] dx \\
        = & \int_{0}^1 \left[\mathbb{P}(\tau_{B-b-1}(x) < F_0(Z) \leq \tau_{a-1}(x) ) - \mathbb{P} (\tau_{B-b-1}(x) < U(0,1) \leq \tau_{a-1}(x) )\right] dx \\
        \geq & - \int_{0}^1 \Delta dx \\
        = & - \Delta. 
    \end{split}
\end{equation*}
Combining all the computations, we have the following lower bound on coverage using \eqref{eq:1,2,defn},
\begin{equation*}
    \begin{split}
        \mathbb{P}(\psi(Z) \in [W_{(a)}, W_{(B-b)}]) \geq & \mathbf{I} + \mathbf{II} \\
        \geq & \frac{B-a-b}{B+1} - \Delta. 
    \end{split}
\end{equation*}
We follow the same steps to obtain an upper bound on the coverage. We have the following from \eqref{eq:equiv}, 
\begin{align*}
    \mathbb{P}(\psi(Z) \in [W_{(a)}, W_{(B-b)}]) = & \mathbb{E}_Z \left[ \mathbb{P}\left( \psi(Z) \in [W_{(a)}, W_{(B-b)}] \mid Z \right) \right]\\
    \leq & \mathbb{E}_Z \left[ \mathbb{P}\left( S_B(Z) \in [a, B-b] \mid Z \right) \right] \\
    = & \mathbb{E}_Z\left[ \mathbb{P}\left( a \leq \text{Bin}(B, F_0(Z)) \leq B - b  \right) \right] \\
    = & \mathbb{E}_Z\left[ \mathbb{P}\left(\text{Bin}(B, F_0(Z))\geq  a  \right) + \mathbb{P}\left(\text{Bin}(B, F_0(Z))\leq  B-b  \right) - 1\right] \\
    = & \mathbb{E}_Z\left[   \mathbb{P}\left(\text{Bin}(B, F_0(Z))\leq  B-b \right) - \mathbb{P}\left(\text{Bin}(B, F_0(Z))\leq  a-1  \right) \right] \\
    = & \mathbb{E}_Z \left[ \sum_{j = 0}^{B-b} \binom{B}{j} F_0(Z)^j (1 - F_0(Z))^{B-j} - \sum_{j = 0}^{a-1} \binom{B}{j} F_0(Z)^j (1 - F_0(Z))^{B-j} \right] \\
    = & \int_{0}^1\left\{  \sum_{j = 0}^{B-b} \binom{B}{j} y^j (1 - y)^{B-j} - \sum_{j = 0}^{a-1} \binom{B}{j} y^j (1 - y)^{B-j} \right\} d\mathbb P(F_0(Z) \leq y) \\
    =&  \mathbf{I}' + \mathbf{II}',
\end{align*}
where the terms $\mathbf{I}', \mathbf{II}'$ are defined as,
\begin{equation*}
    \begin{split}
        \mathbf{I}' = &\int_{0}^1\left\{  \sum_{j = 0}^{B-b} \binom{B}{j} y^j (1 - y)^{B-j} - \sum_{j = 0}^{a-1} \binom{B}{j} y^j (1 - y)^{B-j} \right\} dy, \\
        \mathbf{II}' = &  \int_{0}^1\left\{  \sum_{j = 0}^{B-b} \binom{B}{j} y^j (1 - y)^{B-j} - \sum_{j = 0}^{a-1} \binom{B}{j} y^j (1 - y)^{B-j} \right\} dR(y). 
    \end{split}
\end{equation*}
We can re-write the first term as follows,
\begin{align*}
    \mathbf{I}' = & \int_{0}^1\left\{  \sum_{j = 0}^{B-b} \binom{B}{j} y^j (1 - y)^{B-j} - \sum_{j = 0}^{a-1} \binom{B}{j} y^j (1 - y)^{B-j} \right\} dy \\
    = &   \sum_{j = 0}^{B-b} \binom{B}{j} \int_{0}^1  y^j (1 - y)^{B-j} dy - \sum_{j = 0}^{a-1} \binom{B}{j} \int_{0}^1  y^j (1 - y)^{B-j} dy\\
    = &   \sum_{j = 0}^{B-b}  \binom{B}{j} \frac{j! (B-j)!}{(B+1)!} - \sum_{j = 0}^{a-1}  \binom{B}{j} \frac{j! (B-j)!}{(B+1)!}  \\
    = &    \sum_{j = 0}^{B-b} \frac{1}{B+1} - \sum_{j = 0}^{a-1} \frac{1}{B+1}\\
    = & \frac{B+1-a-b}{B+1}. 
\end{align*}
We obtain an upper bound on $\mathbf{II}'$ using \eqref{eq:r_simplify} as follows,
\begin{equation*}
    \begin{split}
        \mathbf{II}' = & \int_{0}^1\left\{  \sum_{j = 0}^{B-b} \binom{B}{j} y^j (1 - y)^{B-j} - \sum_{j = 0}^{a-1} \binom{B}{j} y^j (1 - y)^{B-j} \right\} dR(y) \\
        = & \int_{0}^1 \left\{  G_{B-b}(y) - G_{a-1}(y) \right\} dR(y) \\
        \stackrel{\eqref{eq:r_simplify}}{=} & \int_{0}^1 \left[ R(\tau_{a-1}(x)) - R(\tau_{B-b}(x))\right] dx \\
        = & \int_0^1 \left[\{\mathbb{P}(F_0(Z) \leq \tau_{a-1}(x)) - \tau_{a-1}(x)  \} - \{ \mathbb{P}(F_0(Z) \leq \tau_{B-b}(x)) - \tau_{B-b}(x) \}\right] dx \\
        = & \int_{0}^1 \left[\mathbb{P}(\tau_{B-b}(x) < F_0(Z) \leq \tau_{a-1}(x) ) - \mathbb{P} (\tau_{B-b}(x) < U(0,1) \leq \tau_{a-1}(x) )\right] dx \\
        \leq & \int_{0}^1 \Delta dx \\
        = & \Delta. 
    \end{split}
\end{equation*}
We get the following upper bound on the coverage by combining all the prior results, 
\begin{equation*}
    \begin{split}
        \mathbb{P}(\psi(Z) \in [W_{(a)}, W_{(B-b)}]) \leq & \mathbf{I}' + \mathbf{II}' \\
        \leq & \frac{B+1-a-b}{B+1} + \Delta. 
    \end{split}
\end{equation*}
This completes the proof of the first statement of the theorem. The second and third statements of the theorem can be proved analogously by observing the following equivalence relations, 
\begin{equation*}
    \begin{split}
        W_{(a)} \leq \psi(Z) < W_{(B-b)} \iff &  a \leq S_B(Z) \leq B-b-1, \\
        W_{(a)} < \psi(Z) \leq W_{(B-b)} \iff & a \leq \widetilde S_B(Z) \leq B-b-1,
    \end{split}
\end{equation*}
where $\widetilde S_B(Z) = \sum_{i = 1}^B \textbf{1}(W_i < \psi(Z))$. The proof is completed by noting that conditioned on $Z$, we have $\widetilde S_B(Z) \sim \mbox{Bin}(B, \widetilde F_0(Z))$.

\section{Proof of Proposition~\ref{prop:ks_connection}}
\label{app:proof_ks_connection}
We define $H(x|Z) = \mathbb P(W_1 \leq x|Z)$ and $G(x) = \mathbb P(\psi(Z) \leq x)$ for any $x \in \mathbb R$. We observe the following,
\begin{equation*}
    \begin{split}
        |F_0(Z) - G(\psi(Z))| = & |\mathbb P(W_1 \leq \psi(Z)|Z)  - G(\psi(Z)) | \\
        = & | H(\psi(Z)|Z) - G(\psi(Z)) | \\
        \leq & \sup_{x \in \mathbb R} |H(x|Z) - G(x) | \\
        = & \Delta^*(Z). 
    \end{split}
\end{equation*}
We fix $\epsilon > 0$. Using the above inequality we have the following for any $ t \in \mathbb R$,
\begin{equation*}
    \begin{split}
        \{F_0(Z) \leq t\} \subset &  \{G(\psi(Z)) \leq t+ \epsilon \} \cup \{\Delta^*(Z) > \epsilon\}, \\
        \{ G(\psi(Z)) \leq t - \epsilon\} \subset & \{F_0(Z) \leq t \} \cup \{\Delta^*(Z) > \epsilon\}.
    \end{split}
\end{equation*}
This implies that for any $t \in \mathbb R$ we have,
\begin{equation*}
    \begin{split}
        \mathbb P(F_0(Z) \leq t ) \leq & \mathbb P (G(\psi(Z)) \leq t+ \epsilon) + \mathbb P( \Delta^*(Z) > \epsilon), \\
        \mathbb P(F_0(Z) \leq t )  \geq &  \mathbb P (G(\psi(Z)) \leq t - \epsilon) - \mathbb P( \Delta^*(Z) > \epsilon). 
    \end{split}
\end{equation*}
In other words we have the following bounds for any $t \in \mathbb R$,
\begin{equation*}
    \begin{split}
        &\{ \mathbb P (G(\psi(Z)) \leq t - \epsilon) - (t - \epsilon)\} - \epsilon - \mathbb P( \Delta^*(Z) > \epsilon) \\
        \leq &  \mathbb P(F_0(Z) \leq t ) - t \\
        \leq & \{ \mathbb P (G(\psi(Z)) \leq t + \epsilon) - (t + \epsilon)\} + \epsilon + \mathbb P( \Delta^*(Z) > \epsilon).
    \end{split}
\end{equation*}
Taking supremum over all $t \in \mathbb R$ we have,
\begin{equation*}
    \begin{split}
        d_{\mathrm{KS}}(F_0(Z) , U(0,1)) =& \sup_{t \in \mathbb R} | \mathbb P(F_0(Z) \leq t ) - t  | \\
        \leq & \sup_{t \in \mathbb R} | \{ \mathbb P (G(\psi(Z)) \leq t ) - t\} | + \epsilon + \mathbb P( \Delta^*(Z) > \epsilon) \\
        = & d_{\mathrm{KS}}(G(\psi(Z)), U(0,1))  +  \epsilon + \mathbb P( \Delta^*(Z) > \epsilon). 
    \end{split}
\end{equation*}
We will show that under the given assumptions, $ d_{\mathrm{KS}}(G(\psi(Z)), U(0,1))  \leq \eta$. Let us define the quantile function of $\psi(Z)$ as $G^{-1}(t)=\inf\{x\in\mathbb{R}\mid G(x)\geq t\}$ for any $t \in \mathbb R$. To see this, we observe that for any $t \in \mathbb R$,
\begin{equation*}
    \begin{split}
        \mathbb P (G(\psi(Z)) \leq t) \leq& \mathbb P(\psi(Z) < G^{-1}(t)) + \mathbb P(\{\psi(Z) = G^{-1}(t) \} \cap \{G(G^{-1}(t)) =t \} ) \\
        = &  \mathbb P(\psi(Z) < G^{-1}(t)) + t - \mathbb P(\psi(Z) < G^{-1}(t)) \\
        = & t, \\
          \mathbb P (G(\psi(Z)) \leq t) \geq & \mathbb{P}(\psi(Z) < G^{-1}(t)) \\
          = & \mathbb{P}(\psi(Z) \leq G^{-1}(t)) - \mathbb{P}(\psi(Z) = G^{-1}(t)) \\
          \geq &t - \eta. 
    \end{split}
\end{equation*}
The last line follows from the following fact,
\begin{equation*}
    \begin{split}
         \mathbb{P}(\psi(Z) = G^{-1}(t)) \leq &\inf_{\mathfrak W \rightarrow 0}  \mathbb{P}(\psi(Z) \in (G^{-1}(t) - (\mathfrak W/2), G^{-1}(t) + (\mathfrak W/2)) ) \\
         \leq & \inf_{\mathfrak W \rightarrow 0} \{ c \mathfrak W + \eta \} \\
         = & \eta. 
    \end{split}
\end{equation*}
Therefore we see that $d_{\mathrm{KS}}(G(\psi(Z)), U(0,1)) = \sup_{t \in \mathbb R} | \mathbb P (G(\psi(Z)) \leq t) - t | \leq \eta$ and hence the following holds,
\[
d_{\mathrm{KS}}(F_0(Z) , U(0,1)) \leq \eta +  \epsilon + \mathbb P( \Delta^*(Z) > \epsilon). 
\]
This completes the proof of the proposition. 
\section{Proof of \Cref{rem:tighter_cov_bound}}
\label{app:proof_tighter_cov}
We begin the proof by defining the following class of random variables, 
\[
\mathcal{A}
:=
\left\{
H|\ H: \Omega \mapsto [0,1]\mbox{ is measurable and } \mathbb E[G_k(H)]=\frac{k+1}{B+1}
\ \text{for every }k\in\{0,1,\dots,B\}
\right\},
\]
where $G_k(w)$ is defined for $k\in\{0,1,\dots,B\}$ as, 
\[
G_k(w) := \mathbb{P}\big(\mathrm{Binomial}(B,w)\le k\big)
= \sum_{j=0}^k \binom{B}{j} w^j(1-w)^{B-j},\qquad w\in[0,1].
\]
Recall that $\mathcal H = \{H | H :\Omega \mapsto \mathbb R \mbox{ is measurable and } \mathbb E[H^r] = 1/(r +1) \mbox{ for } r \in \{1,\ldots,B\}\} $ is the class of univariate random variables whose first $B$ raw moments match with that of $U \stackrel{d}{=}U(0,1)$. We will show that $\mathcal{A} = \mathcal{H}$.
Define the Bernstein basis polynomials of degree $B$ by
\[
b_{j,B}(w) := \binom{B}{j} w^j(1-w)^{B-j},\qquad j=0,1,\dots,B,\ \ w\in[0,1].
\]
Then we have $G_k(w)=\sum_{j=0}^k b_{j,B}(w)$. We first show that if $H\in\mathcal{A}$ then $\mathbb{E}[H^r]=1/(r+1)$ for $r=1,\dots,B$ i.e.\ $W \in \mathcal{W}$. If $W\in\mathcal{A}$ then
\[
\mathbb{P}\big(\mathrm{Binomial}(B,H)\le k\big)
=
\mathbb{E}\big[G_k(H)\big] = \frac{k+1}{B+1},\qquad k=0,1,\dots,B.
\]
For $k\ge 1$, subtract the $k-1$ equation from the $k$ equation to obtain
\[
\mathbb{E}\big[G_k(H)-G_{k-1}(H)\big]
=
\frac{1}{B+1} ,\qquad k=0,1,\dots,B.
\]
Because $G_k-G_{k-1}=b_{k,B}$, this yields
\begin{equation}
\label{eq:bernstein_expectation_constant}
\mathbb{E}\big[b_{k,B}(H)\big]=\frac{1}{B+1},\qquad k=0,1,\dots,B.
\end{equation}
Fix $r\in\{0,1,\dots,B\}$. We next express the monomial $w^r$ in the Bernstein basis. Let $K_w\sim\mathrm{Binomial}(B,w)$ for fixed $w\in[0,1]$. Then the standard binomial identity gives,
\[
\mathbb{E}\!\left[\binom{K_w}{r}\right]=\binom{B}{r}w^r.
\]
On the other hand the following also holds,
\[
\mathbb{E}\!\left[\binom{K_w}{r}\right]
=
\sum_{k=r}^B \binom{k}{r}\,\mathbb{P}(K_w=k)
=
\sum_{k=r}^B \binom{k}{r}\,b_{k,B}(w).
\]
Therefore, for every $w\in[0,1]$ we have,
\begin{equation}
\label{eq:monomial_bernstein_representation}
w^r
=
\sum_{k=r}^B \frac{\binom{k}{r}}{\binom{B}{r}}\,b_{k,B}(w).
\end{equation}
Taking expectations in \eqref{eq:monomial_bernstein_representation} with $w=H$ and using
\eqref{eq:bernstein_expectation_constant} gives
\[
\mathbb{E}[H^r]
=
\sum_{k=r}^B \frac{\binom{k}{r}}{\binom{B}{r}}\,
\mathbb{E}\big[b_{k,B}(H)\big]
=
\frac{1}{B+1}\sum_{k=r}^B \frac{\binom{k}{r}}{\binom{B}{r}}.
\]
Using the identity $\sum_{k=r}^B \binom{k}{r}=\binom{B+1}{r+1}$, we obtain
\[
\mathbb{E}[H^r]
=
\frac{1}{B+1}\cdot \frac{\binom{B+1}{r+1}}{\binom{B}{r}}
=
\frac{1}{B+1}\cdot
\frac{\frac{(B+1)!}{(r+1)!(B-r)!}}{\frac{B!}{r!(B-r)!}}
=
\frac{1}{B+1}\cdot \frac{B+1}{r+1}
=
\frac{1}{r+1}.
\]
Since this holds for any $r \in \{0,1,\ldots,B \}$, we conclude that $ H \in \mathcal H$. We will now show that if $H \in \mathcal H \implies H \in \mathcal A$. Suppose $ H \in \mathcal H$ i.e.\ the following holds,
\[
\mathbb{E}[H^r]=\frac{1}{r+1} = \mathbb E [U^r],\qquad r=0,1,\dots,B.
\]
Fix $k\in\{0,1,\dots,B\}$. Consider $b_{k,B}(w)=\binom{B}{k}w^k(1-w)^{B-k}$.
Expanding $(1-w)^{B-k}$ gives
\[
b_{k,B}(w)
=
\binom{B}{k}\sum_{j=0}^{B-k}\binom{B-k}{j}(-1)^j w^{k+j},
\]
which is a polynomial of degree $B$. Taking expectations, only moments up to order $B$ appear, hence
\[
\mathbb{E}\big[b_{k,B}(H)\big]
=
\binom{B}{k}\sum_{j=0}^{B-k}\binom{B-k}{j}(-1)^j\,\mathbb{E}[H^{k+j}]
=
\binom{B}{k}\sum_{j=0}^{B-k}\binom{B-k}{j}(-1)^j\,\mathbb{E}[U^{k+j}]
=
\mathbb{E}\big[b_{k,B}(U)\big].
\]
We can compute $\mathbb{E}[b_{k,B}(U)]$ directly as, 
\[
\mathbb{E}\big[b_{k,B}(U)\big]
=
\int_0^1 \binom{B}{k}u^k(1-u)^{B-k}\,\mathrm{d}u
=
\binom{B}{k}\,\mathrm{Beta}(k+1,B-k+1),
\]
Using $\mathrm{Beta}(k+1,B-k+1)=\dfrac{k!(B-k)!}{(B+1)!}$ gives,
\[
\mathbb{E}\big[b_{k,B}(U)\big]
=
\binom{B}{k}\frac{k!(B-k)!}{(B+1)!}
=
\frac{B!}{k!(B-k)!}\cdot \frac{k!(B-k)!}{(B+1)!}
=
\frac{1}{B+1}.
\]
Therefore we have,
\[
\mathbb{E}\big[b_{k,B}(H)\big]=\frac{1}{B+1},\qquad k=0,1,\dots,B.
\]
Summing this from $j=0$ to $k$ we have,
\[
\mathbb{E}\big[G_k(H)\big]
=
\mathbb{E}\!\left[\sum_{j=0}^k b_{j,B}(H)\right]
=
\sum_{j=0}^k \mathbb{E}\big[b_{j,B}(H)\big]
=
\sum_{j=0}^k \frac{1}{B+1}
=
\frac{k+1}{B+1}.
\]
Since this holds for every $k\in\{0,1,\dots,B\}$, we have $H\in\mathcal{A}$. Combining both the directions completes the proof of the fact $H \in \mathcal{A} \iff H \in \mathcal H$. 

We wish to bound the probability:
\[
\mathbb{P}\left(\psi(Z) \in [W_{(a)}, W_{(B-b)}] \right), \quad \mbox{for } 0 \leq a < B-b \leq B. 
\]
Following \cite{hall1986number} we observe that
\[
\mathbb{P}\left(\psi(Z) \in [W_{(a)}, W_{(B-b)}] \right) = \mathbb{E}_Z \left[ \mathbb{P}\left( \psi(Z) \in [W_{(a)}, W_{(B-b)}] \mid Z \right) \right].
\]
We consider the following definition,
\[
S_B(Z) := \sum_{i=1}^B \mathbf{1}(W_i \leq \psi(Z)) \sim \text{Bin}(B, F_0(Z)) \quad \mbox{conditioned on }Z.
\]
We also note the following equivalence,
\begin{align}
\label{eq:equiv_remark}
W_{(a)} \leq \psi(Z) \leq W_{(B-b)} \begin{cases}
   &\Rightarrow a \leq S_B(Z) \leq B - b, \\ 
   & \Leftarrow a \leq S_B(Z) \leq B-b-1.
\end{cases} 
\end{align}

Hence:
\begin{align*}
    \mathbb{P}(\psi(Z) \in [W_{(a)}, W_{(B-b)}]) = & \mathbb{E}_Z \left[ \mathbb{P}\left( \psi(Z) \in [W_{(a)}, W_{(B-b)}] \mid Z \right) \right]\\
    \geq & \mathbb{E}_Z \left[ \mathbb{P}\left( S_B(Z) \in [a, B-b-1] \mid Z \right) \right] \\
    = & \mathbb{E}_Z\left[ \mathbb{P}\left( a \leq \text{Bin}(B, F_0(Z)) \leq B - b -1 \right) \right] \\
    = & \mathbb{E}_Z\left[ \mathbb{P}\left(\text{Bin}(B, F_0(Z))\geq  a  \right) + \mathbb{P}\left(\text{Bin}(B, F_0(Z))\leq  B-b-1  \right) - 1\right] \\
    = & \mathbb{E}_Z\left[   \mathbb{P}\left(\text{Bin}(B, F_0(Z))\leq  B-b-1  \right) - \mathbb{P}\left(\text{Bin}(B, F_0(Z))\leq  a-1  \right) \right] \\
    = & \mathbb{E}_Z \left[ \sum_{j = 0}^{B-b-1} \binom{B}{j} F_0(Z)^j (1 - F_0(Z))^{B-j} - \sum_{j = 0}^{a-1} \binom{B}{j} F_0(Z)^j (1 - F_0(Z))^{B-j} \right] \\
    = & \int_{0}^1\left\{  \sum_{j = 0}^{B-b-1} \binom{B}{j} y^j (1 - y)^{B-j} - \sum_{j = 0}^{a-1} \binom{B}{j} y^j (1 - y)^{B-j} \right\} d\mathbb P(F_0(Z) \leq y).
\end{align*}
Fix any $H \in \mathcal{A}$. If we let $R(y) = \mathbb{P}(F_0(Z) \leq y) - \mathbb P(H \leq y)$ we can decompose the above lower bound on coverage into the two following terms,
\begin{equation}
\label{eq:1,2,defn_remark}
    \begin{split}
       \mathbb{P}(\psi(Z) \in [W_{(a)}, W_{(B-b)}]) \geq & \int_{0}^1\left\{  \sum_{j = 0}^{B-b-1} \binom{B}{j} y^j (1 - y)^{B-j} - \sum_{j = 0}^{a-1} \binom{B}{j} y^j (1 - y)^{B-j} \right\} d\mathbb P(F_0(Z) \leq y) \\
    = &\int_{0}^1\left\{  \sum_{j = 0}^{B-b-1} \binom{B}{j} y^j (1 - y)^{B-j} - \sum_{j = 0}^{a-1} \binom{B}{j} y^j (1 - y)^{B-j} \right\} d \mathbb P( H \leq y) \\
    &+ \int_{0}^1\left\{  \sum_{j = 0}^{B-b-1} \binom{B}{j} y^j (1 - y)^{B-j} - \sum_{j = 0}^{a-1} \binom{B}{j} y^j (1 - y)^{B-j} \right\} dR(y) \\
    \coloneqq &  \mathbf{I} + \mathbf{II}.  
    \end{split}
\end{equation}
We observe that the term $\mathbf{I}$ can be simplified as follows, 
\begin{align*}
    \mathbf{I} = & \int_{0}^1\left\{  \sum_{j = 0}^{B-b-1} \binom{B}{j} y^j (1 - y)^{B-j} - \sum_{j = 0}^{a-1} \binom{B}{j} y^j (1 - y)^{B-j} \right\} d \mathbb P( H \leq y) \\
    = & \int_{0}^1G_{B-b-1}(y)  d\mathbb P(H \leq y) - \int_0^1 G_{a-1}(y) d\mathbb{P}(H \leq y)    \\
    = &\mathbb E[G_{B-b-1}(H) - G_{a-1}(H)] \\
    \stackrel{(i)}{=} &     \frac{B-b}{B+1} -  \frac{a}{B+1}\\
    = & \frac{B-a-b}{B+1}. 
\end{align*}
The equality $(i)$ follows from the definition of the class $\mathcal A$. 
Note that $y \mapsto G_k(y) = \mathbb{P}(\mbox{Bin}(B, y) \leq k )$ is a strictly decreasing function as $y$ increases in $(0,1)$. Therefore there exists a unique inverse function $\tau_k(\cdot)$ of $G_k(\cdot)$, such that for any $x \in (0,1)$, $\tau_k(x)$ is the unique solution of $G_k(\tau_k) = x$. We observe the following by Fubini's theorem for any $k \in [n]$,
\begin{equation}
\label{eq:r_simplify_remark}
    \begin{split}
        \int_{0}^1  G_k(y) dR(y) = & \int_{0}^1 dR(y) \int_{0}^{G_k(y)} dx \\
        = & \int_0^1 dx \int_{\tau_k(x)}^{1} dR(y) \\
        \stackrel{(i)}{=} & - \int_0^1 R(\tau_k(x)) dx.
    \end{split}
\end{equation}
The step-$(i)$ follows from the fact that $R(1) = \mathbb{P}(F_0(Z) \leq 1 ) - 1 = 1 - 1 = 0$. Therefore we can simplify the expression for $\mathbf{II}$ as follows,
\begin{equation*}
    \begin{split}
        \mathbf{II} = & \int_{0}^1\left\{  \sum_{j = 0}^{B-b-1} \binom{B}{j} y^j (1 - y)^{B-j} - \sum_{j = 0}^{a-1} \binom{B}{j} y^j (1 - y)^{B-j} \right\} dR(y) \\
        = & \int_{0}^1 \left\{  G_{B-b-1}(y) - G_{a-1}(y) \right\} dR(y) \\
        \stackrel{\eqref{eq:r_simplify_remark}}{=} & \int_{0}^1 \left[ R(\tau_{a-1}(x)) - R(\tau_{B-b-1}(x))\right] dx \\
        = & \int_0^1 \left[\{\mathbb{P}(F_0(Z) \leq \tau_{a-1}(x)) - \mathbb P(H \leq\tau_{a-1}(x))  \} - \{ \mathbb{P}(F_0(Z) \leq \tau_{B-b-1}(x)) - \mathbb P(H \leq \tau_{B-b-1}(x) )\}\right] dx \\
        = & \int_{0}^1 \left[\mathbb{P}(\tau_{B-b-1}(x) < F_0(Z) \leq \tau_{a-1}(x) ) - \mathbb{P} (\tau_{B-b-1}(x) < H \leq \tau_{a-1}(x) )\right] dx \\
        \geq & - \int_{0}^1 \widetilde d_{\mathrm{KS}}(F_0(Z), H) dx \\
        = & - \widetilde d_{\mathrm{KS}}(F_0(Z), H). 
    \end{split}
\end{equation*}
Combining all the computations, we have the following lower bound on coverage using \eqref{eq:1,2,defn_remark},
\begin{equation*}
    \begin{split}
        \mathbb{P}(\psi(Z) \in [W_{(a)}, W_{(B-b)}]) \geq & \mathbf{I} + \mathbf{II} \\
        \geq & \frac{B-a-b}{B+1} - \widetilde d_{\mathrm{KS}}(F_0(Z), H). 
    \end{split}
\end{equation*}
The above lower bound on coverage holds for any $ H \in \mathcal A$. Therefore we have the following stronger coverage guarantee,
\begin{equation*}
    \begin{split}
        \mathbb{P}(\psi(Z) \in [W_{(a)}, W_{(B-b)}]) 
        \geq & \frac{B-a-b}{B+1} - \inf_{H \in \mathcal A}\widetilde d_{\mathrm{KS}}(F_0(Z), H) \\
        = & \frac{B-a-b}{B+1} - \inf_{H \in \mathcal H}\widetilde d_{\mathrm{KS}}(F_0(Z), H).
    \end{split}
\end{equation*}
The last equality follows from $\mathcal A = \mathcal H$. We follow the same steps to obtain an upper bound on the coverage. Fix any $H \in \mathcal A$. We have the following from \eqref{eq:equiv_remark}, 
\begin{align*}
    \mathbb{P}(\psi(Z) \in [W_{(a)}, W_{(B-b)}]) = & \mathbb{E}_Z \left[ \mathbb{P}\left( \psi(Z) \in [W_{(a)}, W_{(B-b)}] \mid Z \right) \right]\\
    \leq & \mathbb{E}_Z \left[ \mathbb{P}\left( S_B(Z) \in [a, B-b] \mid Z \right) \right] \\
    = & \mathbb{E}_Z\left[ \mathbb{P}\left( a \leq \text{Bin}(B, F_0(Z)) \leq B - b  \right) \right] \\
    = & \mathbb{E}_Z\left[ \mathbb{P}\left(\text{Bin}(B, F_0(Z))\geq  a  \right) + \mathbb{P}\left(\text{Bin}(B, F_0(Z))\leq  B-b  \right) - 1\right] \\
    = & \mathbb{E}_Z\left[   \mathbb{P}\left(\text{Bin}(B, F_0(Z))\leq  B-b \right) - \mathbb{P}\left(\text{Bin}(B, F_0(Z))\leq  a-1  \right) \right] \\
    = & \mathbb{E}_Z \left[ \sum_{j = 0}^{B-b} \binom{B}{j} F_0(Z)^j (1 - F_0(Z))^{B-j} - \sum_{j = 0}^{a-1} \binom{B}{j} F_0(Z)^j (1 - F_0(Z))^{B-j} \right] \\
    = & \int_{0}^1\left\{  \sum_{j = 0}^{B-b} \binom{B}{j} y^j (1 - y)^{B-j} - \sum_{j = 0}^{a-1} \binom{B}{j} y^j (1 - y)^{B-j} \right\} d\mathbb P(F_0(Z) \leq y) \\
    =&  \mathbf{I}' + \mathbf{II}',
\end{align*}
where the terms $\mathbf{I}', \mathbf{II}'$ are defined as,
\begin{equation*}
    \begin{split}
        \mathbf{I}' = &\int_{0}^1\left\{  \sum_{j = 0}^{B-b} \binom{B}{j} y^j (1 - y)^{B-j} - \sum_{j = 0}^{a-1} \binom{B}{j} y^j (1 - y)^{B-j} \right\} d\mathbb P(H \leq y), \\
        \mathbf{II}' = &  \int_{0}^1\left\{  \sum_{j = 0}^{B-b} \binom{B}{j} y^j (1 - y)^{B-j} - \sum_{j = 0}^{a-1} \binom{B}{j} y^j (1 - y)^{B-j} \right\} dR(y). 
    \end{split}
\end{equation*}
We can re-write the first term as follows,
\begin{align*}
    \mathbf{I}' = & \int_{0}^1\left\{  \sum_{j = 0}^{B-b} \binom{B}{j} y^j (1 - y)^{B-j} - \sum_{j = 0}^{a-1} \binom{B}{j} y^j (1 - y)^{B-j} \right\} d\mathbb P(H \leq y) \\
    = &  \int_0^1 G_{B-b}(y) d\mathbb P(H \leq y)  - \int_0^1 G_{a-1}(y) d\mathbb P(H \leq y) \\
    =& \mathbb E[ G_{B-b}(H) - G_{a-1}(H) ]  \\
     \stackrel{(i)}{=} &    \frac{B-b+1}{B+1} - \frac{a}{B+1}\\
    = & \frac{B+1-a-b}{B+1}. 
\end{align*}
The step-$(i)$ follows from the definition of the class $\mathcal A$. We obtain an upper bound on $\mathbf{II}'$ using \eqref{eq:r_simplify_remark} as follows,
\begin{equation*}
    \begin{split}
        \mathbf{II}' = & \int_{0}^1\left\{  \sum_{j = 0}^{B-b} \binom{B}{j} y^j (1 - y)^{B-j} - \sum_{j = 0}^{a-1} \binom{B}{j} y^j (1 - y)^{B-j} \right\} dR(y) \\
        = & \int_{0}^1 \left\{  G_{B-b}(y) - G_{a-1}(y) \right\} dR(y) \\
        \stackrel{\eqref{eq:r_simplify}}{=} & \int_{0}^1 \left[ R(\tau_{a-1}(x)) - R(\tau_{B-b}(x))\right] dx \\
        = & \int_0^1 \left[\{\mathbb{P}(F_0(Z) \leq \tau_{a-1}(x)) - \mathbb P(H \leq \tau_{a-1}(x) ) \} - \{ \mathbb{P}(F_0(Z) \leq \tau_{B-b}(x)) -\mathbb P(H \leq  \tau_{B-b}(x) )\}\right] dx \\
        = & \int_{0}^1 \left[\mathbb{P}(\tau_{B-b}(x) < F_0(Z) \leq \tau_{a-1}(x) ) - \mathbb{P} (\tau_{B-b}(x) < H \leq \tau_{a-1}(x) )\right] dx \\
        \leq & \int_{0}^1 \widetilde d_{\mathrm{KS}}(F_0(Z), H) dx \\
        = & \widetilde d_{\mathrm{KS}}(F_0(Z), H). 
    \end{split}
\end{equation*}
We get the following upper bound on the coverage by combining all the prior results and by observing that the above bonds hold for any $H \in \mathcal A$, 
\begin{equation*}
    \begin{split}
        \mathbb{P}(\psi(Z) \in [W_{(a)}, W_{(B-b)}])  \leq & \frac{B+1-a-b}{B+1} + \inf_{H \in \mathcal A}\widetilde d_{\mathrm{KS}}(F_0(Z), H) \\
        = &  \frac{B+1-a-b}{B+1} + \inf_{H \in \mathcal H}\widetilde d_{\mathrm{KS}}(F_0(Z), H).
    \end{split}
\end{equation*}
The last step follows from the fact $\mathcal A = \mathcal H$. This completes the proof of the remark.
\section{Proof of Theorem~\ref{thm:main_result_extension}}
\label{app:proof_main_res_ext}
As in the proof of \Cref{thm:main_result} we observe that, 
\[
\mathbb{P}\left(\psi(Z) \in [W_{(a)}, W_{(B-b)}] \right) = \mathbb{E}_Z \left[ \mathbb{P}\left( \psi(Z) \in [W_{(a)}, W_{(B-b)}] \mid Z \right) \right].
\]
We consider the following definition,
\[
S_B(Z) := \sum_{i=1}^B \mathbf{1}(W_i \leq \psi(Z)) \sim \text{Poi-Bin}(B, F_1(Z),\ldots,F_B(Z)) \quad \mbox{conditioned on }Z,
\]
where $\mbox{Poi-Bin}(B,p_1,\ldots,p_B)$ is a random variable distributed as $\sum_{i = 1}^B \mbox{Bernoulli}(p_i)$ and $F_i(z) = \mathbb P(W_i \leq \psi(Z)| Z = z)$ for all $i \in [B]$. We know from Corollary $2$ of \cite{ehm1991binomial} that, 
\[
C(1 - p^{B+1} - q^{B+1})r \leq d_{\mathrm{TV}}(\mbox{Poi-Bin}(B,p_1,\ldots,p_B), \mbox{Bin}(B, p)) \leq (B/(B+1))(1 - p^{B+1} - q^{B+1})r,
\]
where $C \geq 1/124$ is an universal constant, $p = \sum_{i = 1}^B p_i/B$, $q = 1-p$, $r = 1 - (Bpq)^{-1}\mbox{Var}(\mbox{Poi-Bin}(B,p_1,\ldots,p_B)) =1 - (Bpq)^{-1} \sum_{i = 1}^B p_i(1 - p_i)$. 
We recall the equivalence statements in \eqref{eq:equiv},
\begin{align}
W_{(a)} \leq \psi(Z) \leq W_{(B-b)} \begin{cases}
   &\Rightarrow a \leq S_B(Z) \leq B - b, \\ 
   & \Leftarrow a \leq S_B(Z) \leq B-b-1.
\end{cases} 
\end{align}
Let $I \sim \mbox{Discrete-Uniform}\{1,\ldots,B\}$. We can lower bound the coverage as follows using Theorem $1$ (Lemma $2$) of \cite{ehm1991binomial} where $\overline F(z) = \mathbb E[F_I(z)] =  (1/B)\sum_{i = 1}^BF_i(z)$ plays the role of $p= (1/B)\sum_{i = 1}^B p_i$,
\begin{align*}
    &\mathbb{P}(\psi(Z) \in [W_{(a)}, W_{(B-b)}]) \\
    = & \mathbb{E}_Z \left[ \mathbb{P}\left( \psi(Z) \in [W_{(a)}, W_{(B-b)}] \mid Z \right) \right]\\
    \geq & \mathbb{E}_Z \left[ \mathbb{P}\left( S_B(Z) \in [a, B-b-1] \mid Z \right) \right] \\
    = & \mathbb{E}_Z\left[ \mathbb{P}\left( a \leq \text{Poi-Bin}(B, F_1(Z),\ldots,F_B(Z)) \leq B - b -1 \right) \right] \\
    = & \mathbb{E}_Z\left[ \mathbb{P}\left( a \leq \text{Bin}(B, \overline F(Z)) \leq B - b -1 \right) - d_{\mathrm{TV}}\left(\text{Poi-Bin}(B, F_1(Z),\ldots,F_B(Z)), \text{Bin}(B, \overline F(Z)) \right) \right] \\
    \geq &\mathbb{E}_Z\left[ \mathbb{P}\left( a \leq \text{Bin}(B, \overline F(Z)) \leq B - b -1 \right) \right] - \mathbb{E}_Z\left[ B\Delta(Z) \min \left\{1, \frac{1}{B \overline F(Z) (1 - \overline F(Z))} \right\} \right],
\end{align*}
where with slight abuse of notation, we have $\Delta(z) = \mathrm{Var}(F_I(z))$ in this proof. Following the proof of \Cref{thm:main_result} we know that the first term in the above bound can be lower bounded as follows,
\[
\mathbb{E}_Z\left[ \mathbb{P}\left( a \leq \text{Bin}(B, \overline F(Z)) \leq B - b -1 \right) \right] \geq  1- \frac{a+b +1}{B+1} - \widetilde d_{\mathrm{KS}}(\overline F(Z), U(0,1)). 
\]
Therefore we get the following lower bound on coverage,
\begin{align*}
    \mathbb{P}(\psi(Z) \in [W_{(a)}, W_{(B-b)}])  \geq &  1- \frac{a+b +1}{B+1} - \widetilde d_{\mathrm{KS}}(\overline F(Z), U(0,1)) - \mathbb{E}_Z\left[B \Delta(Z) \min \left\{1, \frac{1}{B \overline F(Z) (1 - \overline F(Z))} \right\} \right] .
\end{align*}
We recall that $\kappa_i = \sup_{u \in \mathcal Z} | F(\psi(u)) - F_i(u)|$ for all $i \in [B]$. This implies that for all $u \in \mathcal Z$ and for all $i \in [B]$ we have $| F(\psi(u)) - F_i(u)| \leq \kappa_i$. Therefore we have the following inequality,
\begin{equation*}
    \begin{split}
        \Delta(Z) = & \mathrm{Var}_I(F_I(Z)) \\
        = & \frac{1}{B}\sum_{i = 1}^B(F_i(Z) - \overline F(Z))^2 \\
        \leq & \frac{1}{B}\sum_{i = 1}^B(F_i(Z) -  F(\psi(Z)))^2  \\
        \leq & \frac{1}{B}\sum_{i = 1}^B \{\sup_{u \in \mathcal Z}|F_i(u) - F(\psi(u))| \}^2 \\
        = &  \frac{1}{B} \sum_{i = 1}^B \kappa_i^2. 
    \end{split}
\end{equation*}
We further observe the following,
\begin{equation*}
    \begin{split}
       & \mathbb E_Z \left[  \min \left\{1, \frac{1}{B \overline F(Z) (1 - \overline F(Z))} \right\} \right] \\
       = & \int_0^1  \left[  \min \left\{1, \frac{1}{B y (1 - y)} \right\} \right] d\mathbb P(\overline F(Z) \leq y) \\
       = &  \int_0^1  \left[  \min \left\{1, \frac{1}{B y (1 - y)} \right\} \right] dy +  \int_0^1  \left[  \min \left\{1, \frac{1}{B y (1 - y)} \right\} \right] d\overline R(y) \\
       \coloneqq & \mathbf I + \mathbf {II},
    \end{split}
\end{equation*}
where $\overline R(y) = \mathbb P(\overline F(Z) \leq y) - y $. We note that if $x \in (0,1)$ then there exists $0 \leq \tau_1(x) \leq \tau_2(x)\leq 1$ such that $1/(By(1-y)) \geq x \iff y \in [0,\tau_1(x)] \cup [\tau_2(x), 1]$. We can bound the second term as follows,
\begin{equation*}
    \begin{split}
        \mathbf {II} = & \int_0^1  \left[  \min \left\{1, \frac{1}{B y (1 - y)} \right\} \right] d\overline R(y) \\
        = & \int_0^1 d\overline R(y) \int_{0}^{\left[  \min \left\{1, \frac{1}{B y (1 - y)} \right\} \right]} dx \\
        = & \int_{0}^1 dx \int_{y \in [0,\tau_1(x)] \cup [\tau_2(x), 1]} d\overline R(y) \\
        = & \int_{0}^1 dx  \left[\int_0^1 d\overline R(y) - \int_{\tau_1(x)}^{\tau_2(x)} d\overline R(y)  \right] \\
        = & \int_{0}^1 \left[0 - \overline R(\tau_2(x))  + \overline R(\tau_1(x)) \right]   dx \\
        \leq & \widetilde d_{\mathrm{KS}}(\overline F(Z), U(0,1)). 
    \end{split}
\end{equation*}
We will now compute the integral in the first term. For $B \leq 4$ we have $By(1-y) \leq B/4 \leq 1$. Therefore $\mathbf I = 1$ for $B \leq 4$. For $B > 4$ the quadratic equation $By(1 - y) = 1$ admits two solutions $a_B = [1 - \sqrt{1 - (4/B)}]/2$ and $b_B = [1 + \sqrt{1 - (4/B)}]/2$. Therefore the following holds,
\begin{equation*}
    \begin{split}
        \mathbf I = & \int_0^1  \left[  \min \left\{1, \frac{1}{B y (1 - y)} \right\} \right] dy \\
        = & \int_{[0,a_B]\cup [b_B,1]} 1 dy + \int_{a_B}^{b_B} \frac{dy}{By(1-y)} \\
        = & a_B + 1 - b_B + \frac{1}{B}\left[\log \left(\frac{y}{1 - y} \right)\right]_{a_B}^{b_B} \\
        = & 1 - \sqrt{1 - \frac{4}{B}} + \frac{2}{B}\log \left( \frac{1 + \sqrt{1 - (4/B)}}{1 -\sqrt{1 - (4/B)}}\right).
    \end{split}
\end{equation*}
Therefore $\mathbf I = I_B$ where,
\begin{equation*}
    I_B = \begin{cases}
        1 & \mbox{if} \quad B \leq 4, \\
       1 - \sqrt{1 - \frac{4}{B}} + \frac{2}{B}\log \left( \frac{1 + \sqrt{1 - (4/B)}}{1 -\sqrt{1 - (4/B)}}\right) & \mbox{if} \quad B > 4.
    \end{cases}
\end{equation*}
Combining these we get the following lower bound on coverage,
\begin{align*}
    \mathbb{P}(\psi(Z) \in [W_{(a)}, W_{(B-b)}])  \geq &  1- \frac{a+b +1}{B+1} - \widetilde d_{\mathrm{KS}}(\overline F(Z), U(0,1)) - \mathbb{E}_Z\left[ B\Delta(Z) \min \left\{1, \frac{1}{B \overline F(Z) (1 - \overline F(Z))} \right\} \right] \\
    \geq & 1- \frac{a+b +1}{B+1} - \widetilde d_{\mathrm{KS}}(\overline F(Z), U(0,1)) - \left(\sum_{i = 1}^B \kappa_i^2\right)\mathbb{E}_Z\left[  \min \left\{1, \frac{1}{B \overline F(Z) (1 - \overline F(Z))} \right\} \right] \\
    = & 1- \frac{a+b +1}{B+1} - \widetilde d_{\mathrm{KS}}(\overline F(Z), U(0,1)) - \left(\sum_{i = 1}^B \kappa_i^2\right)\left[ \mathbf I + \mathbf{II}\right] \\
    \geq & 1- \frac{a+b +1}{B+1} - \widetilde d_{\mathrm{KS}}(\overline F(Z), U(0,1)) - \left(\sum_{i = 1}^B \kappa_i^2\right) [I_B + \widetilde d_{\mathrm{KS}}(\overline F(Z), U(0,1))] .
\end{align*}

Applying the same technique, we obtain the coverage upper bound,
\begin{align*}
    \mathbb{P}(\psi(Z) \in [W_{(a)}, W_{(B-b)}]) \leq &  1- \frac{a+b }{B+1} + \widetilde d_{\mathrm{KS}}(\overline F(Z), U(0,1)) +\left(\sum_{i = 1}^B \kappa_i^2\right)\left[I_B + \widetilde d_{\mathrm{KS}}(\overline F(Z), U(0,1))\right].
\end{align*}
This completes the proof of the theorem. 
\section{Proof of Theorem~\ref{thm:finer_main_thm_ext}}
\label{app:proof_mainthm_finer}
We observe that, 
\begin{equation*}
    \begin{split}
        \mathbb{P}\left(\psi(Z) \in [W_{(a)}, W_{(B-b)}] \right) =& \mathbb{E}_Z \left[ \mathbb{P}\left( \psi(Z) \in [W_{(a)}, W_{(B-b)}] \mid Z \right) \right]\\
        = & \mathbb E_Z[\mathbb P(\psi(Z) \leq W_{(B-b)}|Z)] -  \mathbb E_Z[\mathbb P(\psi(Z) < W_{(a)}|Z)].
    \end{split}
\end{equation*}
As in the proof of \Cref{thm:main_result_extension}, we consider the following definition,
\[
S_B(Z) := \sum_{i=1}^B \mathbf{1}(W_i \leq \psi(Z)) \sim \text{Poi-Bin}(B, F_1(Z),\ldots,F_B(Z)) \quad \mbox{conditioned on }Z.
\]
Suppose $\overline S_B(Z) \sim \mbox{Bin}(B, \overline F(Z))$ conditioned on $Z$, where $\overline F(z) = (1/B) \sum_{i = 1}^B F_i(z)$. We will use the following stochastic ordering property between $\mbox{Poi-Bin}(B, p_1,\ldots,p_B)$ and $\mbox{Bin}(B, \overline p_B)$ (where $\overline p_B = (1/B) \sum_{i = 1}^B p_i$) from \cite{hoeffding1956distribution}, \cite{tang2023poisson},
\begin{equation}
\label{eq:stoch_ord}
    \begin{split}
        \mathbb P(\mbox{Poi-Bin}(B, p_1,\ldots,p_B) \leq k) \leq & \mathbb P(\mbox{Bin}(B, \overline p_B) \leq k) \quad \mbox{for} \quad 0 \leq k \leq B\overline p_B -1, \\
        \mathbb P(\mbox{Poi-Bin}(B, p_1,\ldots,p_B) \leq k) \geq & \mathbb P(\mbox{Bin}(B, \overline p_B) \leq k) \quad \mbox{for} \quad B\overline p_B \leq k \leq B. 
    \end{split}
\end{equation}
We have the following bounds using the above stochastic ordering properties, 
\begin{equation*}
    \begin{split}
\mathbb E_Z[\mathbb P(\psi(Z) < W_{(a)}|Z)] = & \mathbb E_Z[\mathbb P(S_B(Z) \leq a - 1|Z)]\\
=& \mathbb E_Z[\mathbb P(S_B(Z) \leq a - 1|Z)\textbf{1}\{a \leq B \overline F(Z) \}] + \mathbb E_Z[\mathbb P(S_B(Z) \leq a - 1|Z)\textbf{1}\{a > B \overline F(Z) \}] \\
\stackrel{(i)}{\leq} & \mathbb E_Z[\mathbb P(\overline S_B(Z) \leq a - 1|Z)\textbf{1}\{a \leq B \overline F(Z) \}] + \mathbb E_Z[\textbf{1}\{a > B \overline F(Z) \}] \\
= & \mathbb E_Z[\mathbb P(\overline S_B(Z) \leq a - 1|Z)] - \mathbb E_Z[\mathbb P(\overline S_B(Z) \leq a - 1|Z)\textbf{1}\{a > B \overline F(Z) \}] \\
&+\mathbb  E_Z[\textbf{1}\{a > B \overline F(Z) \}] \\
= &  \mathbb E_Z[\mathbb P(\overline S_B(Z) \leq a - 1|Z)] + \mathbb E_Z[\mathbb P(\overline S_B(Z) > a - 1|Z)\textbf{1}\{a > B \overline F(Z) \}]. \\
 \mathbb E_Z[\mathbb P(\psi(Z) \leq W_{(B-b)}|Z)] \geq & \mathbb E_Z[\mathbb P(S_B(Z) \leq B-b - 1|Z)]\\
=& \mathbb E_Z[\mathbb P(S_B(Z) \leq B-b - 1|Z)\textbf{1}\{B-b-1 \geq B \overline F(Z) \}] \\
&+ \mathbb E_Z[\mathbb P(S_B(Z) \leq B-b - 1|Z)\textbf{1}\{B-b-1 < B \overline F(Z) \}] \\
\stackrel{(ii)}{\geq} & \mathbb E_Z[\mathbb P(\overline S_B(Z) \leq B-b - 1|Z)\textbf{1}\{B-b-1 \geq B \overline F(Z) \}] \\
= & \mathbb E_Z[\mathbb P(\overline S_B(Z) \leq B-b - 1|Z)] - \mathbb E_Z[\mathbb P(\overline S_B(Z) \leq B-b - 1|Z)\textbf{1}\{B-b-1 < B \overline F(Z) \}] .
    \end{split}
\end{equation*}
The steps $(i), (ii)$ follow from \eqref{eq:stoch_ord}. We define the following for $0 \leq a, b < B/2$ , 
\begin{equation*}
    \begin{split}
        L_{\overline F(Z)}(a) = & \mathbb E[\mathbb P(\mbox{Bin}(B, \overline F(Z)) > a-1|\overline F(Z)) \textbf{1}\{a > B \overline F(Z)\}], \\
        R_{\overline F(Z)}(b) = & \mathbb E[\mathbb P(\mbox{Bin}(B, \overline F(Z)) \leq B-b-1|\overline F(Z)) \textbf{1}\{B-b-1 < B \overline F(Z)\}].
    \end{split}
\end{equation*}
Using these definitions we have the following lower bound on coverage,
\begin{equation*}
    \begin{split}
         \mathbb{P}\left(\psi(Z) \in [W_{(a)}, W_{(B-b)}] \right)  = & \mathbb E_Z[\mathbb P(\psi(Z) \leq W_{(B-b)}|Z)] -  \mathbb E_Z[\mathbb P(\psi(Z) < W_{(a)}|Z)] \\
         \geq & \mathbb E_Z[\mathbb P(\overline S_B(Z) \leq B-b - 1|Z)]  - \mathbb E_Z[\mathbb P(\overline S_B(Z) \leq a - 1|Z)]  - \left\{  L_{\overline F(Z)}(a) +  R_{\overline F(Z)}(b)\right\} \\
       \stackrel{(i)}{\geq} & 1 - \frac{a+b+1}{B+1} - \widetilde d_{\mathrm{KS}}(\overline F(Z), U(0,1)) -  \left\{  L_{\overline F(Z)}(a) +  R_{\overline F(Z)}(b)\right\}.  
    \end{split}
\end{equation*}
The inequality in step-$(i)$ follows by re-tracing the proof of \Cref{thm:main_result}. The only thing left now is to bound $ L_{\overline F(Z)}(a) +  R_{\overline F(Z)}(b)$. We observe the following, 
\begin{equation*}
    \begin{split}
          R_{\overline F(Z)}(b) = & \mathbb E[\mathbb P(\mbox{Bin}(B, \overline F(Z)) \leq B-b-1|\overline F(Z)) \textbf{1}\{B-b-1 < B \overline F(Z)\}] \\
          = & \mathbb E[\mathbb P(\mbox{Bin}(B, 1-\overline F(Z)) \geq b+1|\overline F(Z)) \textbf{1}\{b+1 > B(1 -\overline F(Z))  \}] \\
          = & L_{1 -\overline F(Z)}(b+1). 
    \end{split}
\end{equation*}
Moreover for any $0 \leq m < B/2$ we can bound the gap between $L_{\overline F(Z)}(m)$ and $L_U(m)$ where $U \stackrel{d}{=}\mbox{Uniform}(0,1)$. We define $g_m(w) = \mathbb P(\mbox{Bin}(B, w)\geq m )\textbf{1}\{m > B w \}$. Let $P_{\overline F(Z)}(\cdot)$ and $P_{U(0,1)}(\cdot)$ be the distribution function of $\overline F(Z), U(0,1)$ respectively. Therefore we have for any $0 \leq m < B/2$, 
\begin{equation*}
    \begin{split}
        |L_{\overline F(Z)}(m) -  L_U(m)| = &\left| \int_0^1 g_m d(P_{\overline F(Z)} - P_{U(0,1)}) \right| \\
        \stackrel{(i)}{=} & \left| \int_0^1 (P_{\overline F(Z)} - P_{U(0,1)}) dg_m\right| \\
        \leq & d_{\mathrm{KS}}(\overline F(Z), U(0,1))  \int_0^1 |dg_m | \\
        \stackrel{(ii)}{\leq} & 2 d_{\mathrm{KS}}(\overline F(Z), U(0,1)). 
    \end{split}
\end{equation*}
Since $g_m(w)$ is a function of bounded variation and vanishes outside $[0,1]$, Stieltjes integration by parts gives us $\int_0^1 g_m d(P_{\overline F(Z)} - P_{U(0,1)}) = - \int_0^1 (P_{\overline F(Z)} - P_{U(0,1)}) dg_m$. This justifies the equality in step $(i)$ of the above derivation. We note that the function $g_m(w)$ is increasing on $(-\infty, m/B)$, then has one downward jump at $w = m/B$ and is $0$ on $[m/B, \infty)$. Therefore we have $\mathrm{\mathrm{TV}}(g_m) = \int_0^1 |dg_m| \leq 1 + 1 = 2$ (inequality $(ii)$) as the increase of the function $g_m(w)$ before $m/B$ is at most $1$, and the downward jump at $m/B $ is also at most $1$. Combining these ideas we have, 
\begin{equation*}
    \begin{split}
          L_{\overline F(Z)}(a) +  R_{\overline F(Z)}(b) = & L_{\overline F(Z)}(a) +  L_{1 -\overline F(Z)}(b+1) \\
          \leq & L_U(a) + 2  d_{\mathrm{KS}}(\overline F(Z), U(0,1)) + L_{1-U}(b+1) +  2d_{\mathrm{KS}}(1-\overline F(Z), 1 - U(0,1)) \\
          = &  L_U(a) + L_U(b+1) + 4  d_{\mathrm{KS}}(\overline F(Z), U(0,1)).
    \end{split}
\end{equation*}
The last equality follows from $U \stackrel{d}{=}1 - U$ and $ d_{\mathrm{KS}}(\overline F(Z), U(0,1)) =  d_{\mathrm{KS}}(1-\overline F(Z), 1-U(0,1))$. We can bound $L_U(m)$ for any $0 \leq m < B/2$,
\begin{equation*}
    \begin{split}
        L_U(m) = & \mathbb E[\mathbb P(\mbox{Bin}(B, U) \geq m|U) \textbf{1}\{m > B U\}] \\
        = & \int_0^1 \mathbb P(\mbox{Bin}(B, u) \geq m)\textbf{1}\{m > B u\} du \\
        = & \int_0^{m/B} \mathbb P(\mbox{Bin}(B, u) \geq m) du \\
        \leq & \frac{1}{2} \int_0^{m/B}  du \\
        = & \frac{m}{2B}. 
    \end{split}
\end{equation*}
Therefore we have the following upper bound, 
\begin{equation*}
    \begin{split}
        L_{\overline F(Z)}(a) +  R_{\overline F(Z)}(b) \leq &  L_U(a) + L_U(b+1) + 4  d_{\mathrm{KS}}(\overline F(Z), U(0,1)) \\
        \leq & \frac{a+ b+1}{2B} + 4  d_{\mathrm{KS}}(\overline F(Z), U(0,1)) . 
    \end{split}
\end{equation*}
Hence we get the lower bound on coverage, 
\begin{equation*}
    \begin{split}
        \mathbb{P}\left(\psi(Z) \in [W_{(a)}, W_{(B-b)}] \right) \geq &  1 - \frac{a+b+1}{B+1} - \widetilde d_{\mathrm{KS}}(\overline F(Z), U(0,1)) -  \left\{  L_{\overline F(Z)}(a) +  R_{\overline F(Z)}(b)\right\} \\
        \geq & 1 - \frac{3(a+b+1)}{2B} - \widetilde d_{\mathrm{KS}}(\overline F(Z), U(0,1))- 4  d_{\mathrm{KS}}(\overline F(Z), U(0,1)) \\
        \geq & 1 - \frac{3(a+b+1)}{2B} - 6 d_{\mathrm{KS}}(\overline F(Z), U(0,1)) . 
    \end{split}
\end{equation*}
This completes the proof of the theorem. 
\section{Proof of Theorem~\ref{thm:beyond_exchangeability}}
\label{app:proof_beyond_exchangeability}
We consider the following definitions for any $\alpha \in (0,1)$,
\begin{equation*}
    \begin{split}
        \widehat r_{\alpha}^R =& \inf \left\{ \theta \in \mathbb R: \frac{1}{B+1}\left[ \sum_{i = 1}^B \textbf{1}\{W _i \leq \theta\}  \right] \geq \alpha \right\} =W_{(\lceil (B+1)\alpha \rceil )}, \\
        \widehat r_{\alpha}^L =& \sup \left\{ \theta \in \mathbb R: \frac{1}{B+1}\left[1 + \sum_{i = 1}^B \textbf{1}\{W _i < \theta\}  \right] \leq \alpha \right\} =  W_{(\lfloor (B+1)\alpha \rfloor -1)}, \\
        r_{\alpha}^R = &\inf \left\{ \theta \in \mathbb R: \frac{1}{B+1}\left[\textbf{1}\{\psi(Z) \leq \theta\} + \sum_{i = 1}^B \textbf{1}\{W _i \leq \theta\}  \right] \geq \alpha \right\}, \\
          r_{\alpha}^L =& \sup \left\{ \theta \in \mathbb R: \frac{1}{B+1}\left[\textbf{1}\{\psi(Z) < \theta\} + \sum_{i = 1}^B \textbf{1}\{W _i < \theta\}  \right] \leq \alpha \right\}.
    \end{split}
\end{equation*}
We note that the objective function $\textbf{1}\{\psi(Z) \leq \theta\} + \sum_{i = 1}^B \textbf{1}\{W _i \leq \theta\} $ is right continuous in $\theta$ and the objective function $\textbf{1}\{\psi(Z) < \theta\} + \sum_{i = 1}^B \textbf{1}\{W _i < \theta\}  $ is left continuous in $\theta$. Therefore by definition, $r_{\alpha}^R, r_{\alpha}^L$ satisfy the following inequalities,
\begin{equation}
\label{eq:r_alpha_conseq}
\begin{split}
    \frac{1}{B+1}\left[\textbf{1}\{\psi(Z) \leq r_{\alpha}^R\} + \sum_{i = 1}^B \textbf{1}\{W _i \leq r_{\alpha}^R\}  \right] \geq & \alpha,  \\
    \frac{1}{B+1}\left[\textbf{1}\{\psi(Z) < r_{\alpha}^L\} + \sum_{i = 1}^B \textbf{1}\{W _i < r_{\alpha}^L\}  \right] \leq & \alpha.
\end{split}
\end{equation}
We observe that $\widehat r_{\alpha}^R$ belongs in the set whose infimum is $r_{\alpha}^R$,
\begin{equation*}
    \begin{split}
        & \frac{1}{B+1}\left[\textbf{1}\{\psi(Z) \leq \widehat r_{\alpha}^R\} + \sum_{i = 1}^B \textbf{1}\{W _i \leq \widehat r_{\alpha}^R\}  \right] \\
        \geq & \frac{1}{B+1}\left[\sum_{i = 1}^B \textbf{1}\{W _i \leq \widehat r_{\alpha}^R\}  \right] \\
        \geq & \alpha. 
    \end{split}
\end{equation*}
Therefore $r_{\alpha}^R \leq \widehat r_{\alpha}^R$ for any $\alpha \in [0,1]$. We can similarly show that $r_{\alpha}^L \geq \widehat r_{\alpha}^L$ for any $\alpha \in [0,1]$. We observe the following for any $0< \gamma, \beta < 1$,
\begin{equation*}
    \begin{split}
         &\mathbb P(\psi(Z) \in [\widehat r_{\gamma/2}^L , \widehat r_{1 - (\beta/2)}^R] ) \\
         \geq & \mathbb P(\psi(Z) \in [r_{\gamma/2}^L, r_{1 - (\beta/2)}^R] ) \\
         = & 1 - \mathbb P(\psi(Z) < r_{\gamma/2}^L) - \mathbb P(\psi(Z) > r_{1 - (\beta/2)}^R) \\
         = &  \mathbb P(\psi(Z) \leq r_{1 - (\beta/2)}^R) - \mathbb P(\psi(Z) < r_{\gamma/2}^L)\\
         = & \mathbb E \left[ \frac{1}{B+1}\left[\textbf{1}\{\psi(Z) \leq r_{1 - (\beta/2)}^R\} + \sum_{i = 1}^B \textbf{1}\{W _i \leq r_{1 - (\beta/2)}^R\}  \right] \right] - \mathbb E \left[  \frac{1}{B+1}\left[\textbf{1}\{\psi(Z) < r_{\gamma/2}^L\} + \sum_{i = 1}^B \textbf{1}\{W _i < r_{\gamma/2}^L\}  \right]\right] \\
         &+ \left(\mathbb P(\psi(Z) \leq r_{1 - (\beta/2)}^R) - \frac{1}{B+1}\left[\mathbb P(\psi(Z) \leq r_{1 - (\beta/2)}^R) + \sum_{i = 1}^B \mathbb P(W_i \leq r_{1 - (\beta/2)}^R)    \right] \right) \\
         &- \left(\mathbb P(\psi(Z) < r_{\gamma/2}^L) - \frac{1}{B+1}\left[\mathbb P(\psi(Z) < r_{\gamma/2}^L) + \sum_{i = 1}^B \mathbb P(W_i < r_{\gamma/2}^L)    \right] \right) \\
         \stackrel{\eqref{eq:r_alpha_conseq}}{\geq} & 1 - (\beta/2) - (\gamma/2) \\
         &+ \left(\mathbb P(\psi(Z) \in [r_{\gamma/2}^L, r_{1 - (\beta/2)}^R]) - \frac{1}{B+1}\left[\mathbb P(\psi(Z) \in [r_{\gamma/2}^L, r_{1 - (\beta/2)}^R]) + \sum_{i = 1}^B \mathbb P(W_i \in [r_{\gamma/2}^L, r_{1 - (\beta/2)}^R])    \right] \right) \\
         \geq & 1 - (\beta/2) - (\gamma/2) - \Gamma(W_1,\ldots,W_B,\psi(Z)).
    \end{split}
\end{equation*}
This completes the proof of the lower bound on the coverage. For obtaining the upper bound on coverage we consider the following definitions for any $\alpha \in (0,1)$,
\begin{equation*}
    \begin{split}
        \widetilde r_{\alpha}^R = &\inf \left\{ \theta \in \mathbb R: \frac{1}{B+1}\left[\textbf{1}\{\psi(Z) \leq \theta\} + \sum_{i = 1}^B \textbf{1}\{W _i \leq \theta\}  \right] -\frac{1}{B+1} \geq \alpha \right\}, \\
         \widetilde  r_{\alpha}^L =& \sup \left\{ \theta \in \mathbb R: \frac{1}{B+1}\left[\textbf{1}\{\psi(Z) < \theta\} + \sum_{i = 1}^B \textbf{1}\{W _i < \theta\}  \right]+ \frac{1}{B+1} \leq \alpha \right\}.
    \end{split}
\end{equation*}
Since the objective function $\textbf{1}\{\psi(Z) \leq \theta\} + \sum_{i = 1}^B \textbf{1}\{W _i \leq \theta\} $ is right continuous in $\theta$ and the objective function $\textbf{1}\{\psi(Z) < \theta\} + \sum_{i = 1}^B \textbf{1}\{W _i < \theta\}  $ is left continuous in $\theta$, $r_{\alpha}^R, r_{\alpha}^L$ satisfy the following inequalities,
\begin{equation}
\label{eq:r_alpha_conseq_upper}
\begin{split}
    \frac{1}{B+1}\left[\textbf{1}\{\psi(Z) \leq \widetilde r_{\alpha}^R\} + \sum_{i = 1}^B \textbf{1}\{W _i \leq \widetilde r_{\alpha}^R\}  \right] -\frac{1}{B+1}\geq & \alpha,  \\
    \frac{1}{B+1}\left[\textbf{1}\{\psi(Z) < \widetilde r_{\alpha}^L\} + \sum_{i = 1}^B \textbf{1}\{W _i < \widetilde r_{\alpha}^L\}  \right] +\frac{1}{B+1} \leq & \alpha.
\end{split}
\end{equation}
Note that $\widetilde r_{\alpha}^R$ belongs to the set whose infimum is $\widehat r_{\alpha}^R$,
\begin{equation*}
    \begin{split}
       & \frac{1}{B+1} \left[\sum_{i = 1}^B \textbf{1}\{W _i \leq \widetilde r_{\alpha}^R\}  \right] \\
        = & \frac{1}{B+1} \left[\textbf{1}\{\psi(Z) \leq \widetilde r_{\alpha}^R\} + \sum_{i = 1}^B \textbf{1}\{W _i \leq \widetilde r_{\alpha}^R\}  \right] - \frac{\textbf{1}\{\psi(Z) \leq \widetilde r_{\alpha}^R\}}{B+1} \\
        \geq & \frac{1}{B+1} \left[\textbf{1}\{\psi(Z) \leq \widetilde r_{\alpha}^R\} + \sum_{i = 1}^B \textbf{1}\{W _i \leq \widetilde r_{\alpha}^R\}  \right] - \frac{1}{B+1}  \\
        \geq & \alpha. 
    \end{split}
\end{equation*}
Therefore we have $\widehat r_{\alpha}^R \leq \widetilde r_{\alpha}^R$ for any $\alpha \in (0,1)$. We can similarly show that $\widehat r_{\alpha}^L > \widetilde r_{\alpha}^L$ for any $\alpha \in (0,1)$. We now observe the following, 
\begin{equation}
\label{eq:beyond_exchange_upper_R}
    \begin{split}
     &\frac{1}{B+1} \left[\textbf{1}\{\psi(Z) \leq \widetilde r_{\alpha}^R\} + \sum_{i = 1}^B \textbf{1}\{W _i \leq \widetilde r_{\alpha}^R\}  \right] \\\
     = & \alpha + \frac{1}{B+1} + \left\{ \frac{1}{B+1} \left[\textbf{1}\{\psi(Z) \leq \widetilde r_{\alpha}^R\} + \sum_{i = 1}^B \textbf{1}\{W _i \leq \widetilde r_{\alpha}^R\}  \right] -\frac{1}{B+1} - \alpha  \right\}\\
     \stackrel{(i)}{\leq} & \alpha + \frac{1}{B+1} + \frac{1}{B+1} \\
     = & \alpha + \frac{2}{B+1}. 
    \end{split}
\end{equation}
The step-$(i)$ follows from the Jump Lemma \cite{angelopoulos2022conformal} as the distributions of $\psi(Z), \{W_i\}_{i = 1}^B$ are continuous. We can similarly show that, 
\begin{equation}
    \label{eq:beyond_exchange_upper_L}
\frac{1}{B+1} \left[\textbf{1}\{\psi(Z) \leq \widetilde r_{\alpha}^L\} + \sum_{i = 1}^B \textbf{1}\{W _i \leq \widetilde r_{\alpha}^L\}  \right] \geq \alpha - \frac{2}{B+1}. 
\end{equation}
Therefore we have the following upper bound on coverage for any $0< \gamma, \beta <1$,
\begin{equation*}
    \begin{split}
         \mathbb P(\psi(Z) \in [\widehat r_{\gamma/2}^L , \widehat r_{1 - (\beta/2)}^R] )  \leq & \mathbb P(\psi(Z) \in [\widetilde r_{\gamma/2}^L, \widetilde r_{1 - (\beta/2)}^R] ) \\
         = & 1 - \mathbb P(\psi(Z) <\widetilde r_{\gamma/2}^L) - \mathbb P(\psi(Z) >\widetilde r_{1 - (\beta/2)}^R) \\
         = &  \mathbb P(\psi(Z) \leq \widetilde r_{1 - (\beta/2)}^R) - \mathbb P(\psi(Z) < \widetilde r_{\gamma/2}^L).
    \end{split}
\end{equation*}
We add and subtract some terms to get,
\begin{equation*}
    \begin{split}
     &\mathbb P(\psi(Z) \in [\widehat r_{\gamma/2}^L , \widehat r_{1 - (\beta/2)}^R] ) \\
 = & \mathbb E \left[ \frac{1}{B+1}\left[\textbf{1}\{\psi(Z) \leq \widetilde r_{1 - (\beta/2)}^R\} + \sum_{i = 1}^B \textbf{1}\{W _i \leq \widetilde r_{1 - (\beta/2)}^R\}  \right] \right] \\
         &- \mathbb E \left[  \frac{1}{B+1}\left[\textbf{1}\{\psi(Z) < \widetilde r_{\gamma/2}^L\} + \sum_{i = 1}^B \textbf{1}\{W _i < \widetilde r_{\gamma/2}^L\}  \right]\right] \\
         &+ \left(\mathbb P(\psi(Z) \leq \widetilde r_{1 - (\beta/2)}^R) - \frac{1}{B+1}\left[\mathbb P(\psi(Z) \leq \widetilde r_{1 - (\beta/2)}^R) + \sum_{i = 1}^B \mathbb P(W_i \leq \widetilde r_{1 - (\beta/2)}^R)    \right] \right) \\
         &- \left(\mathbb P(\psi(Z) < \widetilde r_{\gamma/2}^L) - \frac{1}{B+1}\left[\mathbb P(\psi(Z) < \widetilde r_{\gamma/2}^L) + \sum_{i = 1}^B \mathbb P(W_i < \widetilde r_{\gamma/2}^L)    \right] \right) \\
         \stackrel{\eqref{eq:beyond_exchange_upper_R}, \eqref{eq:beyond_exchange_upper_L}}{\leq} & 1 - (\beta/2) - (\gamma/2) + \frac{4}{B+1}\\
         &+ \left(\mathbb P(\psi(Z) \in [\widetilde r_{\gamma/2}^L, \widetilde r_{1 - (\beta/2)}^R]) - \frac{1}{B+1}\left[\mathbb P(\psi(Z) \in [\widetilde r_{\gamma/2}^L, \widetilde r_{1 - (\beta/2)}^R]) + \sum_{i = 1}^B \mathbb P(W_i \in [\widetilde r_{\gamma/2}^L, \widetilde r_{1 - (\beta/2)}^R])    \right] \right) \\
         \leq & 1 - (\beta/2) - (\gamma/2) + \frac{4}{B+1}+  \Gamma(W_1,\ldots,W_B,\psi(Z)).
    \end{split}
\end{equation*}
This completes the proof of the theorem. 
\section{Proof of Corollary~\ref{cor:boot}}
\label{app:proof_cor_boot}
We will apply the second statement of \Cref{thm:main_result} to prove the coverage guarantee of $\mathrm{CI}^{\mathtt{mod-boot}}_{m,B,\alpha} = \{\theta: S_m(\theta) \in [W_{(\lfloor (B+1)(\alpha/2)\rfloor)},  W_{(\lceil (B+1)(1 - \alpha)\rceil + \lfloor (B+1)(\alpha/2)\rfloor)}) \}$. Here $D_m$ plays the role of $Z$, $S(\tau_m(\widehat \theta_m - \theta_0 ))= \psi(D_m)$ and $\{S(\tau_m(\widehat \theta_m^{*b} - \widehat \theta_m))\}_{b = 1}^B = \{W_b\}_{b = 1}^B$ as in \Cref{thm:main_result}. We note that conditioned on $D_m$, the centered bootstrap estimates $\{W_b\}_{b = 1}^B$ are independent and identically distributed and hence satisfy the conditions of \Cref{thm:main_result}. Moreover we have, 
\[
F_0(D_m) = F_0(Z) =  \mathbb P(W_1 \leq \psi(Z) |Z) = \mathbb P(S(\tau_m(\widehat \theta_m^{*b} - \widehat \theta_m))\leq S(\tau_m(\widehat \theta_m - \theta_0)) | D_m). 
\]
We apply \Cref{thm:main_result} with $a = \lfloor (B+1)(\alpha/2)\rfloor$ and $b = B - \lceil (B+1)(1 - \alpha)\rceil - \lfloor (B+1)(\alpha/2)\rfloor$,
\begin{equation*}
\begin{split}
  &  1 - \frac{\lfloor (B+1)(\alpha/2)\rfloor +  B+1 -  \lceil (B+1)(1 - \alpha)\rceil - \lfloor (B+1)(\alpha/2)\rfloor}{B+1} - \widetilde d_{\mathrm{KS}}(F_0(D_m), U(0,1)) \\
  \leq & \mathbb P \left(S_m(\theta_0) \in [W_{(\lfloor (B+1)(\alpha/2)\rfloor)},  W_{(\lceil (B+1)(1 - \alpha)\rceil + \lfloor (B+1)(\alpha/2)\rfloor)}) \right) \\
  \leq &  1 - \frac{\lfloor (B+1)(\alpha/2)\rfloor +  B+1 -  \lceil (B+1)(1 - \alpha)\rceil - \lfloor (B+1)(\alpha/2)\rfloor}{B+1}  + \widetilde d_{\mathrm{KS}}(F_0(D_m), U(0,1)). 
\end{split} 
\end{equation*}
It is easy to check that, 
\[
1 - \alpha \leq  1 - \frac{ B+1 -  \lceil (B+1)(1 - \alpha)\rceil }{B+1} \leq 1 - \alpha+ \frac{1}{B+1}.
\]
Therefore we have the following coverage guarantee,
\[
- \widetilde d_{\mathrm{KS}}(F_0(D_m), U(0,1)) \leq \mathbb P( \theta_0 \in \mathrm{CI}^{\mathtt{mod-boot}}_{m,B,\alpha} ) - (1 - \alpha) \leq \frac{1}{B+1} +  \widetilde d_{\mathrm{KS}}(F_0(D_m), U(0,1)) . 
\]
This completes the proof of the corollary. 

\section{Proof of \Cref{rem:random_boot}}
\label{app:proof_rand_boot}
In this section, we prove the lower bound on the coverage of the randomized modified confidence interval $\mathrm{CI}^{\mathtt{rand-mod-boot}}_{m,B,\alpha}$. The upper bound can be proved analogously. Let us define $\overline {\mathrm{CI}}^{\mathtt{mod-boot}}_{m,B,\alpha} =  \{\theta: S_m(\theta) \in [W_{(\lfloor (B+1)(\alpha/2)\rfloor)},  W_{(\lfloor (B+1)(1 - \alpha)\rfloor + \lfloor (B+1)(\alpha/2)\rfloor)}) \}$. We note that $\mathrm{CI}^{\mathtt{rand-mod-boot}}_{m,B,\alpha}$ is defined as follows ($U$ is an independently generated Uniform $(0,1)$ random variable), 
 \begin{equation*}
  \begin{split}
      \mathrm{CI}^{\mathtt{rand-mod-boot}}_{m,B,\alpha} =&\begin{cases} \mathrm{CI}^{\mathtt{mod-boot}}_{m,B,\alpha} \quad & \mbox{if} \quad U \leq \tau_{\alpha, B} \\
     \overline {\mathrm{CI}}^{\mathtt{mod-boot}}_{m,B,\alpha}  \quad & \mbox{if} \quad U > \tau_{\alpha, B}
      \end{cases} \quad \mbox{where}, \\
  \tau_{\alpha, B} =& \begin{cases}
        1 \quad & \mbox{if} \quad (B+1)(1-\alpha) \in \mathbb N, \\
        \frac{(1-\alpha) - (\lfloor (B+1)(1-\alpha) \rfloor/ (B+1))}{ (\lceil (B+1)(1-\alpha) \rceil - \lfloor (B+1)(1-\alpha) \rfloor)/ (B+1)} \quad  &\mbox{otherwise}.
    \end{cases}
    \end{split}
  \end{equation*}
  An application of \Cref{thm:main_result} to the setup of non-parametric bootstrap yields the following coverage guarantees,
  \begin{equation*}
      \begin{split}
          \mathbb P(\theta_0 \in \mathrm{CI}^{\mathtt{mod-boot}}_{m,B,\alpha} ) \geq & \frac{\lceil (B+1)(1-\alpha) \rceil}{B+1} -  \widetilde d_{\mathrm{KS}}(F_0(D_m), U(0,1)), \\
          \mathbb P(\theta_0 \in \overline {\mathrm{CI}}^{\mathtt{mod-boot}}_{m,B,\alpha}  ) \geq & \frac{\lfloor (B+1)(1-\alpha) \rfloor}{B+1} -  \widetilde d_{\mathrm{KS}}(F_0(D_m), U(0,1)).
      \end{split}
  \end{equation*}
We use these bounds to obtain the lower bound on the coverage of the randomized modified confidence interval $\mathrm{CI}^{\mathtt{rand-mod-boot}}_{m,B,\alpha}$,
\begin{equation*}
    \begin{split}
     &\mathbb P(\theta_0 \in \mathrm{CI}^{\mathtt{rand-mod-boot}}_{m,B,\alpha} ) \\
     = &      \mathbb P(\theta_0 \in \mathrm{CI}^{\mathtt{rand-mod-boot}}_{m,B,\alpha} | U \leq \tau_{\alpha,B} )\mathbb P(U \leq \tau_{\alpha,B}) +  \mathbb P(\theta_0 \in \mathrm{CI}^{\mathtt{rand-mod-boot}}_{m,B,\alpha} | U > \tau_{\alpha,B}) \mathbb P(U > \tau_{\alpha,B}) \\
     = & \mathbb P(\theta_0 \in \mathrm{CI}^{\mathtt{mod-boot}}_{m,B,\alpha} ) \tau_{\alpha, B} + \mathbb P(\theta_0 \in \overline {\mathrm{CI}}^{\mathtt{mod-boot}}_{m,B,\alpha}  ) (1 - \tau_{\alpha,B}) \\
     \geq & \frac{ \tau_{\alpha,B}\lceil (B+1)(1-\alpha) \rceil}{B+1}  + \frac{(1 - \tau_{\alpha,B})\lfloor (B+1)(1-\alpha) \rfloor}{B+1}  -  \widetilde d_{\mathrm{KS}}(F_0(D_m), U(0,1))   \\
     = & 1 - \alpha - \widetilde d_{\mathrm{KS}}(F_0(D_m), U(0,1))  . 
    \end{split}
\end{equation*}
This completes the proof of the lower bound on the coverage. We can similarly prove that the coverage of the randomized confidence interval is upper bounded by $1 - \alpha + \widetilde d_{\mathrm{KS}}(F_0(D_m), U(0,1)) $. 
\section{Proof of \Cref{cor:permute}}
\label{app:proof_permute}
Let us first prove the calibration result for $B=|G|$. In this scenario, we have $\phi^{\mathtt{mod-permute}}_{m,B,\alpha}(X) =  \textbf{1}\{T(X) \geq T^*_{(\lceil B(1-\alpha) \rceil +2)} \}$. We will apply the second statement of \Cref{thm:main_result} to prove the type-I error control of $\phi^{\mathtt{mod-permute}}_{m,B,\alpha}(X)$. Here $D_m = \{X_1,\ldots,X_m\}$ plays the role of $Z$, $T(X) = \psi(D_m)$ and $\{T(X_{\sigma_b^*})\}_{b = 1}^B$ play the role of $(W_1,\ldots,W_B)$. We note that conditioned on $D_m$, the statistics based on permuted data $\{T(X_{\sigma_b^*})\}_{b = 1}^B$ are independent and identically distributed and hence satisfy the conditions of \Cref{thm:main_result}. Moreover we have, 
\[
F_0(D_m) = F_0(Z) =  \mathbb P(W_b \leq \psi(Z) |Z) = \mathbb P( T(X_{\sigma_b^*}) \leq T(X) | D_m). 
\]
We apply \Cref{thm:main_result} with $a =0$ and $b =\lfloor B\alpha \rfloor -2$,
\begin{equation*}
\begin{split}
  &  1 - \frac{ \lfloor B\alpha \rfloor - 1}{B+1} -  d_{\mathrm{KS}}(F_0(D_m), U(0,1)) \\
  \leq & \mathbb P_{H_0} \left(T(X) < T^*_{(\lceil B(1-\alpha) \rceil +2)} \right) \\
  \leq &  1 - \frac{\lfloor B\alpha \rfloor - 1}{B+1}   +  d_{\mathrm{KS}}(F_0(D_m), U(0,1)). 
\end{split} 
\end{equation*}
Note that since this is a one-sided confidence interval we have $ d_{\mathrm{KS}}(F_0(D_m), U(0,1))$ in the slack instead of $\widetilde  d_{\mathrm{KS}}(F_0(D_m), U(0,1))$. Since there are no ties almost surely, we define $\mbox{rank}(T(X)) \coloneqq \sum_{\sigma \in G} \textbf{1}\{T(X_{\sigma}) \leq T(X) \} $. We observe that under $H_0$, the {KS} distance $d_{\mathrm{KS}}(F_0(D_m), U(0,1))$ can be bounded as follows, 
\begin{equation*}
    \begin{split}
     d_{\mathrm{KS}}(F_0(D_m), U(0,1)) =& \sup_x\left|\mathbb P(F_0(D_m) \leq x) - x \right|   \\
     = & \sup_x\left|\mathbb P \left(\frac{1}{B}\sum_{\sigma \in G}\textbf{1}\{T(X_{\sigma}) \leq T(X) \} \leq x\right) -x \right| \\
     =& \sup_x\left|\mathbb P \left(\mbox{rank}(T(X)) \leq Bx\right) -x \right| \\
     \stackrel{(i)}{\leq} &\frac{1}{B}.
    \end{split}
\end{equation*}
The step $(i)$ follows as under the null hypotheses $H_0$, the statistics based on the permuted data $\{T(X_{\sigma})\}_{\sigma \in G}$ are exchangeable and there are no ties almost surely and consequently we have $\mbox{rank}(T(X)) \sim \mbox{Unif}\{1,\ldots,|G|\}$. Therefore we have the following coverage guarantee,
\[
1 - \frac{ \lfloor B\alpha \rfloor - 1}{B+1} - \frac{1}{B} \leq \mathbb P_{H_0} \left(T(X) < T^*_{(\lceil B(1-\alpha) \rceil +2)} \right) \leq 1 - \frac{ \lfloor B\alpha \rfloor - 1}{B+1} + \frac{1}{B}. 
\]
A straightforward simplification yields the following bounds,
\[
1-\alpha \leq \mathbb P_{H_0} \left(T(X) < T^*_{(\lceil B(1-\alpha) \rceil +2)} \right) \leq 1 - \alpha + \frac{\alpha + 2}{B+1} + \frac{1}{B}.
\]
This implies the following regarding the permutation testing rule $\phi^{\mathtt{mod-permute}}_{m,B,\alpha}(X)$,
\[
\alpha -\frac{\alpha+2}{B+1} -\frac{1}{B} \leq \mathbb E_{H_0}[\phi^{\mathtt{mod-permute}}_{m,B,\alpha}(X)] \leq \alpha.
\]
This completes the proof of the first part of the corollary. For $B < |G|$, we have $\phi^{\mathtt{mod-permute}}_{m,B,\alpha}(X) =  \textbf{1}\{T(X) \geq T^*_{(\lceil (B+1)(1-\alpha) \rceil +1)} \}$. We apply \Cref{thm:main_result} with $a =0$ and $b = B- (\lceil (B+1)(1-\alpha) \rceil -1$,
\begin{equation*}
\begin{split}
  &  1 - \frac{ B- (\lceil (B+1)(1-\alpha) \rceil}{B+1} -  d_{\mathrm{KS}}(F_0(D_m), U(0,1)) \\
  \leq & \mathbb P_{H_0} \left(T(X) < T^*_{(\lceil (B+1)(1-\alpha) \rceil +1)} \right) \\
  \leq &  1 - \frac{B- (\lceil (B+1)(1-\alpha) \rceil}{B+1}   +  d_{\mathrm{KS}}(F_0(D_m), U(0,1)). 
\end{split} 
\end{equation*}
We observe that under $H_0$, the {KS} distance $d_{\mathrm{KS}}(F_0(D_m), U(0,1))$ can be bounded as follows, 
\begin{equation*}
    \begin{split}
     d_{\mathrm{KS}}(F_0(D_m), U(0,1)) =& \sup_x\left|\mathbb P(F_0(D_m) \leq x) - x \right|   \\
     = & \sup_x\left|\mathbb P \left(\frac{1}{|G|}\sum_{\sigma \in G}\textbf{1}\{T(X_{\sigma}) \leq T(X) \} \leq x\right) -x \right| \\
     =& \sup_x\left|\mathbb P \left(\mbox{rank}(T(X) \leq |G|x\right) -x \right| \\
     \stackrel{(i)}{\leq} &\frac{1}{|G|} \\
     \leq & \frac{1}{B+1}.
    \end{split}
\end{equation*}
The step $(i)$ follows as under the null hypotheses $H_0$, the statistics based on the permuted data $\{T(X_{\sigma})\}_{\sigma \in G}$ are exchangeable and there are no ties almost surely and consequently we have $\mbox{rank}(T(X)) \sim \mbox{Unif}\{1,\ldots,|G|\}$. Therefore we have the following coverage guarantee,
\[
1 - \frac{ B- (\lceil (B+1)(1-\alpha) \rceil}{B+1} - \frac{1}{B+1} \leq \mathbb P_{H_0} \left(T(X) < T^*_{(\lceil (B+1)(1-\alpha) \rceil +1)} \right) \leq 1 - \frac{B- (\lceil (B+1)(1-\alpha) \rceil}{B+1} + \frac{1}{|G|}. 
\]
A straightforward simplification yields the following bounds,
\[
1-\alpha \leq \mathbb P_{H_0} \left(T(X) < T^*_{(\lceil (B+1)(1-\alpha) \rceil +1)} \right) \leq 1 - \alpha + \frac{ 2}{B+1} + \frac{1}{|G|}.
\]
This implies the following regarding the permutation testing rule $\phi^{\mathtt{mod-permute}}_{m,B,\alpha}(X)$,
\[
\alpha -\frac{2}{B+1} -\frac{1}{|G|} \leq \mathbb E_{H_0}[\phi^{\mathtt{mod-permute}}_{m,B,\alpha}(X)] \leq \alpha.
\]
This completes the proof of the corollary.

\section{Proof of \Cref{cor:indep_CP}}
\label{app:proof_indep_CP}
The first statement of the corollary follows from a direct application of \Cref{thm:main_result_extension} for one-sided confidence intervals of the form $(0, W_{(B-b)}]$. To see this, we note that $R_1,\ldots,R_m$ play the role of $W_1,\ldots,W_B$ ($B = m$) and $R_{m+1}$ plays the role of $\psi(Z) = Z$. Since $R_1,\ldots,R_{m+1}$ are independent, we have $F_i(z) = \mathbb P(W_i < \psi(Z) |Z=z) = \mathbb P(R_i < R_{m+1}|R_{m+1} = z) = \mathbb P(R_i <z)$ for all $i \in [m]$ and for all $z \in \mathbb R$. For one-sided confidence intervals, we can show that we can replace $\widetilde d_{\mathrm{KS}}(\overline F(Z),U(0,1))$ by the smaller slack $d_{\mathrm{KS}}(\overline F(Z),U(0,1))$. Therefore we get the following coverage bound by setting $a = 0, b = m - \lceil (m+1)(1-\alpha) \rceil$, 
\[
-\delta \leq \mathbb P( Y_{m+1} \in \widehat C_{m,\alpha}(X_{m+1})) -(1 - \alpha) \leq \frac{1}{m+1} + \delta \quad \mbox{where},
 \]
 \begin{equation*}
    \delta = d_{\mathrm{KS}}(\overline F(R_{m+1}), U(0,1)) +\left(\sum_{i = 1}^m \kappa_i^2\right)\left[I_m +  d_{\mathrm{KS}}(\overline F(R_{m+1}), U(0,1))\right].
\end{equation*}
In the above bound, $\kappa_i = d_{\mathrm{KS}}(F_i,F_{m+1})$ for $i \in [m]$ and $\overline F(z) = (1/m) \sum_{i = 1}^m F_i(z)$ for $z \in \mathbb R$. We use \Cref{thm:finer_main_thm_ext} to prove the lower bound on the coverage of the modified prediction interval $\widehat C_{m,\alpha}^{\mathtt{mod-pred}}(X_{m+1}) = \{y: R(X_{m+1}, y) \leq R_{(m+1 - \lfloor 2m \alpha/3 \rfloor )}\}$. For one-sided confidence intervals of the form $(0, W_{(B-b)}]$, the slack $6 d_{\mathrm{KS}}(\overline F(Z), U(0,1))$ can be improved to $3 d_{\mathrm{KS}}(\overline F(Z), U(0,1))$. Setting $a = 0, b = \lfloor 2m \alpha/3 \rfloor -1$ in \Cref{thm:finer_main_thm_ext} yields the following, 
\begin{equation*}
    \begin{split}
        \mathbb P(Y_{m+1} \in \widehat C_{m,\alpha}^{\mathtt{mod-pred}}(X_{m+1}) ) = & \mathbb P(R_{m+1} \in (0, R_{(m+1 - \lfloor 2m \alpha/3 \rfloor )}]) \\
        \geq & 1 - \frac{3\lfloor 2m \alpha/3 \rfloor}{2m} - 3 d_{\mathrm{KS}}(\overline F(R_{m+1}), U(0,1)) \\
        \geq & 1 - \alpha - 3 d_{\mathrm{KS}}(\overline F(R_{m+1}), U(0,1)).
     \end{split}
\end{equation*}
This completes the proof of the corollary. 
\section{Proof of \Cref{cor:indep_CP_same}}
\label{app:proof_indep_CP_same}
The corollary follows from a direct application of \Cref{thm:main_result} for one-sided confidence intervals of the form $(0, W_{(B-b)}]$. To see this, we note that $R_1,\ldots,R_m$ plays the role of $W_1,\ldots,W_B$ ($B = m$) and $R_{m+1}$ play the role of $\psi(Z) = Z$. Since $R_1,\ldots,R_{m+1}$ are independent, $\widetilde F_0(z) = \mathbb P(W_1 < \psi(Z) |Z = z) = \mathbb P(R_1 < R_{m+1}|R_{m+1} = z)= \mathbb P(R_1 < z)$. We apply the third statement of \Cref{thm:main_result} with $a = 0, b =m-  \lceil (m+1)(1 - \alpha) \rceil$,
\begin{equation*}
\begin{split}
  &  1 - \frac{m +1-  \lceil (m+1)(1 - \alpha) \rceil}{m+1} -  d_{\mathrm{KS}}(\widetilde F_0(R_{m+1}), U(0,1)) \\
  \leq &  \mathbb P(Y_{m+1} \in \widehat C_{m,\alpha}(X_{m+1}) ) \\
  \leq & 1 - \frac{m +1-  \lceil (m+1)(1 - \alpha) \rceil}{m+1} + d_{\mathrm{KS}}(\widetilde F_0(R_{m+1}), U(0,1)) .
\end{split} 
\end{equation*}
Note that for one-sided confidence intervals, we can show that we can replace $\widetilde d_{\mathrm{KS}}(\widetilde F_0(R_{m+1}),U(0,1))$ by the smaller slack $d_{\mathrm{KS}}(\widetilde F_0(R_{m+1}),U(0,1))$. It is easy to check that,
\[
1 - \alpha \leq 1 - \frac{m +1-  \lceil (m+1)(1 - \alpha) \rceil}{m+1}  \leq 1 - \alpha + \frac{1}{m+1}. 
\]
Therefore we have the following coverage guarantee,
\[
-d_{\mathrm{KS}}(\widetilde F_0(R_{m+1}), U(0,1)) \leq \mathbb P(Y_{m+1} \in \widehat C_{m,\alpha}(X_{m+1}) ) - (1 - \alpha)  \leq \frac{1}{m+1} + d_{\mathrm{KS}}(\widetilde F_0(R_{m+1}), U(0,1)). 
\]
This completes the proof of the corollary. 
\section{Proof of \Cref{rem:sharper_exchangeability}}
\label{app:proof_rem_exchangeavility}
We define the following for any measurable set $A$, 
\[
\mu(A) \coloneqq \mathbb{P}(S \in A), 
\qquad 
\mu_i(A) \coloneqq \mathbb{P}(S^i \in A), \quad i \in [m+1].
\]
By definition we have,
\[
\Gamma((R_1,\ldots,R_{m+1}))
=
\sup_{A}
\left|
\mu(A)-\frac{1}{m+1}\sum_{i=1}^{m+1}\mu_i(A)
\right|.
\]
Since $S^{m+1}=S$, we have $\mu_{m+1}=\mu$. Therefore for any measurable set $A$,
\begin{align*}
\mu(A)-\frac{1}{m+1}\sum_{i=1}^{m+1}\mu_i(A)
&=
\frac{1}{m+1}\sum_{i=1}^{m+1}\bigl(\mu(A)-\mu_i(A)\bigr) \\
&=
\frac{1}{m+1}\sum_{i=1}^{m}\bigl(\mu(A)-\mu_i(A)\bigr).
\end{align*}
Taking absolute values and applying the triangle inequality gives,
\begin{align*}
\left|
\mu(A)-\frac{1}{m+1}\sum_{i=1}^{m+1}\mu_i(A)
\right|
&\le
\frac{1}{m+1}\sum_{i=1}^{m}\left|\mu(A)-\mu_i(A)\right| \\
&=
\frac{1}{m+1}\sum_{i=1}^{m}
\left|
\mathbb{P}(S\in A)-\mathbb{P}(S^i\in A)
\right|.
\end{align*}
By the definition of total variation distance, we have for each $i\in[m]$,
\[
\left|
\mathbb{P}(S\in A)-\mathbb{P}(S^i\in A)
\right|
\le
d_{\mathrm{TV}}(S,S^i).
\]
Hence we have the following for any measurable set $A$,
\[
\left|
\mu(A)-\frac{1}{m+1}\sum_{i=1}^{m+1}\mu_i(A)
\right|
\le
\frac{1}{m+1}\sum_{i=1}^{m} d_{\mathrm{TV}}(S,S^i).
\]
Since this bound holds uniformly over all measurable sets $A$, taking the supremum over $A$ yields,
\[
\Gamma((R_1,\ldots,R_{m+1}))
\le
\frac{1}{m+1}\sum_{i=1}^{m} d_{\mathrm{TV}}(S,S^i).
\]
This completes the proof of the remark.
\section{Proof of \Cref{cor:randomized_test}}
\label{app:proof_rand_test_invariance}
For the application of \Cref{thm:main_result}, we note the following,
\[
Z:=D_m=\{X_1,\ldots,X_m\},
\qquad
\psi(Z):=T(X),
\qquad
W_b:=T(g_b^*X), \quad b\in[B].
\]
Since \(g_1^*,\ldots,g_B^*\) are sampled independently and uniformly (with replacement) from the finite group \(G\), it follows that conditional on \(Z=D_m\), the random variables \(W_1,\ldots,W_B\) are i.i.d. Indeed, for each \(b\in[B]\) and any Borel set \(A\subset \mathbb R\),
\[
\mathbb P_{H_0}(W_b\in A\mid Z)
=
\frac{1}{|G|}\sum_{g\in G}\mathbf 1\{T(gX)\in A\},
\]
which does not depend on \(b\). Let $W_{(1)}\le \ldots \le W_{(B)}$ denote the order statistics of \(W_1,\ldots,W_B\). We set $k:=\left\lceil (B+1)(1-\alpha)\right\rceil$. Then the modified testing rule can be written as,
\[
\phi^{\mathtt{mod}}(X)=\mathbf 1\{\psi(Z)\ge W_{(k)}\}.
\]
We now apply the one-sided version of the second statement of Theorem~\ref{thm:main_result} with $ a = 0, B- b = k$,
\[
- d_{\mathrm{KS}}(F_0(Z),U(0,1))
\le
\mathbb P_{H_0}\bigl(\psi(Z)<W_{(k)}\bigr)
-
\left(1-\frac{b+1}{B+1}\right)
\le
d_{\mathrm{KS}}(F_0(Z),U(0,1)),
\]
where
\[
F_0(Z)=\mathbb P_{H_0}(W_1\le \psi(Z)\mid Z)
      =\mathbb P_{H_0}(T(g_b^*X)\le T(X)\mid D_m).
\]
We emphasize that in this one-sided setting the slack can be written in terms of the ordinary Kolmogorov--Smirnov distance \(d_{\mathrm{KS}}(\cdot,\cdot)\), rather than \(\widetilde d_{\mathrm{KS}}(\cdot,\cdot)\). This is because when \(a=0\), the proof of Theorem~\ref{thm:main_result} involves only a single lower-tail probability and therefore only the class of half-lines \((-\infty,t]\), not general intervals \((s,t]\). We observe that,
\[
1-\frac{b+1}{B+1}
=
1-\frac{B-k+1}{B+1}
=
\frac{k}{B+1}.
\]
Therefore we have,
\[
- d_{\mathrm{KS}}(F_0(Z),U(0,1))
\le
\mathbb P_{H_0}\bigl(T(X)<W_{(k)}\bigr)-\frac{k}{B+1}
\le
d_{\mathrm{KS}}(F_0(Z),U(0,1)).
\]
Since $\phi^{\mathtt{mod}}(X)=\mathbf 1\{T(X)\ge W_{(k)}\}$, we have the following,
\[
\mathbb E_{H_0}[\phi^{\mathtt{mod}}(X)]
=
\mathbb P_{H_0}(T(X)\ge W_{(k)})
=
1-\mathbb P_{H_0}(T(X)<W_{(k)}).
\]
Therefore we have the following bounds on the type-I error of the test $\phi^{\mathtt{mod}}(X)$,
\[
\left|
\mathbb E_{H_0}[\phi^{\mathtt{mod}}(X)]
-
\left(1-\frac{k}{B+1}\right)
\right|
\le
d_{\mathrm{KS}}(F_0(D_m),U(0,1)).
\]
Since $k=\left\lceil (B+1)(1-\alpha)\right\rceil$, we have the identity,
\[
1-\frac{k}{B+1}
=
\frac{B+1-\lceil (B+1)(1-\alpha)\rceil}{B+1}
=
\frac{\lfloor (B+1)\alpha\rfloor}{B+1}.
\]
Thus we obtain the first statement of the corollary,
\[
\left|
\mathbb E_{H_0}[\phi^{\mathtt{mod}}(X)]
-
\frac{\lfloor (B+1)\alpha\rfloor}{B+1}
\right|
\le
d_{\mathrm{KS}}(F_0(D_m),U(0,1)).
\]

To obtain the more explicit type-I error bounds, we note that,
\[
\alpha-\frac{1}{B+1}
\le
\frac{\lfloor (B+1)\alpha\rfloor}{B+1}
\le
\alpha.
\]
Combining this with the previous display gives,
\[
\alpha-\frac{1}{B+1}-d_{\mathrm{KS}}(F_0(D_m),U(0,1))
\le
\mathbb E_{H_0}[\phi^{\mathtt{mod}}(X)]
\le
\alpha+d_{\mathrm{KS}}(F_0(D_m),U(0,1)).
\]
This completes the proof of the corollary.
\end{document}